%% file: neutral_shields.tex
\documentclass[11pt,letterpaper,twoside,reqno,nosumlimits]{amsart}
\usepackage{newclude}
\usepackage[margin=1.05in]{geometry}
\input{macro-3}
\usepackage{xurl}
\hypersetup{pdftitle={A sharp uniqueness criterion for neutral vortex arrays},
  pdfsubject={Planar Euler evolution and prescribed vorticity envelopes}}
\allowdisplaybreaks[1]
\title{A Sharp Uniqueness-Nonuniqueness Dichotomy\\ for Neutral Vortex Arrays in $2$D Euler Equations}
\author{Thomas Y. Hou}
\address{Applied and Computational Mathematics,
California Institute of Technology, Pasadena, CA, USA}
\email{hou@cms.caltech.edu}

\author{Xiang Qin}
\address{Applied and Computational Mathematics,
California Institute of Technology, Pasadena, CA, USA}
\email{xqin2@caltech.edu}

\author{Peicong Song}
\address{Applied and Computational Mathematics,
California Institute of Technology, Pasadena, CA, USA}
\email{psong2@caltech.edu}

\author{Wenhan Zhang}
\address{Applied and Computational Mathematics,
California Institute of Technology, Pasadena, CA, USA}
\email{wzhang9@caltech.edu}
\date{}
\begin{document}
\maketitle
\begin{abstract}
We establish a sharp uniqueness--nonuniqueness criterion for
the two-dimensional incompressible Euler equations within a class
of stationary arrays of disjoint neutral vortices.
Let $\Theta$ be any finite nondecreasing envelope for the
vorticity $L^p$ norms, and set
\(
    \mathcal I_\Theta=\int_{p_0}^\infty \Theta(p)^{-1}\idiff p.
\)
If $\mathcal I_\Theta=\infty$, every array in the corresponding class
admits a unique Euler evolution, even among bounded
energy-admissible weak velocity solutions.
If $\mathcal I_\Theta<\infty$, there exists a compactly supported array
in the same class admitting two distinct unforced Lagrangian
Euler evolutions from the same initial datum.
The criterion holds for arbitrary nondecreasing envelopes $\Theta$.
\end{abstract}
\tableofcontents

\section{Introduction}\label{sec:introduction}

\subsection{Background}

We consider the two-dimensional incompressible Euler equations in
vorticity form,
\begin{equation}\label{eq:euler}
 \partial_t\omega+\Div(u\omega)=0,\qquad
 u=K*\omega,
 \qquad (t,x)\in(0,T)\times\R^2,
\end{equation}
with the conventions
\begin{equation}\label{eq:conventions}
 x^\perp=(-x_2,x_1),\qquad
 \nabla^\perp=(-\partial_2,\partial_1),\qquad
 K(x)=\frac{x^\perp}{2\pi|x|^2}\quad(x\neq0).
\end{equation}
Here \(\omega=\curl u=\partial_1u_2-\partial_2u_1\) is the
scalar vorticity.  The term \emph{canonical velocity} will always mean
the convolution \(K*\omega\), with no additional harmonic vector field.
The corresponding momentum equations are
\begin{equation}\label{eq:momentum}
 \partial_tu+\Div(u\otimes u)+\nabla P=0,\qquad \Div u=0,
\end{equation}
where \(P\) is the pressure.  All the evolutions considered in this
paper are unforced.

The transport of vorticity is central to the global theory in two
dimensions.  For smooth compactly supported initial vorticity on the plane,
the canonical Euler evolution exists uniquely for all time.
The global classical theory goes back to Wolibner and Kato
\cite{Wolibner1933,Kato1967}; systematic accounts are given in
\cite{MajdaBertozzi2002,MarchioroPulvirenti1994}.
Yudovich's theorem \cite{Yudovich1963} extends global existence and
uniqueness to bounded integrable vorticity.  At this regularity the
velocity need not be Lipschitz, but its logarithmic modulus of
continuity still controls the particle flow and the nonlinear
evolution.  Yudovich subsequently allowed unbounded vorticity with
controlled growth of its finite Lebesgue norms \cite{Yudovich1995}.
Localized versions of this theory are developed by Crippa and
Stefani \cite{CrippaStefani2024}.

These results distinguish two questions about unbounded vorticity.
One can prescribe membership in \(L^p\) for a fixed finite exponent,
or require membership in every finite \(L^p\).  A more quantitative
question prescribes the rate at which \(\norm{\omega}_{L^p}\) may
increase as \(p\to\infty\).  Exact transport preserves all these
norms, but uniqueness depends on whether their growth permits an
effective estimate for the difference of two evolutions.  The
purpose of this paper is to identify an exact growth threshold for
a class of stationary arrays of neutral vortices.

Weak existence reaches substantially rougher classes.  The
compactness theory of DiPerna and Majda \cite{DiPernaMajda1987}
and Delort's existence theorem for distinguished-sign vortex-sheet
data \cite{Delort1991} illustrate the role of concentration and
sign in this problem.  Uniqueness requires additional information.
There is a further distinction at the Lagrangian level: the
Sobolev and bounded-variation transport theories of DiPerna--Lions
and Ambrosio \cite{DiPernaLions1989,Ambrosio2004} concern a
prescribed velocity field.  In Euler, the transported scalar also
determines that velocity.  The solutions constructed below each
have a unique regular Lagrangian flow for their own velocity,
although their common initial vorticity generates two different
coupled evolutions.

{The central question motivating this work is whether the growth restrictions on the vorticity \(L^p\) norms underlying the classical Yudovich uniqueness theory are sharp. More generally, we seek to understand how the growth of these norms as \(p\to\infty\) determines whether the Euler evolution is unique. We address this question for a class of compactly supported stationary Euler flows consisting of disjoint radial vortices with zero circulation. A remarkable feature of these vortices is that their velocity vanishes exactly outside their supports. This geometric cancellation yields a bound on the strain of the stationary velocity with a constant independent of \(p\), allowing us to improve the classical relative-energy estimate. We prove uniqueness whenever the integral of the reciprocal vorticity-growth envelope diverges.

We then show that this condition is sharp within the same class. When the integral converges, we construct an initial vorticity admitting two distinct unforced Euler evolutions, one stationary and the other nonstationary. The construction exploits an unstable second angular mode of a two-layer radial vortex. A small deformation of this vortex generates an exterior strain that excites the unstable mode of a neighboring larger vortex. By arranging such vortices at successively smaller spatial scales and shorter intrinsic time scales, we obtain a cascade of unstable transitions propagating from fine to coarse scales in finite time. At every finite cutoff, the perturbation is confined initially to the finest scale. As the cutoff tends to infinity, this initial discrepancy vanishes, while its dynamically amplified effect persists at a fixed larger scale. The resulting limiting solutions have the same initial vorticity, exactly transport the same distribution of vorticity, and conserve kinetic energy and all initially finite vorticity norms.

The significance of this result is that the same integral condition governs both the uniqueness theorem and the construction of nonunique solutions. Thus, the geometry of the stationary Euler flow determines a precise boundary between uniqueness and nonuniqueness within the specified initial-data class. Although this does not settle the sharpness of the classical Yudovich criterion for general unbounded vorticities, it provides a concrete setting in which the competing mechanisms of uniqueness and instability can be understood quantitatively.}

\subsection{Statement of the results}

Fix \(p_0>2\) and a finite, nondecreasing function
\(\Theta:[p_0,\infty)\to[1,\infty)\).  Its vorticity class and the
integral governing our result are
\begin{equation}\label{eq:envelope}
 Y^\Theta
 =\left\{f:\sup_{p\geq p_0}
             \frac{\norm{f}_{L^p(\R^2)}}{\Theta(p)}<\infty\right\},
 \qquad
 \mathcal I_\Theta=\int_{p_0}^\infty\frac{\idiff p}{\Theta(p)}.
\end{equation}
The envelope need not be continuous.  Its monotonicity and the integral
in \eqref{eq:envelope} are the only growth assumptions.

The basic vortex is the two-layer scalar
\begin{equation}\label{eq:unit-shield}
 a=\frac3\pi,\qquad q_{\mathrm c}=\frac1{\sqrt3},\qquad
 \Omega=\frac{3a}{2}\ind_{B_{q_{\mathrm c}}}
               -\frac a2\ind_{B_1}.
\end{equation}
Here \(B_R=B(0,R)\), \(B(c,R)\) is the open Euclidean disk,
and \(\ind_D\) is the indicator of \(D\).
The values of \(\Omega\) are \(a\), \(-a/2\), and zero.
Its positive and negative masses both have magnitude one.
In particular, its circulation vanishes, and its canonical velocity
vanishes outside \(B_1\).  We call a compactly supported radial vortex
with zero circulation a \emph{neutral shield}; the term refers to
its exact cancellation of the exterior velocity.
The nonuniqueness construction uses translates and rescalings of
the two-layer profile \(\Omega\).  We call one such two-layer vortex,
together with its subsequently transported domains, a \emph{cell}.

\begin{definition}[The stationary data class]\label{def:shields}
The class \(\mathscr S\) consists of all vorticities of the form
\begin{equation}\label{eq:general-array}
 \omega_0(x)=\sum_{m\geq0}\Omega_m(|x-c_m|),
\end{equation}
where \(c_m\in\R^2\), \(d_m>0\), and the measurable real-valued
profiles \(\Omega_m:[0,\infty)\to\R\) satisfy
\begin{equation}\label{eq:radial-neutrality}
 \Omega_m(\varrho)=0\quad\text{for almost every }\varrho\geq d_m,
 \qquad
 \int_0^{d_m}\varrho\Omega_m(\varrho)\idiff\varrho=0.
\end{equation}
The disks \(B(c_m,d_m)\) are pairwise disjoint and lie in one bounded
set, and the sum satisfies
\(\omega_0\in\bigcap_{1\leq p<\infty}L^p(\R^2)\).
A finite array is included by allowing all but finitely many
profiles to vanish.

\end{definition}

In particular, the fixed-profile arrays of the following form
\[
\omega_0(x)=\sum_{m\geq0}
      \alpha_m\Omega\left(\frac{x-c_m}{d_m}\right)
\]
are included by taking
\(\Omega_m(\varrho)=\alpha_m\Omega(\varrho/d_m,0)\), with
\(\alpha_m\in\R\). 
These are the arrays used in the
nonuniqueness construction.
Every datum in Definition~\ref{def:shields} has zero circulation
and a compactly supported canonical velocity.
Lemma~\ref{uniq:stationary-lemma} proves its stationarity,
including at accumulation points of the disks.

\begin{definition}[Lagrangian Euler evolution]\label{def:lagrangian}
A Lagrangian Euler evolution on \([0,T]\) is a triple
\((\omega,u,X)\) with
\begin{equation}\label{eq:solution-class}
 \omega\in C\left([0,T];L^p(\R^2)\right)
 \quad(1\leq p<\infty),\qquad
 u=K*\omega\in C\left([0,T];L^2(\R^2)\right)
                  \cap L^\infty\left((0,T)\times\R^2\right),
\end{equation}
satisfying \eqref{eq:euler} and \eqref{eq:momentum} in distributions.
The vorticities are supported in one compact set on each finite time
interval.  The map \(X(t,x_0)\) is an almost-everywhere invertible,
measure-preserving regular Lagrangian flow of \(u\): for almost every
\(x_0\), its time path is absolutely continuous and satisfies
\begin{equation}\label{eq:flow}
 X(0,x_0)=x_0,\qquad
 X(t,x_0)=x_0+\int_0^t u\left(s,X(s,x_0)\right)\idiff s.
\end{equation}
At every time it transports the initial vorticity exactly, in the sense
that
\begin{equation}\label{eq:transport}
 \omega\left(t,X(t,x_0)\right)=\omega_0(x_0)
 \quad\text{for almost every }x_0\in\R^2.
\end{equation}
\end{definition}

The null set in \eqref{eq:transport} may depend on the time.
All uniqueness assertions for flows are understood up to null sets.
The strong continuity in \eqref{eq:solution-class} includes the initial
time, so the common datum in the theorem below is a strong trace in
every finite Lebesgue space.

The main result gives the complete dichotomy for this data class.

\begin{theorem}[Uniqueness and bifurcation]\label{thm:main}
Let \(p_0>2\) and let \(\Theta\) satisfy the assumptions preceding
\eqref{eq:envelope}.
\begin{enumerate}[label=\textup{(\roman*)},leftmargin=2.2em]
\item\label{main:unique}
If \(\mathcal I_\Theta=\infty\), then for every
\(\omega_0\in\mathscr S\cap Y^\Theta\), the stationary velocity
\(U=K*\omega_0\) is the unique bounded energy-admissible solution
with initial velocity \(U\), in the sense of
Definition~\ref{uniq:admissible-definition}, on every finite interval.
In particular, \(\omega_0\) has exactly one Lagrangian Euler
evolution, with stationary vorticity \(\omega(t)=\omega_0\).
\item\label{main:nonunique}
If \(\mathcal I_\Theta<\infty\), there is a nonzero
\(\omega_0\in\mathscr S\cap Y^\Theta\) with two distinct Lagrangian
Euler evolutions on \([0,1]\).  With branch labels
\(j\in\{\mathrm{quiet},\mathrm{trig}\}\), they satisfy
\begin{equation}\label{eq:main-branches}
 \omega^{\mathrm{quiet}}(t)=\omega_0,\qquad
 \omega^{\mathrm{trig}}(0)=\omega_0,\qquad
 \norm{\omega^j(t)}_{L^p}=\norm{\omega_0}_{L^p}
                  \leq C_\Theta\Theta(p)\quad(p\geq p_0).
\end{equation}
There are a real test \(\phi_*\in C_c^\infty(\R^2)\), a time
\(t_*\in(0,1)\), and a number \(\gamma_*>0\) such that
\begin{equation}\label{eq:main-gap}
 \int_{\R^2}\phi_*(x)
       \left(\omega^{\mathrm{trig}}(t_*,x)
                         -\omega_0(x)\right)\idiff x
       \geq\gamma_*.
\end{equation}
Both branches conserve kinetic energy and every absolutely integrable
function of the vorticity.
\end{enumerate}
The constant \(C_\Theta\) and the datum in part~\ref{main:nonunique}
may depend on \(\Theta\).
\end{theorem}

Here a \emph{bifurcating datum} means a datum with the stationary and
distinct triggered evolutions in part~\ref{main:nonunique}.
Theorem~\ref{thm:main} is an exact criterion for uniqueness of the
\emph{whole data class} \(\mathscr S\cap Y^\Theta\):
the divergent integral gives uniqueness for every member, whereas
the convergent integral provides a member with nonunique evolution.

For either branch, exact transport also gives
\begin{equation}\label{eq:intro-value-conservation}
 \int_{\R^2}\beta(\omega^j(t,x))\idiff x
 =\int_{\R^2}\beta(\omega_0(x))\idiff x
\end{equation}
for every Borel function \(\beta:\R\to\R\) with \(\beta(0)=0\) and
\(\beta(\omega_0)\in L^1\).  Thus the two branches have the same
distribution of vorticity values; their distinction is a difference
of spatial evolution, detected by \eqref{eq:main-gap}.

In part~\ref{main:unique}, a bounded weak solution of
\eqref{eq:momentum} is called \emph{energy-admissible} if
\(u\in C_{\mathrm w}([0,T];L^2)\cap L^\infty_{t,x}\),
it has initial value \(U=K*\omega_0\), and
\(\norm{u(t)}_2\leq\norm{U}_2\) for almost every \(t\).
Here \(C_{\mathrm w}\) denotes weak continuity in \(L^2\).
Theorem~\ref{uniq:admissible-theorem} proves uniqueness in this
class without a vorticity assumption on the competing weak solution.

The power-law boundary is \(\Theta(p)=p\): this envelope gives
uniqueness, while \(\Theta(p)=p^{1+\varepsilon}\) admits a
bifurcating datum for every \(\varepsilon>0\).
At the next boundary, \(\Theta(p)=p\log(e+p)\) gives uniqueness,
while \(\Theta(p)=p(\log(e+p))^{1+\varepsilon}\) admits a
bifurcating datum for every \(\varepsilon>0\).
The same integral test treats further iterated logarithms.
These examples concern one datum and one pair of evolutions
for the chosen envelope, with simultaneous control of all
finite exponents.

\begin{remark}[Motivation and scope of the sharp dichotomy]
\label{rem:scope-sharpness}
One of the main motivations for this work is to understand the
sharpness of the classical Yudovich--Osgood uniqueness criterion
\eqref{eq:intro-general-osgood}
for more general classes of unbounded vorticity.
Theorem~\ref{thm:main} does not resolve this broader question.
Instead, it establishes a sharp uniqueness--nonuniqueness
dichotomy within the geometrically restricted class
$\mathscr S\cap Y^\Theta$.
The improved uniqueness threshold in this class relies on the
radial structure and exact shielding of the stationary vortices,
which are not available for general initial data.
Thus, our result does not rule out the sharpness of the
classical Yudovich--Osgood criterion in broader vorticity classes.
We view the present result as an intermediate step toward
understanding this general sharpness problem.
\end{remark}
\begin{remark}
The sharpness asserted here concerns the initial-data class $\mathscr S\cap Y^\Theta$, not arbitrary vorticities with the same norm growth.
\end{remark}

\subsection{Relation to earlier work}

\paragraph{$\quad\circ$\ \textbf{Growth conditions and the role of the reference field.}}
Yudovich's energy method \cite{Yudovich1995} relates vorticity
growth to an optimized differential inequality.  For a general
vorticity field, the estimate of the velocity gradient carries
a factor comparable to \(p\) at large exponents.  A familiar
sufficient uniqueness condition is consequently
\begin{equation}\label{eq:intro-general-osgood}
 \int_{p_0}^\infty\frac{\idiff p}{p\Theta(p)}=\infty.
\end{equation}
The stationary shield geometry permits the larger range described
by \(\mathcal I_\Theta=\infty\).  The distinction comes from an
estimate for the particular reference velocity, and the
nonuniqueness part proves sharpness within the corresponding
initial-data class.  The examples following
Theorem~\ref{thm:main} display the resulting borderline explicitly.

\paragraph{$\quad\circ$\ \textbf{Instability and forced nonuniqueness.}}
Vishik \cite{Vishik2018I,Vishik2018II} established nonuniqueness
for forced planar Euler in finite-exponent vorticity classes.
An unstable vortex is placed in similarity variables, so that a
perturbation approaching the background as logarithmic time tends
to \(-\infty\) produces a second solution with the same initial
trace.  The forcing makes that background an exact solution.
The monograph \cite{VishikAccount2024} gives a detailed treatment.
Castro, Faraco, Mengual, and Solera
\cite{CastroFaracoMengualSolera2025} simplify the construction of
the unstable vortex by first analyzing a piecewise constant
two-layer profile and then regularizing it.  Our unit shield is
a specific member of their piecewise constant family; the
nonlinear contour response and the interaction of infinitely
many such cells are the issues addressed here.
Albritton, Colombo, and Mescolini
\cite{AlbrittonColomboMescolini2025} obtain nonunique vanishing
viscosity limits near the forced Euler scenario, allowing initial
perturbations that depend on viscosity and vanish with it.

The instability method also connects Euler to nonuniqueness for
Navier--Stokes.  Albritton, Bru\`e, and Colombo
\cite{AlbrittonBrueColombo2022} construct distinct Leray solutions
of the forced three-dimensional equation.  For the unforced
equation, Jia and \v Sver\'ak
\cite{JiaSverak2014,JiaSverak2015} develop forward self-similar
solutions and a conditional route from their spectral instability
to nonuniqueness.  Guillod and \v Sver\'ak
\cite{GuillodSverak2023} provide numerical bifurcation evidence.
Hou, Wang, and Yang \cite{HouWangYang2026} present a
computer-assisted nonuniqueness proof for unforced Leray--Hopf
solutions, whereas Ionescu, Jia, and Palasek
\cite{IonescuJiaPalasek2026} construct numerical unstable profiles
in the axisymmetric swirl-free class.  The relevance of these
works is the passage from an unstable direction to distinct
nonlinear evolutions: the initial trace and the interaction with
the remaining modes are part of that passage.

\paragraph{$\quad\circ$\ \textbf{Euler similarity profiles.}}
There is a substantial theory of self-similar Euler flows.
Elling constructs algebraic spiral solutions and profiles with
mixed-sign vorticity \cite{Elling2013,Elling2016}; Shao, Wei,
and Zhang establish further algebraic spiral solutions
\cite{ShaoWeiZhang2025}.  Recent constructions include Choi's
large-data profiles \cite{Choi2026SelfSimilar}, the homogeneous
steady states of Choi and Coiculescu \cite{ChoiCoiculescu2026}
whose similarity dynamics has several attracting stagnation
points, and the Lagrangian construction of Hou, Qin, Song, and Zhang
\cite{HouQinSongZhang2026}.  The last work treats homogeneous
angular data without a smallness or sign assumption in its
stated class and retains exact backward particle transport.
Such profiles supply possible backgrounds for a similarity
approach.  In the present construction, the accumulation is in
spatially separated cells with increasingly short intrinsic
time scales.  The complete initial datum is compactly supported,
and its stationary branch is already determined by shielding.

\paragraph{$\quad\circ$\ \textbf{Weak nonuniqueness and positive-time vorticity.}}
The constructions of Scheffer and Shnirelman
\cite{Scheffer1993,Shnirelman1997}, and the differential-inclusion
formulation of De Lellis and Sz\'ekelyhidi
\cite{DeLellisSzekelyhidi2009}, establish nonuniqueness at low
velocity regularity.  Preserving positive-time vorticity
integrability is a further constraint.  Bru\`e and Colombo
\cite{BrueColombo2023} construct nonunique solutions with
vorticity in a Lorentz space.  Bru\`e, Colombo, and Kumar
\cite{BrueColomboKumar2026} obtain weak Euler solutions on the
torus with vorticity in \(C_tL^p_x\) for some \(p>1\), using
dipole-based convex integration.  Their low-exponent vorticity equation is interpreted through
the symmetrized Biot--Savart bilinear form.  The exact
Lagrangian transport required here is an additional property.

Mengual \cite{Mengual2024} constructs nonunique admissible
velocities from truncated power-law initial vortices with
\(L^1\cap L^p\) integrability for any prescribed finite \(p>1\).  In the evolving mixing region, however, the
positive-time vorticity is controlled as a distribution.
Theorem~\ref{thm:main} keeps the prescribed growth of all
finite vorticity norms throughout the evolution.  Its two
solutions are rearrangements of the same scalar, conserve
kinetic energy, and are distinguished by their spatial
organization through the fixed test in \eqref{eq:main-gap}.

\paragraph{$\quad\circ$\ \textbf{Recent unforced constructions.}}
Choi \cite{Choi2026} constructs a compactly supported,
zero-circulation common datum in every finite \(L^p\), with
two unforced evolutions that remain strongly continuous in
each of those spaces.  The mechanism is a cascade of
transitions from rotating quadrupoles to separating dipoles.
The quantitative norm calculation in that work gives
\(\norm{\omega_0}_p\asymp p^{4/\alpha}\) for each fixed
\(0<\alpha<1\), with comparison constants independent of
\(p\).  Theorem~\ref{thm:main} addresses a prescribed
envelope: its scale selection applies whenever
\(\mathcal I_\Theta<\infty\), including envelopes arbitrarily
close to the divergent-integral boundary.  The stationary
shield class also supplies the matching uniqueness theorem.

Coiculescu \cite{Coiculescu2026} studies departure from a
radial nonnegative datum through a threefold-symmetric
multiscale construction.  Part~I establishes a compactness
criterion and the linear analysis; a companion draft and
interval-arithmetic certificate materials are available in
the repository linked from that paper.  This approach uses
a radial reference and spectral certification of its building
blocks.  The Hexagon manuscript \cite{Hexagon2026} proposes,
for each finite \(p>2\), two strong vanishing-viscosity
accumulation paths from one compactly supported datum, using
concentric unstable profiles and their interacting histories.
Its repository identifies it as an AI-generated manuscript
whose proof is being studied by a reading group.
These works formulate different forms of unforced
nonuniqueness.  The assertion proved here is the integral
dichotomy for the disjoint shield arrays of
Definition~\ref{def:shields}.

\paragraph{$\quad\circ$\ \textbf{Relation to our companion work}} {In our companion paper \cite{HQPZ26}, we establish nonuniqueness for the unforced two-dimensional Euler equations using unstable four-vortex configurations that evolve into separating dipoles. The present construction uses a different local mechanism, based on the unstable second angular mode of a stationary two-layer radial vortex. Although both constructions propagate instability from increasingly fine to larger scales, their building blocks and reference evolutions are fundamentally different. More importantly, the exact shielding property in the present setting allows us to establish a matching uniqueness theorem and thereby identify a sharp uniqueness–nonuniqueness threshold in terms of the growth of the vorticity ($L^p$) norms}.

\paragraph{$\quad\circ$\ \textbf{Cascades and the realization of vortex dynamics.}}
Successive amplification across scales appears in the dyadic
Euler model of Katz and Pavlovi\'c \cite{KatzPavlovic2005}
and in Tao's averaged Navier--Stokes equation
\cite{Tao2016}.  The manuscripts
\cite{OpenAIEuler2026,OpenAINavierStokes2026} develop
amplification mechanisms for three-dimensional singularity
formation: the Euler construction transfers growth between
localized oscillations, while the Navier--Stokes construction
uses a collapsing background, oscillatory corrections, and
an external force.  Across these settings, a local growing
mode becomes useful only when its output supplies a controlled
input at the next scale and the accumulated errors remain
summable.  Here the transfer travels from smaller, faster
cells to larger, slower cells and produces a fixed-scale
difference between two limits.

The realization of finite-dimensional vortex dynamics by
transported vorticity is developed in the concentration
theory of Turkington, Marchioro--Pulvirenti, and
Butt\`a--Marchioro
\cite{Turkington1987,MarchioroPulvirenti1993,ButtaMarchioro2018}.
Rearrangement inequalities such as
\cite{BrascampLiebLuttinger1974} and the singular-integral
estimates in \cite{Stein1970} provide tools for controlling
interaction energies and velocity fields.  The global patch
regularity theorems of Chemin and Bertozzi--Constantin
\cite{Chemin1993,BertozziConstantin1993} give another
important part of this background.  Our local building block
is an exact finite-size Euler vortex.  We work directly with
its two interfaces while its shape is small, and use
bounded-vorticity evolution and moment bounds after it leaves
that neighborhood.  Quantitative flow stability, as developed
by Crippa and De Lellis \cite{CrippaDeLellis2008}, is a
related route to Lagrangian compactness.  In Section~\ref{sec:limit}
we first identify the limiting velocity and then use
renormalization to recover exact transport and strong
vorticity convergence.

\subsection{Strategy of the proof}

{We now outline the main analytical ingredients behind the two directions of our theorem. The uniqueness argument exploits an exponent-uniform strain estimate for stationary neutral-vortex arrays. The nonuniqueness construction combines an exact nonlinear growing orbit of a two-layer vortex with quantitative control of its response to the exterior fields generated by other vortices. The main challenge is to arrange the interactions across infinitely many scales while keeping the approximation errors and total transition time under control.}

\paragraph{$\quad\circ$\ \textbf{A relative-energy estimate without an exponent loss.}}
For the stationary velocity \(U=K*\omega_0\),
Lemma~\ref{uniq:stationary-lemma} proves
\begin{equation}\label{eq:intro-strain}
 \norm{S(U)}_{L^p}
 \leq\frac{2p-1}{2(p-1)}\norm{\omega_0}_{L^p}
 \leq\frac32\norm{\omega_0}_{L^p},
 \qquad 2\leq p<\infty,
\end{equation}
where \(S(U)=(\nabla U+\nabla U^{\mathsf T})/2\) is measured in
the operator norm.
The estimate follows from a radial Hardy inequality and passes
to the infinite array because the velocity supports are disjoint.  If \(u\) is another canonical
evolution from the same initial datum and
\(\mathcal E(t)=\norm{u(t)-U}_2^2/2\), the relative-energy
identity and interpolation give an inequality of the form
\begin{equation}\label{eq:intro-energy}
 \mathcal E'(t)\leq C\Theta(p)\mathcal E(t)^{1-1/p}
 \qquad(p\geq p_0),
\end{equation}
with \(C\) independent of \(p\).  Optimizing over \(p\) leads
to the concave modulus
\begin{equation}\label{eq:intro-modulus}
 \mu_\Theta(y)=\inf_{p\geq p_0}
                   \Theta(p)y^{1-1/p},\qquad 0<y\leq1.
\end{equation}
Lemma~\ref{uniq:modulus-equivalence} proves that the reciprocal
integral of this modulus diverges exactly when
\(\mathcal I_\Theta=\infty\), even for discontinuous
envelopes.  Osgood's argument then forces
\(\mathcal E=0\).  For energy-admissible weak momentum
solutions, Theorem~\ref{uniq:admissible-theorem} applies the
same argument to an absolutely continuous majorant of the
relative energy.

\paragraph{$\quad\circ$\ \textbf{An exact nonlinear orbit of one shield.}}
We describe a cell by its two boundary curves, keeping their
enclosed areas fixed and separating the common and relative
centroids from the genuine shape variables.  The growing
second angular mode has eigenvalue
\begin{equation}\label{eq:intro-rate}
 \lambda=\sigma-\frac{\iunit a}{4},\qquad
 \sigma=\frac{a}{4\sqrt3}.
\end{equation}
Its positive real part gives amplification, while its
imaginary part records the rotation of the perturbation.
Proposition~\ref{cell:donor} constructs an analytic contour
family parametrized by a complex growing coordinate \(z\);
its associated patch vorticity is denoted by \(G(z)\).  For a sufficiently small fixed section amplitude
\(a_*>0\), the orbit
\begin{equation}\label{eq:intro-orbit}
 W_*(s)=G(a_*e^{\lambda s}),\qquad s\leq0,
\end{equation}
solves Euler exactly and approaches \(\Omega\) as
\(s\to-\infty\).  Here \(s\) is unit-cell time.  The
analyticity belongs to the contours; the corresponding
vorticities are compared through the areas of symmetric
differences.  The bounded-vorticity theorem continues
\(W_*\) for positive time.

This exact orbit is essential on a long waiting interval.
Near it, the local equation splits into a growing complex
coordinate and a complementary shape component.  Flattening
the nonlinear family removes the autonomous error from the
growing equation.  Proposition~\ref{cell:response} controls
the response to an exterior harmonic velocity while tracking
the analytic radius, the centroids, and the length of the
waiting interval.  Proposition~\ref{cell:first-exit} then
locates the first crossing of the section.

\paragraph{$\quad\circ$\ \textbf{The complete outgoing moment determines the next input.}}
Identify the plane with \(\C\) for this calculation, and
define the barred quadratic moment of the isolated orbit by
\(Q_*(s)=\int_{\R^2}\overline y^{\,2}W_*(s,y)\idiff y\).
The deformation creates an exterior strain whose projection
onto a larger cell's growing mode depends on this complex
moment.  For a time-scale ratio \(r\in(1/2,1)\), the relevant
coefficient is the complete transform
\begin{equation}\label{eq:intro-transform}
 L=\int_{\R}e^{-r\lambda s}Q_*(s)\idiff s.
\end{equation}
Corollary~\ref{cell:transform} proves absolute convergence
and a nonzero leading term as \(r\) approaches one after
\(a_*\) is fixed.  Using the complete orbit includes both
the approach to the section and its subsequent Euler
evolution.  The complex weight retains the phase that a
scalar size estimate would discard.

Cells are placed in pairs related by \(x\mapsto-x\).  This
symmetry cancels the total impulse, meaning the signed first
moment of vorticity, of each pair.  The
paired multipole expansion controls the smaller impulse
generated in an individual cell.  Alternating reflection,
together with the corresponding sign change of vorticity,
makes the outgoing complex phase enter the next growing
coordinate without conjugation.  These identities yield
the transfer coefficient in \eqref{cas:clocks}.

\begin{figure}[t]
    \centering
    \includegraphics[width=0.95\textwidth]
        {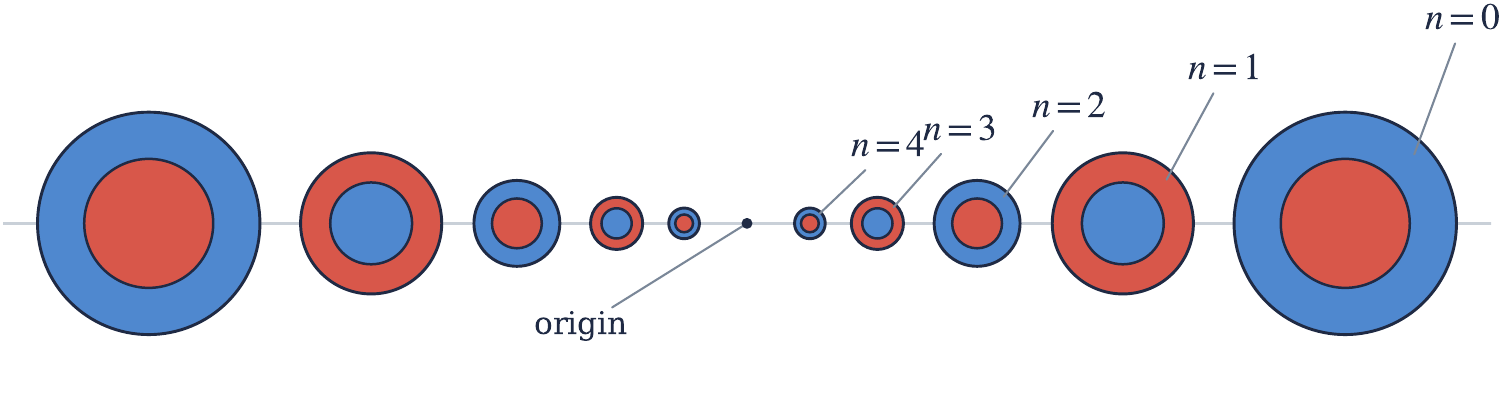}
    \caption{
        Schematic arrangement of the multiscale neutral shields.
        At each level \(n\), a pair of radial shields of radius
        \(d_n\) is centered at \(\pm R_ne_1\).
        The pairs accumulate at the origin, with the signs of
        the inner and outer layers alternating between successive
        levels.
        Only the first five levels are shown; the relative sizes
        and separations are not drawn to scale.
    }
    \label{fig:intro-neutral-shields}
\end{figure}

\paragraph{$\quad\circ$\ \textbf{One actual Euler solution at each finite cutoff.}}
At level \(n\), the centers are \(\pm R_ne_1\), where
\(e_1=(1,0)\), the physical radius is \(d_n=h_nR_n\),
and the vorticity amplitude is \(A_n=r^{-n}\).
Its intrinsic time scale is therefore \(\tau_n=A_n^{-1}=r^n\).
For a logarithmic scale sequence \((S_n)\), we take
\begin{equation}\label{eq:intro-scales}
 R_n=\exp\left(-(S_n-S_0)\right),\qquad h_n=S_n^{-M},
\end{equation}
where \(M>0\) is chosen in Section~\ref{sec:envelopes}.
The geometric arrangement is illustrated in
Figure~\ref{fig:intro-neutral-shields}.
As \(n\) increases, the paired shields approach the origin
and evolve on progressively shorter time scales.
The leading interaction transfers perturbations from
smaller, faster cells to larger, slower ones through their
induced exterior velocity fields, without material exchange
between the cells.

\begin{figure}[t]
    \centering
    \includegraphics[width=\textwidth]
        {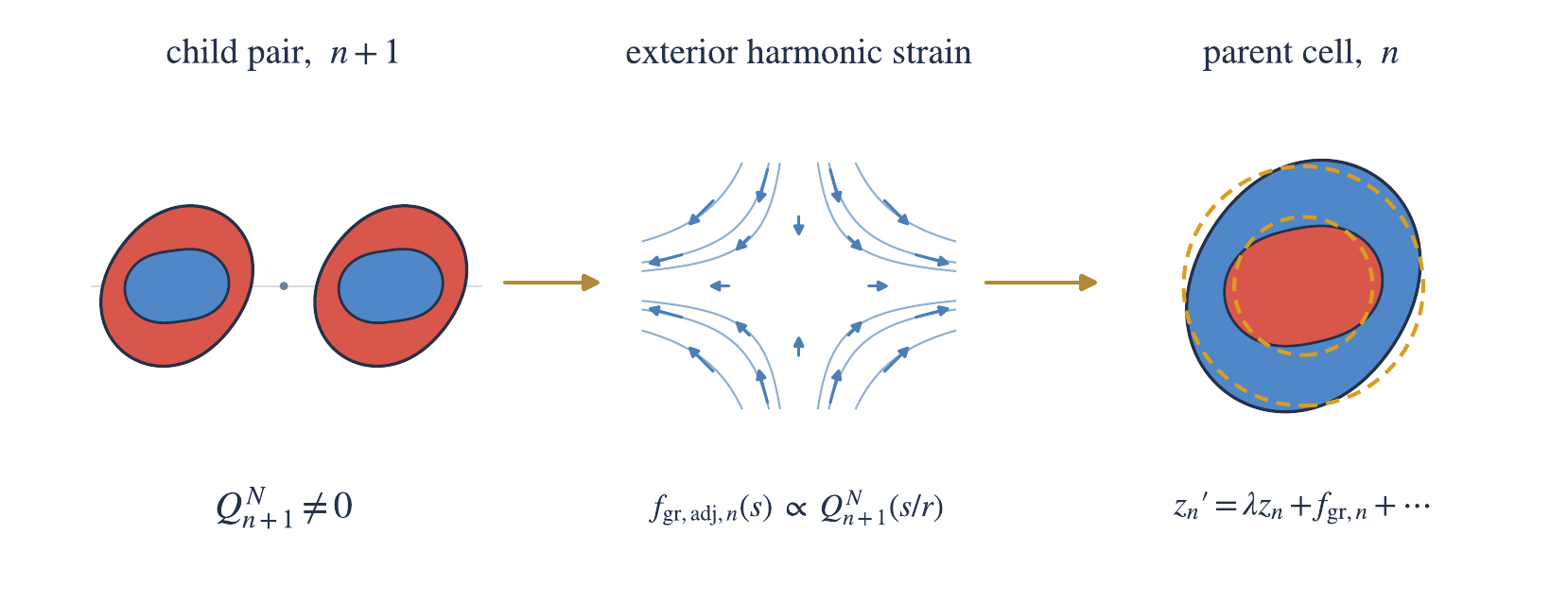}
    \caption{
        Schematic of the transfer mechanism between adjacent levels.
        A deformation of the child pair at level \(n+1\) generates
        an exterior field that is harmonic near the parent cell at
        level \(n\).
        Its strain component forces the growing second angular mode
        of the parent cell, with the effective input determined by
        the weighted quadrupole history of the child pair.
        Dashed circles indicate the unperturbed radial interfaces.
        The geometry and deformations are schematic.
    }
    \label{fig:intro-transfer-mechanism}
\end{figure}

The geometric arrangement alone does not explain the cascade.
What matters dynamically is how a perturbation created at a
smaller scale is transmitted to the next larger one.
Figure~\ref{fig:intro-transfer-mechanism} illustrates this
adjacent-level mechanism schematically.
Once the child pair at level \(n+1\) develops a nonzero
quadrupole, its exterior velocity field is harmonic near the
parent cell at level \(n\).
The corresponding strain forces the growing second angular mode
of the parent cell and initiates its departure from the radial
shield.
The effective input is determined by the weighted time history
of the child's quadrupole, rather than by a single instantaneous
value.
Throughout this process, the material cells remain disjoint;
the transfer is mediated by the induced exterior field, not by
material exchange between scales.

For a finite cutoff \(N\), every cell except the innermost
pair is initially circular.  That pair receives a small
perturbation from \(G\).  We then evolve the complete finite
datum by its bounded-vorticity Euler solution.  In particular,
every exterior field acting on a cell is generated by the
other material components of the same evolution.

The main difficulty is to keep the larger cells close to
their circular states while they wait for the smaller cells
to act.  A bound on their total possible strain, multiplied
by the full waiting time, is insufficient.  Section~\ref{sec:finite-array}
instead estimates their response through a causal integral
inequality, using only forcing already produced at each time.
The spectral gap between neighboring amplitudes yields a
small feedback coefficient.  An integrated input is then
recorded before the receiver can reach its section; this
determines its first departure and the phase passed to the
next level.  Uniform estimates of the complete analytic
exterior field, including its impulse and higher-moment
remainders, close the induction in Theorem~\ref{cas:finite}.

\paragraph{$\quad\circ$\ \textbf{The reciprocal-envelope integral.}}
For an adjacent pair, the transfer coefficient \(D_n\) in
\eqref{cas:clocks} has size comparable to
\(h_{n+1}^4(R_{n+1}/R_n)^4\), with constants fixed by the
local orbit.  The nominal delay \(\ell_n\) between successive
section crossings consequently satisfies
\begin{equation}\label{eq:intro-delay}
 \sigma\ell_n
 =4(S_{n+1}-S_n)+4M\log S_{n+1}+O(1).
\end{equation}
The \(O(1)\) term is uniform in \(n\).  The corresponding physical
delay is \(\tau_n\ell_n\), up to the exit-time errors controlled
in Theorem~\ref{cas:finite}.
At the same time, disjointness gives the exact norm sum
\begin{equation}\label{eq:intro-norm-sum}
 \norm{\omega_0}_{L^p}^p
 =2\norm{\Omega}_{L^p}^p\sum_{n\geq0}A_n^pd_n^2.
\end{equation}
Thus larger amplitudes accelerate the evolution, whereas
smaller support areas preserve the prescribed norm growth;
the logarithm of the spatial separation measures the delay
needed to amplify the remote input.

For an arbitrary envelope with \(\mathcal I_\Theta<\infty\),
Lemma~\ref{env:minorant} constructs a smooth positive function
\(g\leq\Theta\), with controlled positive logarithmic slope
and \(\int^\infty\idiff p/g(p)<\infty\).
We choose \(s_0\) in its domain and define \(s_n\) and \(S_n\) by
\begin{equation}\label{eq:intro-envelope-scales}
 g(s_n)=g(s_0)r^{-n},\qquad S_n=\Lambda s_n,
\end{equation}
where \(\Lambda\) is the final large separation parameter.
Monotonicity and the fixed ratio
\(g(s_{n+1})/g(s_n)=r^{-1}\) give the comparison
\begin{equation}\label{eq:intro-clock-integral}
 \sum_{n\geq0}r^n(s_{n+1}-s_n)
 \asymp g(s_0)\int_{s_0}^\infty\frac{\idiff s}{g(s)}.
\end{equation}
The controlled logarithmic slope also makes the logarithmic
terms in \eqref{eq:intro-delay} summable.  Proposition~\ref{env:norms}
proves \(\norm{\omega_0}_p\asymp g(p)\) for large \(p\)
and the required upper bound by \(\Theta(p)\).
This is how the same scale choice reconciles a finite total
transition time with the prescribed growth of the vorticity
norms.

\paragraph{$\quad\circ$\ \textbf{A common datum and a fixed observable in the limit.}}
The initial perturbation is inserted ever deeper as
\(N\to\infty\), and the finite initial data converge to
the complete quiet array in every finite \(L^p\).
Nevertheless, Theorem~\ref{cas:finite} transports its effect
to the outermost cell in a uniformly bounded positive time.
One test from a fixed finite family detects a nonzero
quadratic moment there.  A subsequence fixes that test,
producing the observable \(\phi_*\) in \eqref{eq:main-gap}.

Section~\ref{sec:limit} first obtains strong velocity
compactness in \(C_tL^2_x\), using uniform \(H^1\) bounds,
an \(L^2\) time derivative, and zero-circulation decay at
infinity.  After passing the vorticity equation to the
limit, DiPerna--Lions theory identifies its solution with
exact transport by the limiting velocity.  Conservation
of the limiting norms then upgrades weak vorticity
convergence to strong convergence in every finite
\(C_tL^p_x\).  The fixed outer observable survives this
limit, while shielding supplies the stationary evolution
from the same datum.  This completes both the dynamical
separation and the Lagrangian assertions of
Theorem~\ref{thm:main}.

\subsection{Notation and organization}

Physical time is \(t\); \(s\) denotes unit-cell time in
the local dynamics and a scalar integration variable when
its meaning is specified by the formula.  The indices
\(n\) and \(N\) denote a scale level and a finite cutoff,
while angular Fourier modes are indexed by \(k\).
The labels \(\mathrm{quiet}\) and \(\mathrm{trig}\) refer
to the stationary and triggered branches in
\eqref{eq:main-branches}.  We identify \(\R^2\) with
\(\C\) only in contour and moment calculations; bars then
denote complex conjugation.  Spatial or Fr\'echet derivatives
are written as \(\mathrm D\).  Unless another domain is
specified, spatial norms are over \(\R^2\).  The notation
\(f\asymp g\) means two-sided comparison by positive
constants independent of the variable under discussion.
Generic constants may increase between estimates; their
allowed dependence is stated where it matters, and the
finite-cutoff estimates are uniform in \(N\).

The remainder of the paper is organized as follows.
Section~\ref{sec:uniqueness} proves the uniqueness direction
and its energy-admissible extension.
Section~\ref{sec:local-cell} develops the exact local contour
evolution, the nonlinear growing orbit, and its response to
an exterior field.  Section~\ref{sec:finite-array} proves
the finite-array transfer theorem with uniform control of
the waiting intervals.  Section~\ref{sec:envelopes} selects
the scales from the prescribed envelope.
Section~\ref{sec:limit} constructs the limiting Lagrangian
evolution and proves Theorem~\ref{thm:main}.

\section{Uniqueness for stationary radial arrays}
\label{uniq:section}\label{sec:uniqueness}

The uniqueness argument, following Yudovich's energy method
\cite{Yudovich1995}, uses a feature of the initial configuration that is
absent from a general vorticity field.  A radial Hardy inequality bounds
the strain of each neutral vortex by its vorticity, with a constant
uniform at large exponents.  Exact exterior shielding and disjointness
then give the same estimate for the entire array.
We first establish this estimate and then derive the
relative-energy identity in the class of solutions used here.  The resulting
Osgood condition contains $\Theta(p)$ without an additional factor of $p$.
Lemma~\ref{uniq:modulus-equivalence} also identifies the integral criterion
for the optimized energy modulus without any regularity assumption on the
growth envelope beyond monotonicity.

\subsection{The stationary field}
\label{uniq:stationary-subsection}

Let \(\omega_0\in\mathscr S\), with the representation
\eqref{eq:general-array}.  The following calculation supplies the
stationary solution and the estimate needed for uniqueness.

We use the canonical velocity $U=K*\omega_0$, with
$K(x)=x^\perp/(2\pi|x|^2)$, and the strain matrix
$S(U)=(\nabla U+\nabla U^{\mathsf T})/2$.  The norm of a matrix in a strain
estimate is its operator norm; the Frobenius norm will be indicated explicitly.

\begin{lemma}[Exact shielding and strain estimate]\label{uniq:stationary-lemma}
Let $\omega_0\in\mathscr S$ have the representation
\eqref{eq:general-array}.  Its canonical velocity $U=K*\omega_0$ is a
bounded, compactly supported stationary Euler velocity and belongs to
$W^{1,p}(\R^2)$ for every $1<p<\infty$.  It vanishes outside the union
of the disks $B(c_m,d_m)$, and its strain satisfies
\begin{equation}\label{uniq:strain-Lp}
 \begin{aligned}
 \norm{S(U)}_{L^p}
   &\leq\frac{2p-1}{2(p-1)}\norm{\omega_0}_{L^p}
       &&(1<p<\infty),\\
 \norm{S(U)}_{L^p}
   &\leq\frac32\norm{\omega_0}_{L^p}
       &&(2\leq p<\infty).
 \end{aligned}
\end{equation}
\end{lemma}

\begin{proof}
Write $\omega_m(x)=\Omega_m(|x-c_m|)$ and $U_m=K*\omega_m$, with each
$\Omega_m$ extended by zero for $\varrho\geq d_m$.  For a radial profile
$f$, define its radial average by
\begin{equation}\label{uniq:radial-average}
 (\mathsf H f)(\varrho)
   =\frac{2}{\varrho^2}\int_0^\varrho s f(s)\idiff s
   =2\int_0^1\tau f(\tau\varrho)\idiff\tau,
 \qquad \varrho>0.
\end{equation}
The radial Biot--Savart formula gives
\begin{equation}\label{uniq:unit-velocity}
 U_m(x)=v_{\theta,m}(\varrho)e_\theta,
 \qquad
 v_{\theta,m}(\varrho)
   =\frac1\varrho\int_0^\varrho s\Omega_m(s)\idiff s
   =\frac{\varrho}{2}(\mathsf H\Omega_m)(\varrho),
 \qquad
 \varrho=|x-c_m|>0,
\end{equation}
where $e_\theta=(x-c_m)^\perp/\varrho$.  Zero circulation makes both
$v_{\theta,m}$ and $\mathsf H\Omega_m$ vanish for $\varrho\geq d_m$.
For any fixed $q>2$, H\"older's inequality gives
$|v_{\theta,m}(\varrho)|\leq C_q\norm{\omega_m}_{L^q}
\varrho^{1-2/q}$ near zero.  Thus $U_m$ extends continuously by zero
at its center and across its outer circle.

Differentiating \eqref{uniq:unit-velocity} away from the center gives
$v_{\theta,m}'=\Omega_m-(\mathsf H\Omega_m)/2$ almost everywhere.
In the orthonormal polar frame, the gradient and strain therefore satisfy
\begin{equation}\label{uniq:pointwise-gradient}
 \begin{aligned}
 |\nabla U_m|_{\mathrm F}^2
   &=\left|\Omega_m-\frac12\mathsf H\Omega_m\right|^2
       +\frac14|\mathsf H\Omega_m|^2,\\
 |S(U_m)|&=\frac12|\Omega_m-\mathsf H\Omega_m|.
 \end{aligned}
\end{equation}
All radial quantities in this formula are evaluated at $\varrho$.
Minkowski's inequality applied to \eqref{uniq:radial-average} yields
the Hardy bound
\begin{equation}\label{uniq:hardy-bound}
 \begin{aligned}
 \norm{(\mathsf H\Omega_m)(|\cdot|)}_{L^p}
 &\leq2\int_0^1\tau
        \norm{\Omega_m(\tau|\cdot|)}_{L^p}\idiff\tau\\
 &=2\norm{\omega_m}_{L^p}
         \int_0^1\tau^{1-2/p}\idiff\tau
   =\frac p{p-1}\norm{\omega_m}_{L^p},
 \qquad 1<p<\infty.
 \end{aligned}
\end{equation}
The gradient given by \eqref{uniq:pointwise-gradient} consequently belongs
to $L^p$.
Integration by parts on punctured disks identifies it as the weak
gradient: the outer boundary term vanishes because
$v_{\theta,m}(d_m)=0$, and the inner boundary term is bounded by
$C\varrho|v_{\theta,m}(\varrho)|=O(\varrho^{2-2/q})$.
There is therefore no singular derivative at either boundary.
Equations~\eqref{uniq:pointwise-gradient} and \eqref{uniq:hardy-bound}
give
\begin{equation}\label{uniq:single-cell-estimates}
 \norm{\nabla U_m}_{L^p}\leq C_p\norm{\omega_m}_{L^p},
 \qquad
 \norm{S(U_m)}_{L^p}
   \leq\frac{2p-1}{2(p-1)}\norm{\omega_m}_{L^p},
 \qquad 1<p<\infty,
\end{equation}
where one may take $C_p=1+\frac{p}{\sqrt2(p-1)}$.

Disjointness gives the exact norm identity
\begin{equation}\label{uniq:array-norm}
 \norm{\omega_0}_{L^p}^p
   =\sum_{m\geq0}\norm{\omega_m}_{L^p}^p
   =2\pi\sum_{m\geq0}\int_0^{d_m}
                  \varrho|\Omega_m(\varrho)|^p\idiff\varrho,
 \qquad 1\leq p<\infty.
\end{equation}
Let $\omega_{0,N}=\sum_{m=0}^N\omega_m$ and
$U_N=K*\omega_{0,N}$.  Equation~\eqref{uniq:array-norm} implies
$\omega_{0,N}\to\omega_0$ in every finite $L^p$.
The near--far estimate
\begin{equation}\label{uniq:velocity-bound}
 \norm{K*f}_{L^\infty}
   \leq C_q\left(\norm{f}_{L^1}+\norm{f}_{L^q}\right),
 \qquad q>2,
\end{equation}
then gives uniform convergence $U_N\to U$.  Each $U_N$ vanishes
outside the union of the disks, so $U$ has the same property and
belongs to $L^2$.  Since all these velocities have support in one
fixed bounded set, their convergence also holds in every finite $L^p$.
The weak gradients of different vortices have disjoint supports; hence
\eqref{uniq:single-cell-estimates} gives, for $N>M$,
\begin{equation}\label{uniq:sobolev-tail}
 \norm{\nabla(U_N-U_M)}_{L^p}
   \leq C_p\left(\sum_{m=M+1}^N
                      \norm{\omega_m}_{L^p}^p\right)^{1/p}
   =C_p\norm{\omega_{0,N}-\omega_{0,M}}_{L^p}.
\end{equation}
Thus $U_N$ converges to $U$ strongly in $W^{1,p}$ for every $1<p<\infty$.
Summing the strain estimate in \eqref{uniq:single-cell-estimates}
over the disjoint disks and taking this limit proves
\eqref{uniq:strain-Lp}, because $(2p-1)/(2(p-1))\leq3/2$ for $p\geq2$.

For each $\varphi\in C_c^\infty(\R^2)$, integration in the angular
variable gives
\begin{equation}\label{uniq:radial-stationarity}
 \begin{aligned}
 \int_{\R^2}U_m\omega_m\cdot\nabla\varphi\idiff x
 &=\int_0^{d_m}v_{\theta,m}(\varrho)\Omega_m(\varrho)
       \left(\int_0^{2\pi}
          \partial_\theta\Bigl(\varphi\bigl(c_m+\varrho(\cos\theta,\sin\theta)\bigr)\Bigr)
          \idiff\theta\right)\idiff\varrho\\
 &=0.
 \end{aligned}
\end{equation}
Since the velocity of each vortex vanishes outside its disk,
\eqref{uniq:radial-stationarity} also gives
$\operatorname{div}(U_N\omega_{0,N})=0$ for every $N$.
The product difference satisfies
\begin{equation}\label{uniq:stationary-product-limit}
 \begin{aligned}
 \norm{U_N\omega_{0,N}-U\omega_0}_{L^1}
 &\leq\norm{U_N-U}_{L^\infty}\norm{\omega_{0,N}}_{L^1}
       +\norm{U}_{L^\infty}\norm{\omega_{0,N}-\omega_0}_{L^1}
   \longrightarrow0.
 \end{aligned}
\end{equation}
Passing to the limit using \eqref{uniq:stationary-product-limit} gives
$\operatorname{div}(U\omega_0)=0$ distributionally on the whole plane,
including the accumulation points of the disks.
Since $U\in H^1\cap L^\infty$, the identity
$\curl\bigl((U\cdot\nabla)U\bigr)=\operatorname{div}(U\omega_0)=0$ shows that
$(U\cdot\nabla)U$ is a distributional gradient on the simply connected
plane.  Thus $U$ also solves
the stationary momentum equation \eqref{eq:momentum}.

The stationary flow is explicit.  Let $\operatorname{Rot}(\vartheta)$
denote counterclockwise rotation through the angle $\vartheta$.
For $0<|x-c_m|<d_m$, set
\begin{equation}\label{uniq:quiet-flow}
 X^{\mathrm{quiet}}(t,x)
   =c_m+\operatorname{Rot}\left(
       \frac t2(\mathsf H\Omega_m)(|x-c_m|)\right)(x-c_m),
\end{equation}
and set $X^{\mathrm{quiet}}(t,c_m)=c_m$ and
$X^{\mathrm{quiet}}(t,x)=x$ outside the union of the disks.
Formula~\eqref{uniq:quiet-flow} preserves radii and area, has inverse obtained by
replacing $t$ with $-t$, and differentiates to $U$.
It therefore defines a measure-preserving regular Lagrangian flow
and transports $\omega_0$ exactly.  The angular velocity need not
have a finite limit at a center; the center itself remains fixed.
The Sobolev regularity of $U$ and its bounded support also give
uniqueness of this flow up to null sets by the Sobolev transport theory
of DiPerna and Lions \cite{DiPernaLions1989}.
\end{proof}

\subsection{The relative-energy identity}
\label{uniq:relative-subsection}

We next justify the energy calculation for distributional solutions.  The
argument uses only one fixed vorticity exponent above $2$ for the competing
solution.  In particular, the competitor does not need the full family of
bounds in \eqref{eq:envelope}.

\begin{lemma}[Canonical velocities in the energy class]
\label{uniq:energy-regularity}
Let $\omega$ solve
\begin{equation}\label{uniq:euler}
 \partial_t\omega+\operatorname{div}(u\omega)=0,
 \qquad u=K*\omega,
\end{equation}
on $(0,T)\times\R^2$.  Suppose that its vorticity has support in one fixed
disk, has zero total circulation, and satisfies
\begin{equation}\label{uniq:competitor-bounds}
 \omega\in L^\infty\left((0,T);L^1(\R^2)\cap L^q(\R^2)\right),
 \qquad q>2.
\end{equation}
Then the velocity satisfies
\begin{equation}\label{uniq:energy-spaces}
 u\in L^\infty\left((0,T);H^1(\R^2)\cap L^\infty(\R^2)\right)
          \cap W^{1,\infty}\left((0,T);L^2(\R^2)\right),
\end{equation}
and its momentum equation holds as the $L^2$ identity
\begin{equation}\label{uniq:projected-momentum}
 \partial_tu=-\mathbb P\left(u\cdot\nabla u\right),
\end{equation}
where $\mathbb P$ is the orthogonal projection onto divergence-free fields in
$L^2(\R^2)$.  If its initial trace is $\omega_0$ in distributions, with the
same support and integrability properties, then its $L^2$ velocity trace is
$K*\omega_0$.
\end{lemma}

\begin{proof}
Interpolation gives $\omega\in L^\infty_tL^2_x$, and
\eqref{uniq:velocity-bound} bounds $u$ uniformly.  If $B(0,R)$ contains the
vorticity support, zero circulation gives, for $|x|\geq2R$,
\begin{equation}\label{uniq:velocity-tail}
 u(t,x)=\int_{\R^2}\left(K(x-y)-K(x)\right)\omega(t,y)\idiff y,
 \qquad
 |u(t,x)|\leq C R\norm{\omega(t)}_{L^1}|x|^{-2}.
\end{equation}
The near-field bound and \eqref{uniq:velocity-tail} give
$u\in L^\infty_tL^2_x$.  The fixed-exponent Fourier multiplier estimate
$\norm{\nabla u}_{L^2}\leq C\norm{\omega}_{L^2}$ then gives its $H^1$ bound.
Applying $K*$ to \eqref{uniq:euler} yields
\begin{equation}\label{uniq:time-derivative}
 \partial_tu=-\sum_{j=1}^2(\partial_jK)*(u_j\omega)
                       \in L^\infty\left((0,T);L^2(\R^2)\right).
\end{equation}
Indeed, each $\partial_jK$ is a bounded Fourier multiplier of order zero on
$L^2$, while $u_j\omega\in L^\infty_tL^2_x$.  This proves
\eqref{uniq:energy-spaces}.

Both sides of \eqref{uniq:projected-momentum} belong to $L^2$.  Their
difference is divergence-free and has zero curl by \eqref{uniq:euler} and
the identity $\curl(u\cdot\nabla u)=\operatorname{div}(u\omega)$.
A divergence-free, curl-free $L^2$ vector field on $\R^2$ vanishes, as follows
immediately from its Fourier transform.  This proves the momentum equation.
Finally, the $L^2$ time trace supplied by \eqref{uniq:energy-spaces} is
divergence-free and has curl $\omega_0$ by the assumed distributional trace.
Its difference from $K*\omega_0$ is again a harmonic $L^2$ vector field and
therefore vanishes.
\end{proof}

\begin{lemma}[Relative energy about the quiet field]
\label{uniq:relative-energy}
Let $\omega_0\in\mathscr S$ and $U=K*\omega_0$.
Suppose that $\omega$ satisfies the hypotheses of
Lemma~\ref{uniq:energy-regularity}.  For $v=u-U$, the relative energy
$\mathcal E(t)=\norm{v(t)}_{L^2}^2/2$ is absolutely continuous and satisfies
\begin{equation}\label{uniq:relative-identity}
 \mathcal E'(t)=-\int_{\R^2}v(t,x)^{\mathsf T}S(U)(x)v(t,x)\idiff x
 \quad\text{for almost every }t\in(0,T).
\end{equation}
\end{lemma}

\begin{proof}
Lemma~\ref{uniq:stationary-lemma} and
Lemma~\ref{uniq:energy-regularity}, also applied to the stationary solution,
allow us to subtract the projected momentum equations and obtain
\begin{equation}\label{uniq:difference-momentum}
 \partial_tv=-\mathbb P\left(u\cdot\nabla v+v\cdot\nabla U\right)
                  \quad\text{in }L^2.
\end{equation}
The Hilbert-space chain rule and orthogonality of $\mathbb P$ give
\begin{equation}\label{uniq:chain-rule}
 \mathcal E'(t)=-\int_{\R^2}v\cdot(u\cdot\nabla v)\idiff x
       -\int_{\R^2}v\cdot(v\cdot\nabla U)\idiff x.
\end{equation}
The first integral is zero.  To justify this on the plane, let $\chi_R$ be a
smooth cutoff equal to $1$ on $B(0,R)$, supported in $B(0,2R)$, and satisfying
$\norm{\nabla\chi_R}_{L^\infty}\leq C/R$.  Since $u,v\in H^1\cap L^\infty$,
Sobolev approximation and incompressibility give
\begin{equation}\label{uniq:transport-cancellation}
 \int_{\R^2}\chi_R v\cdot(u\cdot\nabla v)\idiff x
       =-\frac12\int_{\R^2}|v|^2u\cdot\nabla\chi_R\idiff x,
 \qquad
 \left|\int_{\R^2}|v|^2u\cdot\nabla\chi_R\idiff x\right|
       \leq\frac{C}{R}\norm{u}_{L^\infty}\norm{v}_{L^2}^2.
\end{equation}
Letting $R\to\infty$ proves the cancellation.  The skew-symmetric part of
$\nabla U$ contributes zero to the second integral in
\eqref{uniq:chain-rule}, which proves \eqref{uniq:relative-identity}.
\end{proof}

\subsection{The Osgood criterion}
\label{uniq:osgood-subsection}

For a possibly discontinuous monotone envelope, it is convenient to use the
optimized modulus
\begin{equation}\label{uniq:optimized-modulus}
 \mu_\Theta(y)=\inf_{p\geq p_0}\Theta(p)y^{1-1/p}
       \quad(0<y\leq1),
 \qquad \mu_\Theta(0)=0.
\end{equation}
This modulus is positive, nondecreasing, concave, and continuous on $(0,1)$;
continuity at zero follows from the estimate
$\mu_\Theta(y)\leq\Theta(p_0)y^{1-1/p_0}$.  Concavity follows because each
function under the infimum is concave.  If $0<y\leq e^{-p_0}$, the admissible
choice $p=\log(1/y)$ gives
\begin{equation}\label{uniq:log-modulus}
 \mu_\Theta(y)\leq e\,y\Theta\left(\log(1/y)\right).
\end{equation}
Thus no continuity or doubling assumption on $\Theta$ is needed below.

\begin{lemma}[The optimized integral criterion]\label{uniq:modulus-equivalence}
For the modulus in \eqref{uniq:optimized-modulus}, one has
\begin{equation}\label{uniq:integral-comparison}
 e^{-1}\int_{p_0}^{\infty}\frac{\idiff p}{\Theta(p)}
 \leq\int_0^1\frac{\idiff y}{\mu_\Theta(y)}
 \leq\frac{2p_0}{\Theta(p_0)}
         +4\int_{p_0}^{\infty}\frac{\idiff p}{\Theta(p)}.
\end{equation}
Thus the two integrals in \eqref{uniq:integral-comparison} diverge or converge
together.
\end{lemma}

\begin{proof}
Define $H_\Theta(s)=\inf_{p\geq p_0}\Theta(p)e^{s/p}$ for $s\geq0$.
The substitution $y=e^{-s}$ gives
\begin{equation}\label{uniq:modulus-transform}
 \int_0^1\frac{\idiff y}{\mu_\Theta(y)}
      =\int_0^\infty\frac{\idiff s}{H_\Theta(s)}.
\end{equation}
For $s\geq p_0$, choosing $p=s$ gives $H_\Theta(s)\leq e\Theta(s)$,
which proves the lower bound.  For the upper bound, put $p_j=2^jp_0$.
If $p_j\leq p<p_{j+1}$, monotonicity of $\Theta$ gives
$e^{-s/p}/\Theta(p)\leq e^{-s/p_{j+1}}/\Theta(p_j)$.
Taking the supremum over $p$ and then enlarging that supremum to a sum yields
\begin{equation}\label{uniq:dyadic-modulus}
 \frac1{H_\Theta(s)}
     \leq\sum_{j\geq0}\frac{e^{-s/p_{j+1}}}{\Theta(p_j)}.
\end{equation}
Tonelli's theorem and monotonicity give
\begin{equation}\label{uniq:dyadic-modulus-integral}
 \begin{aligned}
 \int_0^\infty\frac{\idiff s}{H_\Theta(s)}
   \leq\sum_{j\geq0}\frac{p_{j+1}}{\Theta(p_j)}
   \leq\frac{2p_0}{\Theta(p_0)}
           +4\sum_{j\geq1}\int_{p_{j-1}}^{p_j}
                              \frac{\idiff p}{\Theta(p)}.
 \end{aligned}
\end{equation}
Together with \eqref{uniq:modulus-transform}, this proves
\eqref{uniq:integral-comparison}.
\end{proof}

\begin{theorem}[Uniqueness for stationary radial arrays]
\label{uniq:uniqueness-theorem}
Let $\omega_0\in\mathscr S\cap Y^\Theta$, with $\mathscr S$ as in
Definition~\ref{def:shields}.  Suppose that
\begin{equation}\label{uniq:divergent-integral}
 \int_{p_0}^\infty\frac{\idiff p}{\Theta(p)}=\infty.
\end{equation}
Let $\omega$ be any canonical distributional solution of
\eqref{uniq:euler} with initial trace $\omega_0$, uniformly bounded vorticity
support on each finite time interval, and
\eqref{uniq:competitor-bounds} for one fixed $q>2$.  Then
$\omega(t)=\omega_0$ for all times in its distributionally continuous
representative.  In particular, the quiet solution is the unique canonical
Lagrangian Euler evolution in the compactly supported all-finite-$p$ class.
\end{theorem}

\begin{proof}
Testing \eqref{uniq:euler} against a spatial cutoff that equals $1$ on the
common support preserves the circulation, which is initially zero by
\eqref{eq:radial-neutrality} and the $L^1$ convergence in
Lemma~\ref{uniq:stationary-lemma}.  Lemma~\ref{uniq:energy-regularity} therefore applies,
and $\mathcal E(0)=0$ for the relative energy in
Lemma~\ref{uniq:relative-energy}.  Choose $M_{\mathrm v}\geq1$ such that
$\norm{u-U}_{L^\infty_{t,x}}\leq M_{\mathrm v}$.  By
\eqref{uniq:strain-Lp}, \eqref{eq:envelope}, H\"older's
inequality, and interpolation, one has, for every $p\geq p_0$,
\begin{equation}\label{uniq:energy-p-estimate}
 \begin{aligned}
 \mathcal E'(t)\leq\norm{S(U)}_{L^p}
            \norm{v(t)}_{L^{2p/(p-1)}}^2
      \leq C_0\Theta(p)M_{\mathrm v}^{2/p}(2\mathcal E(t))^{1-1/p}.
 \end{aligned}
\end{equation}
Here $C_0$ is independent of $p$ and $t$.  At each time where
\eqref{uniq:relative-identity} holds and $v(t)\in L^2\cap L^\infty$,
H\"older's inequality applies to every $p\geq p_0$ simultaneously.

Choose a constant $\mathcal E_{\mathrm{ref}}$ larger than $2M_{\mathrm v}^2+2\sup_t\mathcal E(t)+1$ and put $y(t)=\mathcal E(t)/\mathcal E_{\mathrm{ref}}$.
Then $0\leq y(t)\leq1/2$, and \eqref{uniq:energy-p-estimate} implies
\begin{equation}\label{uniq:energy-optimized}
 y'(t)\leq2C_0\mu_\Theta(y(t))
            \quad\text{for almost every }t.
\end{equation}
The substitution $s=\log(1/y)$ in \eqref{uniq:log-modulus} gives
\begin{equation}\label{uniq:osgood-integral}
 \int_0^{e^{-p_0}}\frac{\idiff y}{\mu_\Theta(y)}
 \geq e^{-1}\int_{p_0}^\infty\frac{\idiff s}{\Theta(s)}
 =\infty.
\end{equation}
For $0<\varepsilon<1/2$, monotonicity of $\mu_\Theta$ and the chain rule
applied to $y+\varepsilon$ yield
\begin{equation}\label{uniq:epsilon-osgood}
 \int_\varepsilon^{y(t)+\varepsilon}
                \frac{\idiff s}{\mu_\Theta(s)}\leq2C_0t.
\end{equation}
If $y(t)>0$, the left-hand side tends to infinity as
$\varepsilon\downarrow0$ by \eqref{uniq:osgood-integral}, a contradiction.
Consequently $\mathcal E=0$, hence $u=U$ and $\omega=\curl u=\omega_0$.
\end{proof}

\subsection{Bounded energy-admissible momentum solutions}
\label{uniq:admissible-subsection}

The relative-energy argument also applies when no vorticity integrability is
imposed on the competing solution.  The additional hypothesis needed in that
case is an energy inequality.  We make the solution class explicit.

\begin{definition}[Bounded energy-admissible solution]
\label{uniq:admissible-definition}
A bounded energy-admissible solution with initial velocity $U\in L^2(\R^2)$
is a divergence-free vector field
\begin{equation}\label{uniq:admissible-spaces}
 u\in C_{\mathrm w}\left([0,T];L^2(\R^2)\right)
       \cap L^\infty\left((0,T)\times\R^2\right)
\end{equation}
with $u(0)=U$ that satisfies
\begin{equation}\label{uniq:admissible-weak-form}
 \int_0^T\int_{\R^2}
   \left(u\cdot\partial_t\varphi
          +(u\otimes u):\nabla\varphi\right)\idiff x\idiff t
       +\int_{\R^2}U\cdot\varphi(0)\idiff x=0
\end{equation}
for every compactly supported smooth divergence-free vector test field
$\varphi$ on $[0,T)\times\R^2$, and obeys
\begin{equation}\label{uniq:admissible-energy}
 \norm{u(t)}_{L^2}^2\leq\norm{U}_{L^2}^2
           \quad\text{for almost every }t\in(0,T).
\end{equation}
Here $C_{\mathrm w}$ denotes continuity for the weak topology of $L^2$.
\end{definition}

\begin{theorem}[Uniqueness among bounded admissible solutions]
\label{uniq:admissible-theorem}
Let $\omega_0\in\mathscr S\cap Y^\Theta$ and $U=K*\omega_0$.
Suppose that \eqref{uniq:divergent-integral} holds.
Then, for every $T>0$, each bounded energy-admissible solution
with initial velocity $U$ equals $U$ on $[0,T]$.
\end{theorem}

\begin{proof}
Let $U_\varepsilon$ be a standard spatial mollification of $U$.
Lemma~\ref{uniq:stationary-lemma} shows that $U$ has bounded support; hence
$U_\varepsilon$ is an admissible divergence-free spatial test field.
For each fixed $p>1$ with $U\in W^{1,p}$, boundedness and finite energy give
$u\otimes u\in L^\infty_tL^{p/(p-1)}_x$.  We may therefore insert
$\varphi(t,x)=\eta(t)U_\varepsilon(x)$ into
\eqref{uniq:admissible-weak-form} and pass to the limit
$\varepsilon\downarrow0$.  The resulting scalar distributional identity has
an integrable right-hand side.  Weak time continuity then identifies its
absolutely continuous representative and gives
\begin{equation}\label{uniq:cross-identity}
 \int_{\R^2}u(t)\cdot U\idiff x-\norm{U}_{L^2}^2
       =\int_0^t\int_{\R^2}(u\otimes u):\nabla U\idiff x\idiff s
       \quad\text{for every }t\in[0,T].
\end{equation}
Put $w=u-U$ and $\mathcal E(t)=\norm{w(t)}_{L^2}^2/2$.
Combining \eqref{uniq:cross-identity} with
\eqref{uniq:admissible-energy} yields
\begin{equation}\label{uniq:admissible-relative-pre}
 \mathcal E(t)\leq-\int_0^t\int_{\R^2}(u\otimes u):\nabla U\idiff x\idiff s
          \quad\text{for almost every }t.
\end{equation}
We now remove the constant and linear terms in $w$ from its right-hand side.
The stationary projected momentum equation gives
$\mathbb P(U\cdot\nabla U)=0$ in $L^2$, while $\mathbb Pw=w$.
Consequently, we obtain
\begin{equation}\label{uniq:admissible-mixed}
 \int_{\R^2}w\cdot(U\cdot\nabla U)\idiff x=0,
 \qquad
 \int_{\R^2}U\cdot(w\cdot\nabla U)\idiff x
   =\int_{\R^2}w\cdot\nabla\left(\frac{|U|^2}{2}\right)\idiff x=0.
\end{equation}
The second identity follows because $w$ is divergence-free and
$|U|^2/2$ is compactly supported and belongs to $H^1$.
The constant term vanishes by the same argument with $w$ replaced by $U$.
Thus \eqref{uniq:admissible-relative-pre} becomes
\begin{equation}\label{uniq:admissible-relative}
 \mathcal E(t)\leq-\int_0^t\int_{\R^2}(w\otimes w):\nabla U\idiff x\idiff s
       \leq F(t),
 \qquad
 F(t)=\int_0^t\left|\int_{\R^2}w^{\mathsf T}S(U)w\idiff x\right|\idiff s.
\end{equation}
Since $S(U)\in L^2$ and $w\otimes w\in L^\infty_tL^2_x$,
the function $F$ in \eqref{uniq:admissible-relative} is finite,
nonnegative, and absolutely continuous, with $F(0)=0$.
If $M_{\mathrm v}\geq1$ bounds $|w|$, the same H\"older estimate as in
\eqref{uniq:energy-p-estimate} gives
\begin{equation}\label{uniq:admissible-F}
 F'(t)\leq C_0\Theta(p)M_{\mathrm v}^{2/p}(2\mathcal E(t))^{1-1/p}
        \leq C_0\Theta(p)M_{\mathrm v}^{2/p}(2F(t))^{1-1/p}.
\end{equation}
After normalizing $F$ by a sufficiently large constant $\mathcal E_{\mathrm{ref}}$, the proof of
\eqref{uniq:energy-optimized}--\eqref{uniq:epsilon-osgood} applies without
change and gives $F=0$.  Therefore $\mathcal E=0$ almost everywhere, and the weakly
continuous representative of $u$ equals $U$ at every time.
\end{proof}

Every canonical solution satisfying
Lemma~\ref{uniq:energy-regularity} is energy admissible, in fact with equality.
Indeed, pairing \eqref{uniq:projected-momentum} with $u$ and using the cutoff
argument in \eqref{uniq:transport-cancellation} gives
$\idiff\norm{u(t)}_{L^2}^2/\idiff t=0$.  Thus
Theorem~\ref{uniq:admissible-theorem} also yields the conclusion of
Theorem~\ref{uniq:uniqueness-theorem}, while covering additional bounded weak
momentum solutions.

% Local analytic module.  This fragment uses macro-3.tex and no additional macros.
\section{The neutral cell and its response to an exterior field}
\label{cell:section}\label{sec:local-cell}

The cascade requires a cell that can remain close to a circular steady state for an arbitrarily long interval and then follow a prescribed growing orbit.  This section constructs that orbit and proves a response estimate whose constants do not depend on the length of the waiting interval.  The exterior field is prescribed only while proving the local statement.  In the finite array it will be the velocity generated by the other transported cells.

Throughout this section, $s$ denotes time in unit-cell coordinates.  The fixed constants and circular vorticity are
\begin{equation}\label{cell:unit}
 a=\frac{3}{\pi},\qquad q_{\mathrm c}=\frac1{\sqrt3},\qquad
 R_{\mathrm{in}}=q_{\mathrm c},\quad R_{\mathrm{out}}=1,\qquad
 \Omega=\frac{3a}{2}\mathbf1_{B_{q_{\mathrm c}}}-\frac a2\mathbf1_{B_1}.
\end{equation}
Thus $\Omega$ takes the values $a$, $-a/2$, and $0$.  Its circulation is zero, and its Biot--Savart velocity vanishes outside $B_1$.

\subsection{The contour equation and centered coordinates}
\label{cell:coordinates}

Write $\Delta_{\mathrm{in}}=3a/2$ and $\Delta_{\mathrm{out}}=-a/2$.  A nearby cell has the form $W=\Delta_{\mathrm{in}}\mathbf1_{D_{\mathrm{in}}}+\Delta_{\mathrm{out}}\mathbf1_{D_{\mathrm{out}}}$, where the two domains are nested and have areas $\pi/3$ and $\pi$.  Positively oriented parametrizations of their boundaries satisfy the exact velocity formula
\begin{equation}\label{cell:contour}
 U_j(\theta)=-\frac1{2\pi}\sum_{\ell\in\{\mathrm{in},\mathrm{out}\}}
 \Delta_\ell\int_0^{2\pi}
 \log\bigl|Z_j(\theta)-Z_\ell(\vartheta)\bigr|
 \partial_\vartheta Z_\ell(\vartheta)\idiff\vartheta.
\end{equation}
Here $U_j=(K*W)(Z_j)$; integration by parts over $D_\ell$ gives \eqref{cell:contour}.  The diagonal logarithm is integrable.

Let $\xi_j$ be the centroid of $D_j$ and represent the boundaries by
\begin{equation}\label{cell:radial-chart}
 Z_j(\theta)=\xi_j+R_j\bigl(1+\eta_j(\theta)\bigr)e_r(\theta),
 \qquad e_r(\theta)=(\cos\theta,\sin\theta),\qquad e_\theta=e_r^\perp.
\end{equation}
For a Fourier series $f(\theta)=\sum_k f_ke^{ik\theta}$, set
\begin{equation}\label{cell:fourier-norm}
 \|f\|_{\rho,d}=\sum_{k\in\mathbb Z}(1+|k|)^d e^{\rho|k|}|f_k|,
 \qquad d\in\{0,1,2\}.
\end{equation}
The area and centroid requirements are precisely
\begin{equation}\label{cell:constraints}
 \bigl((1+\eta_j)^2\bigr)_0=1,\qquad
 \bigl((1+\eta_j)^3\bigr)_1=
 \bigl((1+\eta_j)^3\bigr)_{-1}=0.
\end{equation}
Indeed, polar integration gives the area from the zeroth Fourier coefficient of $(1+\eta_j)^2$ and the first moment from the first Fourier coefficients of $(1+\eta_j)^3$.

Let $g_j$ contain only modes $|k|\geq2$.  The derivative of the three equations in \eqref{cell:constraints}, with respect to the missing modes $0,1,-1$, is diagonal at zero, with diagonal entries $2,3,3$.  The analytic implicit function theorem therefore gives a chart $\eta_j=H_c(g_j)$ satisfying
\begin{equation}\label{cell:centered-chart}
 P_{\geq2}H_c(g_j)=g_j,\qquad
 H_c(g_j)=g_j+O\bigl(\|g_j\|_{\rho,1}^2\bigr).
\end{equation}
The chart preserves real functions.  On the subspace of half-turn-even functions it is given explicitly by
\begin{equation}\label{cell:even-chart}
 H_c(g)=g-1+\sqrt{1-(g^2)_0},
\end{equation}
with the square-root branch equal to one at $g=0$.

Put $\xi=\xi_{\mathrm{out}}$ and $b=\xi_{\mathrm{in}}-\xi_{\mathrm{out}}$.  The reduced state is $\mathbf h=(b,g_{\mathrm{in}},g_{\mathrm{out}})$, equipped with the norms
\begin{align}
 \|\mathbf h\|_{\mathcal X_\rho}&=|b|+\sum_j\|g_j\|_{\rho,1},
 &\|\mathbf h\|_{\mathcal Z_\rho}&=|b|+\sum_j\|g_j\|_{\rho,0},\label{cell:state-norms}\\
 \|\mathbf h\|_{\mathcal Y_\rho}&=\sum_j\sum_{|k|\geq2}|k|(1+|k|)e^{\rho|k|}|(g_j)_k|.
 \label{cell:derivative-norm}
\end{align}
The last expression is a seminorm; the finite-dimensional centroid variable is not charged an angular derivative.  When complexifying these spaces, use independent coordinates $b_+=b_1+ib_2$ and $b_-=b_1-ib_2$, with norm $|b_+|+|b_-|$.

For an exterior field $E$, define the centroid velocities by
\begin{equation}\label{cell:mean-velocities}
 \mathcal V_j=\frac1{|D_j|}\int_{D_j}(K*W+E)\idiff y,
 \qquad \xi'=\mathcal V_{\mathrm{out}},\qquad
 b'=\mathcal V_{\mathrm{in}}-\mathcal V_{\mathrm{out}}.
\end{equation}
The normal equation for $\eta_j$ is $\eta_j'=V_j$, where
\begin{equation}\label{cell:normal}
 V_j=\frac1{R_j}
 \left(e_r-\frac{\partial_\theta\eta_j}{1+\eta_j}e_\theta\right)
 \cdot\bigl(U_j+E(Z_j)-\mathcal V_j\bigr).
\end{equation}
Incompressibility gives the derivative of the area constraint in \eqref{cell:constraints}, and \eqref{cell:mean-velocities} gives the derivatives of the two centroid constraints.  Since $P_{\geq2}\mathrm D H_c$ is the identity, these three differentiated constraints imply the exact identity
\begin{equation}\label{cell:lift}
 V_j=\mathrm D H_c(g_j)P_{\geq2}V_j.
\end{equation}
Consequently $g_j'=P_{\geq2}V_j$, together with \eqref{cell:mean-velocities}, is equivalent to the full normal equation.

\subsection{Analytic estimates for the exact vector field}
\label{cell:analytic-estimates}

We record the derivative loss explicitly because the response estimate will be used on intervals of unbounded length.  The relevant bounds hold uniformly for $0\leq\rho\leq1/16$ on one fixed small ball in $\mathcal X_\rho$.

\begin{lemma}\label{cell:tame-lemma}
The autonomous reduced equation can be written as $\mathbf h'=\mathcal A \mathbf h+\mathcal N(\mathbf h)$, where $\mathcal A$ is its derivative at the circular state.  The nonlinear map is holomorphic from a fixed ball of $\mathcal X_\rho$ to $\mathcal Z_\rho$ and has symmetric multilinear coefficients $\mathcal N_m$, $m\geq2$, satisfying
\begin{equation}\label{cell:coefficient-bound}
 \|\mathcal N_m(f_1,\ldots,f_m)\|_{\mathcal Z_\rho}
 \leq a C_0C_1^m\prod_{j=1}^m\|f_j\|_{\mathcal X_\rho}.
\end{equation}
For $\mathfrak X=\|\mathbf h\|_{\mathcal X_\rho}$, $\mathfrak Y=\|\mathbf h\|_{\mathcal Y_\rho}$, and the analogous quantities $\widetilde{\mathfrak X},\widetilde{\mathfrak Y}$ associated with $\widetilde{\mathbf h}$, one has
\begin{align}
 \|\mathcal N(\mathbf h)\|_{\mathcal X_\rho}
 &\leq C\mathfrak X(\mathfrak X+\mathfrak Y),\label{cell:tame-self}\\
 \|\mathcal N(\mathbf h)-\mathcal N(\widetilde{\mathbf h})\|_{\mathcal X_\rho}
 &\leq C(\mathfrak X+\widetilde{\mathfrak X})\|\mathbf h-\widetilde{\mathbf h}\|_{\mathcal Y_\rho}
 +C(\mathfrak X+\widetilde{\mathfrak X}+\mathfrak Y+\widetilde{\mathfrak Y})\|\mathbf h-\widetilde{\mathbf h}\|_{\mathcal X_\rho}.
 \label{cell:tame-difference}
\end{align}
For each fixed finite Fourier projection $P$, the right side of the corresponding projected difference estimate is bounded by $C_P(\mathfrak X+\widetilde{\mathfrak X})\|\mathbf h-\widetilde{\mathbf h}\|_{\mathcal X_\rho}$.
\end{lemma}

\begin{proof}
The norm $\|\cdot\|_{\rho,0}$ is a Banach algebra.  The analogous two-variable algebra has weight $e^{\rho(|k|+|\ell|)}$, and integration in either angular variable is bounded in this algebra.  For a self-interaction, factor the two complex chords to obtain
\begin{align}
 \log|Z_j(\theta)-Z_j(\vartheta)|
 &=\log R_j+\log|e^{i\theta}-e^{i\vartheta}|
 +\frac12\log(1+Q_+\eta_j)+\frac12\log(1+Q_-\eta_j),
 \label{cell:kernel-factor}\\
 Q_\pm\eta(\theta,\vartheta)
 &=\frac{e^{\pm i\theta}\eta(\theta)-e^{\pm i\vartheta}\eta(\vartheta)}
 {e^{\pm i\theta}-e^{\pm i\vartheta}}.
 \label{cell:divided-difference}
\end{align}
For the input $e^{ik\theta}$, the quotient $Q_+$ is the sum of $|k+1|$ monomials, and $Q_-$ is the sum of $|k-1|$ monomials.  The sum of the absolute frequencies of every nonzero monomial is $|k|$.  This gives the estimate
\begin{equation}\label{cell:divided-bound}
 \|Q_\pm\eta\|_{\rho,0;2\mathrm{var}}\leq\|\eta\|_{\rho,1}.
\end{equation}
One further angular derivative costs at most $C\|\eta\|_{\rho,2}$.  Thus the two logarithmic series in \eqref{cell:kernel-factor} converge in the two-variable algebra on a fixed small ball.  Their holomorphic extensions use the same complex input in both factors; they do not conjugate that input.

The unperturbed logarithm in \eqref{cell:kernel-factor} has nonzero Fourier coefficients $-1/(2|k|)$.  Its convolution with a differentiated contour is therefore a bounded Fourier multiplier.  For the interaction between distinct interfaces, the unperturbed chord has a reciprocal geometric expansion with ratio $q_{\mathrm c}$, convergent in the two-variable algebra whenever $q_{\mathrm c}e^{2\rho}<1$.  Relative centroid displacements enter this separated denominator analytically.

Every term of the resulting normal equation is a product of divided differences, a differentiated contour, and at most one normal-projection factor $\partial_\theta\eta/(1+\eta)$.  Each occurrence of an input is charged at most one angular derivative.  The centroid averages have the same property: the self-average vanishes by antisymmetry, and the mutual average can be computed by integrating the outer boundary formula over the inner domain.  The integration points remain a positive distance from the outer boundary.  The chart corrections involve only analytic products and fixed Fourier projections.  These facts prove holomorphy from $\mathcal X_\rho$ to $\mathcal Z_\rho$ on a radius independent of $\rho$.  Cauchy's inequality and polarization then give \eqref{cell:coefficient-bound}; the factor $m^m/m!\leq e^m$ is absorbed into $C_1^m$.

To control an additional output derivative, use rotational covariance.  Physical rotation through $\varphi$ acts on the complexified state by
\begin{equation}\label{cell:rotation-action}
 S_\varphi(b_+,b_-,g_{\mathrm{in}},g_{\mathrm{out}})
 =(e^{i\varphi}b_+,e^{-i\varphi}b_-,
 g_{\mathrm{in}}(\cdot-\varphi),g_{\mathrm{out}}(\cdot-\varphi)).
\end{equation}
The generator is $\mathcal J\mathbf h=(ib_+,-ib_-,-\partial_\theta g_{\mathrm{in}},-\partial_\theta g_{\mathrm{out}})$.  Equivariance gives, first for finite Fourier series,
\begin{equation}\label{cell:rotation-differentiate}
 \mathcal J\mathcal N_m(f_1,\ldots,f_m)
 =\sum_{j=1}^m\mathcal N_m(f_1,\ldots,\mathcal J f_j,\ldots,f_m).
\end{equation}
For a reduced state $\mathbf f$, the inequalities
$\|\mathcal J\mathbf f\|_{\mathcal X_\rho}\leq
C(\|\mathbf f\|_{\mathcal X_\rho}+\|\mathbf f\|_{\mathcal Y_\rho})$
and
$\|\mathbf f\|_{\mathcal X_\rho}\leq
C(\|\mathbf f\|_{\mathcal Z_\rho}+\|\mathcal J\mathbf f\|_{\mathcal Z_\rho})$,
combined with \eqref{cell:coefficient-bound}, prove a bound with exactly one additional input derivative.  Summing over $m$ gives \eqref{cell:tame-self}.  Telescoping the multilinear differences gives \eqref{cell:tame-difference}.  A fixed finite projection is bounded from $\mathcal Z_\rho$ to $\mathcal X_\rho$, so its difference estimate follows directly from \eqref{cell:coefficient-bound}.
\end{proof}

\subsection{Centroids and exterior fields}
\label{cell:exterior}

The signed impulse $I=\int yW(y)\idiff y$ satisfies
\begin{equation}\label{cell:impulse}
 I=\frac{a\pi}{2}b,\qquad
 b'=\frac2{a\pi}\int EW\idiff y.
\end{equation}
The second identity follows from $\int(K*W)W=0$.  In particular, the autonomous relative-centroid equation is identically zero.  The linear part of the outer centroid's self-induced velocity is $-ab^\perp/4$.  Indeed, the velocity of the unit disk in its interior is $y^\perp/2$, so antisymmetry and integration over the inner disk give
\begin{equation}\label{cell:centroid-drift}
 \frac1\pi\int_{D_{\mathrm{out}}}K*W\idiff y
 =-\frac a4b^\perp+O\bigl(\|\mathbf h\|_{\mathcal X_\rho}^2\bigr).
\end{equation}

Suppose $E$ is divergence free and curl free in the radius-four tube about a reference trajectory $\xi_{\mathrm{ref}}'=E(s,\xi_{\mathrm{ref}})$.  Assume that, uniformly for $|\xi-\xi_{\mathrm{ref}}|\leq1$, its complex representation has the expansion
\begin{equation}\label{cell:exterior-coefficients}
 E_{\mathbb C}(s,\xi+y)=E_{\mathbb C}(s,\xi)
 +\sum_{\ell\geq1}v_\ell(s,\xi)\overline y^{\,\ell},
 \qquad \sum_{\ell\geq1}2^\ell|v_\ell(s,\xi)|\leq\mathfrak a(s).
\end{equation}
Here $\mathfrak a\geq0$ is locally integrable, and the constant translation is locally integrable in time.  A gradient bound on a slightly larger harmonic tube implies \eqref{cell:exterior-coefficients} by Cauchy's estimate.

At the circular state the exterior shape forcing is
\begin{equation}\label{cell:circle-force}
 F_j(\theta)=\real\sum_{\ell\geq1}v_\ell R_j^{\ell-1}e^{-i(\ell+1)\theta},
 \qquad F_b=0.
\end{equation}
The relative-centroid component vanishes because a harmonic function has the same mean on either concentric disk.  If $\mathcal F_E(\mathbf h)$ denotes the entire exterior reduced field and $F=\mathcal F_E(0)$, then we have
\begin{equation}\label{cell:exterior-tame}
 \|F\|_{\mathcal X_\rho}\leq C\mathfrak a,\qquad
 \|\mathcal F_E(\mathbf h)-F\|_{\mathcal X_\rho}
 \leq C\mathfrak a\bigl(\|\mathbf h\|_{\mathcal X_\rho}+\|\mathbf h\|_{\mathcal Y_\rho}\bigr).
\end{equation}
To verify the second estimate, expand each term of \eqref{cell:exterior-coefficients} at the chart evaluation points.  These points have modulus at most a fixed number smaller than $3/2$ on a sufficiently small complex chart.  The degree-$m$ coefficient is bounded in $\mathcal Z_\rho$ by
\begin{equation}\label{cell:exterior-degree}
 C^{m+1}\sum_{\ell\geq1}\ell(3/2)^\ell|v_\ell|
 \prod_{j=1}^m\|f_j\|_{\mathcal X_\rho}.
\end{equation}
The rotational weight of $v_\ell$ is $\ell+1$.  Applying \eqref{cell:rotation-differentiate} to this joint coefficient--state map charges either an extra state derivative or the coefficient sum with one extra factor $\ell+1$.  Both coefficient sums are bounded by $C\mathfrak a$ because of \eqref{cell:exterior-coefficients}.  Summation gives \eqref{cell:exterior-tame}.  Before applying the rotation generator, \eqref{cell:exterior-degree} also gives the derivative-free fixed-mode bound
\begin{equation}\label{cell:exterior-projection}
 \|P(\mathcal F_E(\mathbf h)-F)\|_{\mathcal X_\rho}
 \leq C_P\mathfrak a\|\mathbf h\|_{\mathcal X_\rho}.
\end{equation}

The center estimates require no time derivative of $E$.  Subtracting the circular shield centered at $\xi$ in \eqref{cell:impulse}, and using its harmonic mean cancellation, gives
\begin{align}
 |b'|&\leq C\mathfrak a\|\mathbf h\|_{\mathcal X_\rho},\label{cell:relative-center-bound}\\
 |(\xi-\xi_{\mathrm{ref}})'|
 &\leq C\mathfrak a|\xi-\xi_{\mathrm{ref}}|+C|b|
 +C\|\mathbf h\|_{\mathcal X_\rho}^2+C\mathfrak a\|\mathbf h\|_{\mathcal X_\rho}.
 \label{cell:reference-center-bound}
\end{align}
For later use, a curvature estimate improves \eqref{cell:relative-center-bound}.  Let $\mathfrak c(s)=\sup_{|y-\xi|\leq2}|\mathrm D^2 E(s,y)|$.  Since the circular shield has zero circulation and impulse, Taylor's formula yields
\begin{equation}\label{cell:impulse-curvature}
 I'=\mathrm D E(s,\xi)I+R_I,\qquad
 |R_I|\leq C\mathfrak c(s)\|\mathbf h\|_{\mathcal X_\rho}.
\end{equation}
More explicitly, subtract $E(\xi)+\mathrm D E(\xi)(y-\xi)$ and integrate the remainder against $W-\Omega(\cdot-\xi)$; its $L^1$ norm is at most $C\|\mathbf h\|_{\mathcal X_\rho}$.  If $I(0)=0$, $\int_0^H\mathfrak a(v)\idiff v\leq1$, $\mathfrak c\leq\mathfrak c_*$, and $\|\mathbf h(s)\|_{\mathcal X_\rho}\leq CAe^{\sigma s}$, then we have
\begin{equation}\label{cell:small-impulse}
 |I(s)|\leq C\int_0^s\mathfrak c(v)Ae^{\sigma v}\idiff v
 \leq C\mathfrak c_*Ae^{\sigma s}.
\end{equation}
After leaving the small chart, suppose the support stays within radius $C(1+s)$ of an exterior reference particle, the accumulated exterior strain is bounded, and $\mathfrak c\leq\mathfrak c_*$ on this support.  Taylor's formula without circular subtraction gives
\begin{equation}\label{cell:large-time-impulse}
 |I(s)|\leq C\bigl(|I(0)|+\mathfrak c_*(1+s)^3\bigr).
\end{equation}
The cubic power follows by integrating the squared support radius.

\subsection{The growing orbit}
\label{cell:orbit}

For a positive Fourier mode $k$, a fractional radial displacement of the
circle $R_j$ has the following streamfunction at radial distance
$\varrho>0$:
\begin{equation}\label{cell:linear-stream}
 \psi_k(\varrho)=-\frac{\Delta_jR_j^2}{2k}(\eta_j)_k
 \left(\frac{\min(\varrho,R_j)}{\max(\varrho,R_j)}\right)^k.
\end{equation}
Continuity and the radial derivative jump verify \eqref{cell:linear-stream}.  The radial velocity is $-\varrho^{-1}\partial_\theta\psi$.  Since the two base angular velocities are $a/2$ and $0$, the shape block of $\mathcal A$ is
\begin{equation}\label{cell:linear-block}
 \mathcal A_k=\frac i2
 \begin{pmatrix}
 3a/2-ka&-(a/2)q_{\mathrm c}^{k-2}\\
 (3a/2)q_{\mathrm c}^{k+2}&-a/2
 \end{pmatrix},\qquad k\geq2,
 \qquad \mathcal A_{-k}=\overline{\mathcal A_k}.
\end{equation}
The relative-centroid block is zero.  The discriminant of the real matrix inside $i/2$ is $a^2\bigl((k-2)^2-3^{1-k}\bigr)$.  At $k=2$, the growing eigenpair is
\begin{equation}\label{cell:eigenpair}
 \lambda=\sigma-\frac{ia}{4},\qquad
 \sigma=\frac{a}{4\sqrt3},\qquad
 v_*=(1,i/\sqrt3)^{\mathsf T}.
\end{equation}
The other eigenvalue has real part $-\sigma$.  For $k\geq3$, both eigenvalues are distinct and imaginary.  The spectral projections are uniformly bounded in $k$: their diagonal gap grows linearly, their off-diagonal entries decay geometrically, and the remaining finitely many blocks are nondegenerate.

Let $P_{k,+},P_{k,-}$ be these spectral projections, and write
$g_k=((g_{\mathrm{in}})_k,(g_{\mathrm{out}})_k)^{\mathsf T}$
for the two-component Fourier coefficient.
Replacing the shape summands in \eqref{cell:state-norms} by
\begin{equation}\label{cell:spectral-norm}
 \sum_{|k|\geq2}(1+|k|)e^{\rho|k|}
 \bigl(|P_{k,+}g_k|+|P_{k,-}g_k|\bigr)
\end{equation}
defines an equivalent norm, denoted by $\|\cdot\|_{\mathcal X_\rho,\mathrm{sp}}$.  Define its derivative seminorm by inserting $|k|$.  The linear upper Dini derivative is at most $\sigma$ times this norm.  On the complement of the two real growing coordinates, it is nonpositive.  Thus norm equivalence does not enlarge the linear growth rate.

\begin{proposition}[The unit growing orbit]\label{cell:donor}
There are $a_0>0$ and a real-analytic contour map $\mathcal G(z)$ for $|z|<a_0$, with $\mathcal G(0)=0$, whose associated vorticity is denoted by $G(z)$.  The vorticity $G(z)$ has the two values and areas in \eqref{cell:unit}, is half-turn symmetric, and has zero circulation and impulse.  For every $0<a_*\leq a_0$, the family
\begin{equation}\label{cell:ancient-orbit}
 W_*(s)=G(a_*e^{\lambda s}),\qquad s\leq0,
\end{equation}
is an exact Euler solution and has a unique bounded-vorticity continuation for $s>0$.  Its barred quadrupole $Q_*(s)=\int\overline y^{\,2}W_*(s,y)\idiff y$ satisfies
\begin{align}
 Q_*(s)&=(1-i\sqrt3)a_*e^{\lambda s}
 +O(a_*^2e^{2\sigma s}),&&s\leq0,\label{cell:quadrupole-ancient}\\
 |Q_*(s)|&\leq C(1+s)^2,&&s\geq0.\label{cell:quadrupole-future}
\end{align}
The constants are uniform for $0<a_*\leq a_0$.  The solution from $G(a_*\zeta)$, $|\zeta|=1$, has barred quadrupole $\zeta Q_*(s)$ at every time.
\end{proposition}

\begin{proof}
Restrict the exact reduced equation to half-turn-even shapes and $b=0$.  With $z,\widetilde z$ temporarily independent complex variables, seek
\begin{equation}\label{cell:manifold-series}
 \mathcal G_{\mathrm{sh}}(z,\widetilde z)
 =zv_*e^{2i\theta}+\widetilde z\overline{v_*}e^{-2i\theta}
 +\sum_{\mu+\nu\geq2}z^\mu \widetilde z^\nu g_{\mu\nu}.
\end{equation}
At degree $m=\mu+\nu$, the invariance equation determines $g_{\mu\nu}$ by applying $\bigl((\mu\lambda+\nu\overline\lambda)-\mathcal A\bigr)^{-1}$ to a polynomial in lower-degree coefficients.  The spectral calculation gives
\begin{equation}\label{cell:resolvent}
 \left\|\bigl((\mu\lambda+\nu\overline\lambda)-\mathcal A\bigr)^{-1}f\right\|_{\rho,d}
 \leq\frac C{(m-1)\sigma}\|f\|_{\rho,d}.
\end{equation}
Rotational equivariance confines $g_{\mu\nu}$ to angular frequency $2(\mu-\nu)$.  Recovering one derivative therefore costs at most $1+2m$, while $(1+2m)/(m-1)\leq5$ for $m\geq2$.  The inverse in \eqref{cell:resolvent} consequently compensates for the derivative loss in \eqref{cell:coefficient-bound}.

Let $B_m$ be the sum of the $\mathcal X_\rho$ norms of the degree-$m$ coefficients.  Their recursion is majorized coefficientwise by the analytic scalar equation
\begin{equation}\label{cell:majorant}
 B(\zeta_0)=B_1\zeta_0+
 \frac{5CC_0a}{\sigma}\frac{\bigl(C_1B(\zeta_0)\bigr)^2}{1-C_1B(\zeta_0)}.
\end{equation}
The derivative of the right nonlinear term with respect to $B$ vanishes at zero.  The implicit function theorem gives a convergent power series with nonnegative coefficients.  Writing $B(\zeta_0)=\sum_{m\geq1}\beta_m\zeta_0^m$, coefficient induction gives $B_m\leq\beta_m$.  Hence \eqref{cell:manifold-series} converges in $\mathcal X_\rho$.  On a smaller amplitude disk, angular and amplitude derivatives introduce only polynomial factors in $m$, absorbed by the strict ratio of radii.  The differentiated series therefore satisfies the full invariance equation.  Conjugating that recursion gives reality when $\widetilde z=\overline z$.

Set $\mathcal G(z)=\mathcal G_{\mathrm{sh}}(z,\overline z)$ and recover the interfaces through \eqref{cell:even-chart}.  Equation \eqref{cell:lift} gives the exact normal equation along \eqref{cell:ancient-orbit}.  Section~\ref{cell:material} below identifies such analytic contours with their transported vorticity.  Half-turn symmetry fixes both centroids at zero.  The bounded-vorticity Euler theorem \cite{Yudovich1963,CrippaStefani2024} supplies the positive-time continuation.  Conservation of the $L^1$ and $L^\infty$ norms gives a uniform velocity bound and support radius at most $C(1+s)$, which proves \eqref{cell:quadrupole-future}.

Polar integration gives
\begin{equation}\label{cell:polar-quadrupole}
 \int_{D_j}\overline y^{\,2}\idiff y
 =2\pi R_j^4(\eta_j)_2+O(\|\eta_j\|_{\rho,1}^2).
\end{equation}
Substituting \eqref{cell:eigenpair} and the jumps in \eqref{cell:unit} gives \eqref{cell:quadrupole-ancient}.  It also gives $\|G(z)-\Omega\|_1\leq C|z|$ and $\|G(z)-\Omega\|_p\leq C_p|z|^{1/p}$.  The analyticity asserted in the proposition concerns the contour map; no analyticity as an $L^p$-valued vorticity map is needed.

Finally, physical rotation through $\varphi$ replaces $z$ by $e^{-2i\varphi}z$ and multiplies the barred quadrupole by $e^{-2i\varphi}$.  The coefficient construction is equivariant, and uniqueness extends this covariance to every positive time.  Taking $e^{-2i\varphi}=\zeta$ proves the phase assertion.
\end{proof}

The transfer coefficient must include the whole donor history.  The following consequence also fixes the order in which its constants can be chosen.

\begin{corollary}[The complete donor transform]\label{cell:transform}
For $1/2\leq r<1$ and $0<a_*\leq a_0$, define
\begin{equation}\label{cell:transform-definition}
 L(r,a_*)=\int_{\R}e^{-r\lambda s}Q_*(s)\idiff s.
\end{equation}
This integral is absolutely convergent.  With $\delta=1-r$, one has
\begin{align}
 \left|L(r,a_*)-\frac{a_*}{\sigma\delta}\right|
 &\leq C_0(1+a_*^2),\label{cell:transform-asymptotic}\\
 \int_{\R}e^{-r\sigma s}|Q_*(s)|\idiff s
 &\leq\frac{2a_*}{\sigma\delta}+C_0(1+a_*^2).
 \label{cell:absolute-transform}
\end{align}
There is $c_{\mathrm{tr}}>0$ such that $0<\delta\leq c_{\mathrm{tr}}a_*$ implies
\begin{equation}\label{cell:transform-ratio}
 |L(r,a_*)|\geq\frac{a_*}{2\sigma\delta},\qquad
 \int_{\R}e^{-r\sigma s}|Q_*(s)|\idiff s
 \leq5|L(r,a_*)|.
\end{equation}
\end{corollary}

\begin{proof}
Since $\lambda=\sigma(1-i\sqrt3)$, the leading term in \eqref{cell:quadrupole-ancient} contributes
\begin{equation}\label{cell:transform-leading}
 (1-i\sqrt3)a_*\int_{-\infty}^0e^{(1-r)\lambda s}\idiff s
 =\frac{a_*}{\sigma\delta}.
\end{equation}
The ancient error is bounded by $Ca_*^2\int_{-\infty}^0e^{(2-r)\sigma s}\idiff s\leq Ca_*^2$.  The positive-time contribution is bounded by $C\int_0^\infty e^{-(\sigma/2)s}(1+s)^2\idiff s$.  This proves \eqref{cell:transform-asymptotic}; taking absolute values gives \eqref{cell:absolute-transform}.  Choose $c_{\mathrm{tr}}$ so that $C_0(1+a_0^2)\leq a_*/(2\sigma\delta)$ whenever $\delta\leq c_{\mathrm{tr}}a_*$.  The two inequalities in \eqref{cell:transform-ratio} follow.  In particular, the absolute-to-net input constant can be fixed before $a_*$ and $\delta$ are chosen.
\end{proof}

\subsection{Uniform response during the waiting interval}
\label{cell:forced}

The next proposition does not assume that the net growing input is nonzero.  This distinction allows its use while an input coefficient is still being assembled.

\begin{proposition}[Uniform forced response]\label{cell:response}
Fix $K_{\mathrm{in}}\geq1$ and $0<\rho_0\leq1/32$.  There are $\eta_0,\varepsilon_0>0$ and $C_{\mathrm{in}}<\infty$, independent of $H$, with the following property.  Suppose the exterior satisfies \eqref{cell:exterior-coefficients} on $0\leq s\leq H$.  Start from the quiet cell or from $G(z_0)$, with both centroids at $\xi_{\mathrm{ref}}(0)$.  Assume, for some $A>0$, that
\begin{equation}\label{cell:response-assumptions}
 |z_0|+\int_0^He^{-\sigma s}\mathfrak a(s)\idiff s\leq K_{\mathrm{in}}A,
 \qquad \int_0^H\mathfrak a(s)\idiff s\leq\varepsilon_0,
 \qquad Ae^{\sigma H}\leq\eta_0.
\end{equation}
For the quiet cell, take $z_0=0$.  The actual transported cell has analytic contours throughout this interval and an absolutely continuous decreasing radius satisfying
\begin{equation}\label{cell:exact-radius}
 \rho(0)=\rho_0,\qquad
 \rho'=-C_{\mathrm{rad}}\bigl(\|\mathbf h\|_{\mathcal X_\rho,\mathrm{sp}}+\mathfrak a\bigr),
 \qquad \rho\geq\rho_0/2.
\end{equation}
Writing $\mathfrak X(s)=\|\mathbf h(s)\|_{\mathcal X_{\rho(s)},\mathrm{sp}}$, one has
\begin{equation}\label{cell:response-bound}
 \mathfrak X(s)\leq C_{\mathrm{in}}Ae^{\sigma s},\qquad
 \int_0^H\mathfrak X(s)\idiff s\leq C_{\mathrm{in}}\eta_0/\sigma,
 \qquad \sup_{0\leq s\leq H}|\xi(s)-\xi_{\mathrm{ref}}(s)|<\tfrac12.
\end{equation}
There are unique coordinates $\mathbf h=\mathcal G(z)+w$ with $P_+w=0$, where
\begin{equation}\label{cell:growing-projection}
 P_+\mathbf h=(1/2,-i\sqrt3/2)(g_{\mathrm{in}},g_{\mathrm{out}})_2.
\end{equation}
They satisfy
\begin{align}
 \|w(s)\|_{\mathcal X_{\rho(s)},\mathrm{sp}}
 &\leq C_{\mathrm{in}}\int_0^s\mathfrak a(v)\idiff v,
 \label{cell:complement-bound}\\
 \left|e^{-\lambda s}z(s)-z_0
 -\int_0^se^{-\lambda v}f_{\mathrm{gr}}(v)\idiff v\right|
 &\leq C_{\mathrm{in}}A(1+s)\int_0^s\mathfrak a(v)\idiff v,
 \label{cell:relative-response}
\end{align}
where the exact growing projection of the circle forcing is
\begin{equation}\label{cell:growing-force}
 f_{\mathrm{gr}}(s)=\frac12e^{-i\pi/3}\overline{v_1(s,\xi(s))}.
\end{equation}
\end{proposition}

\begin{proof}
All norms in this proof are the adapted norms of \eqref{cell:spectral-norm}.  Differentiating the Fourier weights in \eqref{cell:exact-radius} produces $-C_{\mathrm{rad}}(\mathfrak X+\mathfrak a)\mathfrak Y$.  Lemma~\ref{cell:tame-lemma} and \eqref{cell:exterior-tame} show that a fixed sufficiently large $C_{\mathrm{rad}}$ absorbs every term with one additional derivative.  The remaining upper Dini derivative obeys
\begin{equation}\label{cell:norm-differential}
 \mathrm D^+\mathfrak X\leq\sigma \mathfrak X+C\mathfrak X^2+C\mathfrak a\mathfrak X+C\mathfrak a.
\end{equation}
Set $\mathfrak Z=e^{-\sigma s}\mathfrak X$ and bootstrap $\mathfrak Z\leq C_*A$.  The two nonlinear integrals in the resulting inequality are bounded by
\begin{equation}\label{cell:bootstrap-errors}
 C C_*^2A^2\int_0^se^{\sigma v}\idiff v
 \leq \frac{CC_*^2}{\sigma}A\eta_0,
 \qquad
 CC_*A\int_0^s\mathfrak a(v)\idiff v
 \leq CC_*A\varepsilon_0.
\end{equation}
Choose $C_*$ proportional to $K_{\mathrm{in}}$, and then choose $\eta_0,\varepsilon_0$ small.  Equations \eqref{cell:response-assumptions} and \eqref{cell:bootstrap-errors} strictly improve the bootstrap and give the first two bounds in \eqref{cell:response-bound}.  Shrinking these two small constants once more bounds the total radius loss by $\rho_0/2$.

The center tube has the same uniformity in $H$.  From $b(0)=0$ and \eqref{cell:relative-center-bound}, Fubini gives
\begin{align}
 \int_0^H|b(s)|\idiff s
 &\leq C_{\mathrm{in}} A\int_0^H(H-v)\mathfrak a(v)e^{\sigma v}\idiff v\notag\\
 &\leq C_{\mathrm{in}}\eta_0\int_0^H(H-v)e^{-\sigma(H-v)}\mathfrak a(v)\idiff v
 \leq\frac{C_{\mathrm{in}}\eta_0\varepsilon_0}{e\sigma}.
 \label{cell:center-waiting}
\end{align}
Also, $\int_0^H\mathfrak X(s)^2\idiff s\leq C_{\mathrm{in}}\eta_0^2$ and $\int_0^H\mathfrak a(s)\mathfrak X(s)\idiff s\leq C_{\mathrm{in}}\eta_0\varepsilon_0$.  Gronwall's inequality in \eqref{cell:reference-center-bound} gives
\begin{equation}\label{cell:center-closure}
 \sup_{0\leq s\leq H}|\xi-\xi_{\mathrm{ref}}|
 \leq C_{\mathrm{in}}(\eta_0^2+\eta_0\varepsilon_0)<\tfrac12.
\end{equation}
The factor $(H-v)e^{-\sigma(H-v)}$ in \eqref{cell:center-waiting} is essential to the independence from $H$.

We next construct a solution for which these estimates are valid.  Begin with radius $2\rho_0\leq1/16$, and choose the small constants so that the preceding radius loss is at most $\rho_0/4$.  Let $\Pi_J$ retain the shape modes $2\leq|k|\leq J$ and both relative-center coordinates.  Solve the finite-dimensional Carath\'eodory system
\begin{equation}\label{cell:galerkin}
 g_{j,J}'=\Pi_JP_{\geq2}V_j,
 \qquad \xi_{j,J}'=\mathcal V_j,
 \qquad \overline\rho_J'
 =-C_{\mathrm{rad}}\bigl(\|\mathbf h_J\|_{\mathcal X_{\overline\rho_J},\mathrm{sp}}+\mathfrak a\bigr),
 \quad \overline\rho_J(0)=2\rho_0,
\end{equation}
with $g_{j,J}(0)=(\Pi_J\mathcal G(z_0))_j$ and $\xi_{j,J}(0)=\xi_{\mathrm{ref}}(0)$, and with the exact chart $\eta_{j,J}=H_c(g_{j,J})$ imposed at every state.  For the quiet datum, set all initial shape coefficients to zero.  The projections commute with $\mathcal A$ and are contractions in the adapted norms.  The center identities remain exact because the represented domains have their prescribed areas and centroids, and antisymmetry gives $I_J' =\int EW_J$.  By \eqref{cell:lift}, the normal residual is precisely
\begin{equation}\label{cell:galerkin-residual}
 V_j-\partial_s\eta_{j,J}
 =\mathrm D H_c(g_{j,J})(\Id-\Pi_J)P_{\geq2}V_j.
\end{equation}
Thus the finite approximation has no extra area or centroid drift.

The estimates already proved prevent finite-dimensional exit up to $H$, and $\overline\rho_J\geq7\rho_0/4$.  On the fixed smaller radius $3\rho_0/2$, every fixed number of angular derivatives is uniformly bounded.  The time derivatives have a common integrable majorant, with the locally integrable reference translation added for the common center.  Compact weighted Fourier embeddings and Arzel\`a--Ascoli give convergence on that smaller radius.  The discarded tail in \eqref{cell:galerkin-residual} tends to zero exponentially, so the limit solves the exact reduced and normal equations.  In particular, it has the fixed-radius bound $\|\mathbf h(s)\|_{\mathcal X_{3\rho_0/2},\mathrm{sp}}\leq C_{\mathrm{in}}Ae^{\sigma s}$.

For this exact solution, now solve the scalar equation \eqref{cell:exact-radius} with initial radius $\rho_0$.  Its right side is locally Lipschitz in $\rho$ because of the unused analytic margin, and its total loss is at most $\rho_0/4$.  All Dini inequalities can therefore be applied at this exact radius.  No convergence of the approximate norms at their largest radii is needed.  Uniqueness follows by comparing two solutions at a smaller radius with derivative $-C(\mathfrak X_1+\mathfrak X_2+\mathfrak a)$.  The difference estimate \eqref{cell:tame-difference} absorbs the single derivative loss.  Cauchy's estimate on the unused exterior tube controls the difference of the coefficients at two centers by $C\mathfrak a|\xi_1-\xi_2|$.  Adding the center difference to the norm and applying Gronwall gives uniqueness on every finite $H$.

It remains to prove the relative estimate.  The map $z\mapsto P_+\mathcal G(z)=z+O(|z|^2)$ is a local real-analytic diffeomorphism.  It provides the stated decomposition with $P_+w=0$.  Let $H_+(z)=P_+\mathrm D\mathcal G(z)$, interpreted as a map between two-dimensional real spaces.  Subtracting the exact invariance equation gives
\begin{align}
 z'&=\lambda z+H_+(z)^{-1}P_+\mathcal R,\notag\\
 w'&=\mathcal Aw+
 \bigl(\Id-\mathrm D\mathcal G(z)H_+(z)^{-1}P_+\bigr)\mathcal R,
 \label{cell:flattened-equations}\\
 \mathcal R&=\mathcal N(\mathcal G(z)+w)-\mathcal N(\mathcal G(z))+\mathcal F_E(\mathbf h).
 \notag
\end{align}
Put $\mathfrak W=\|w\|_{\mathcal X_\rho,\mathrm{sp}}$ and $\mathfrak W_1=\|w\|_{\mathcal Y_\rho,\mathrm{sp}}$.  The manifold's unused analytic radius gives $\|\mathcal G(z)\|_{\mathcal Y_\rho}\leq C|z|$, while $|z|+\mathfrak W\leq C\mathfrak X$.  Hence \eqref{cell:tame-difference} gives
\begin{equation}\label{cell:flattened-tame}
 \|\mathcal N(\mathcal G(z)+w)-\mathcal N(\mathcal G(z))\|_{\mathcal X_\rho}
 \leq C\mathfrak X \mathfrak W_1+C\mathfrak X \mathfrak W.
\end{equation}
The exterior contributes at most $C\mathfrak a\mathfrak W_1+C\mathfrak a$.  The finite-mode correction in \eqref{cell:flattened-equations} is bounded by $C(\mathfrak X \mathfrak W+\mathfrak a)$ using the derivative-free bounds in Lemma~\ref{cell:tame-lemma} and \eqref{cell:exterior-projection}.  The radius derivative absorbs $C(\mathfrak X+\mathfrak a)\mathfrak W_1$, and the linear part is nonexpanding on $P_+w=0$.  Therefore, we obtain
\begin{equation}\label{cell:complement-differential}
 \mathrm D^+\mathfrak W\leq C(\mathfrak X+\mathfrak a)\mathfrak W+C\mathfrak a,
 \qquad \mathfrak W(0)=0.
\end{equation}
Its Gronwall factor is bounded independently of $H$, proving \eqref{cell:complement-bound}.

The $k=2$ Fourier coefficient of \eqref{cell:circle-force} is $\overline{v_1}/2$ on both circles.  Applying \eqref{cell:growing-projection} gives \eqref{cell:growing-force}.  The same fixed-mode estimates in the first equation of \eqref{cell:flattened-equations} give
\begin{equation}\label{cell:amplitude-remainder}
 z'=\lambda z+f_{\mathrm{gr}}+R_z,
 \qquad |R_z|\leq C\bigl((|z|+\mathfrak W)\mathfrak W+\mathfrak a(|z|+\mathfrak W)\bigr).
\end{equation}
Fubini and \eqref{cell:response-assumptions} imply
\begin{equation}\label{cell:weighted-complement}
 \int_0^He^{-\sigma s}\mathfrak W(s)\idiff s
 \leq\frac{C_{\mathrm{in}}}{\sigma}\int_0^He^{-\sigma v}\mathfrak a(v)\idiff v
 \leq C_{\mathrm{in}}A.
\end{equation}
Using $|z(s)|\leq C_{\mathrm{in}}Ae^{\sigma s}$, the four terms in \eqref{cell:amplitude-remainder} are estimated by
\begin{align}
 \int_0^se^{-\sigma v}|z(v)|\mathfrak W(v)\idiff v
 &\leq C_{\mathrm{in}}As\int_0^s\mathfrak a(v)\idiff v,\notag\\
 \int_0^se^{-\sigma v}\mathfrak W(v)^2\idiff v
 &\leq C_{\mathrm{in}}A\int_0^s\mathfrak a(v)\idiff v,\notag\\
 \int_0^se^{-\sigma v}\mathfrak a(v)\bigl(|z(v)|+\mathfrak W(v)\bigr)\idiff v
 &\leq C_{\mathrm{in}}A\int_0^s\mathfrak a(v)\idiff v.
 \label{cell:weighted-remainders}
\end{align}
Integrating the amplitude equation proves \eqref{cell:relative-response}.  All autonomous nonlinear terms were removed through the exact manifold; no forcing error independent of the exterior remains.
\end{proof}

For an exterior generated by a real vorticity $\omega_{\mathrm{ext}}$, the complex Biot--Savart formula determines the transfer phase without an adjoint normalization.  The complex velocity and its first coefficient are
\begin{equation}\label{cell:force-biot-savart}
 E_{\mathbb C}(x)=\frac i{2\pi}\int\frac{\omega_{\mathrm{ext}}(y)}{\overline x-\overline y}\idiff y,
 \qquad
 \overline{v_1(x)}=\frac i{2\pi}\int\frac{\omega_{\mathrm{ext}}(y)}{(x-y)^2}\idiff y,
\end{equation}
and the growing forcing is
\begin{equation}\label{cell:force-exact}
 f_{\mathrm{gr}}(x)=\frac{e^{i\pi/6}}{4\pi}
 \int\frac{\omega_{\mathrm{ext}}(y)}{(x-y)^2}\idiff y.
\end{equation}
A change from physical time to unit receiver time multiplies \eqref{cell:force-exact} by that receiver's time scale.  If a source is supported in $|y|\leq|d|/2$ and $M_m=\int y^m\omega_{\mathrm{ext}}(y)\idiff y$, then expansion of $(d-y)^{-2}$ gives
\begin{equation}\label{cell:multipole-force}
 \int\frac{\omega_{\mathrm{ext}}(y)}{(d-y)^2}\idiff y
 =\frac{M_0}{d^2}+\frac{2M_1}{d^3}+\frac{3M_2}{d^4}+R_3,
 \qquad |R_3|\leq\frac C{|d|^5}\int|y|^3|\omega_{\mathrm{ext}}(y)|\idiff y.
\end{equation}
For the reflected Euler scalar $\widetilde W(y)=-W(\overline y)$, one has $M_2(\widetilde W)=-\int\overline y^{\,2}W(y)\idiff y$.  Its quadrupole therefore enters \eqref{cell:force-exact} with coefficient $-3e^{i\pi/6}/(4\pi)$ times the barred quadrupole.  In particular, the incoming phase is $\zeta$, not $\overline\zeta$.

\subsection{A first-exit statement with an early input certificate}
\label{cell:exit}

The net input must be determined before the receiver leaves its small neighborhood.  The following statement includes that timing requirement.

\begin{proposition}[The first growing exit]\label{cell:first-exit}
In Proposition~\ref{cell:response}, fix a certificate time $F\geq0$ and define
\begin{equation}\label{cell:certificate}
 B_{\mathrm{in}}=z_0+\int_0^Fe^{-\lambda s}f_{\mathrm{gr}}(s)\idiff s,
 \qquad T_{\mathrm{pred}}=\sigma^{-1}\log\frac{a_*}{|B_{\mathrm{in}}|}.
\end{equation}
Assume $0<|B_{\mathrm{in}}|<a_*$, put $H=T_{\mathrm{pred}}+2$, and impose \eqref{cell:response-assumptions} with $A=|B_{\mathrm{in}}|$.  In particular, require $a_*e^{2\sigma}\leq\eta_0$.  Define
\begin{equation}\label{cell:exit-errors}
 \varepsilon_{\mathrm{res}}=(1+H)\int_0^H\mathfrak a(s)\idiff s,
 \qquad
 \varepsilon_{\mathrm{tail}}=|B_{\mathrm{in}}|^{-1}
 \int_F^He^{-\sigma s}\mathfrak a(s)\idiff s.
\end{equation}
Choose $L_{\mathrm{in}}$ so that $C_{\mathrm{in}}e^{-\sigma L_{\mathrm{in}}}\leq1/4$, enlarging $C_{\mathrm{in}}$ so that $|z(s)|\leq C_{\mathrm{in}}Ae^{\sigma s}$.  Assume
\begin{equation}\label{cell:early-certificate}
 F\leq T_{\mathrm{pred}}-L_{\mathrm{in}},
\end{equation}
and, for a fixed sufficiently small $c_*>0$, assume
\begin{equation}\label{cell:exit-pointwise}
 \int_0^H\mathfrak a(s)\idiff s+
 \|\mathfrak a\|_{L^\infty((T_{\mathrm{pred}}-1,T_{\mathrm{pred}}+2)\cap(0,H))}
 \leq c_*a_*.
\end{equation}
If $\varepsilon_{\mathrm{res}}+\varepsilon_{\mathrm{tail}}$ is sufficiently small, the first time $T_{\mathrm{ev}}$ at which $|z|=a_*$ exists and satisfies
\begin{equation}\label{cell:exit-time-bound}
 |T_{\mathrm{ev}}-T_{\mathrm{pred}}|
 \leq C_{\mathrm{in}}(\varepsilon_{\mathrm{res}}+\varepsilon_{\mathrm{tail}}).
\end{equation}
The crossing is transverse in the sense that $\partial_s|z|\geq\sigma a_*/4$ almost everywhere while $a_*/2\leq|z|\leq2a_*$ near this event.  The pre-event history obeys
\begin{equation}\label{cell:whole-ancient-envelope}
 \mathfrak X(s)\leq C_{\mathrm{in}}a_*e^{\sigma(s-T_{\mathrm{ev}})},\qquad 0\leq s\leq T_{\mathrm{ev}}.
\end{equation}
For every fixed $L>0$ such that $(T_{\mathrm{ev}}-L,T_{\mathrm{ev}})\subset(F,H)$, the state on that interval is within
\begin{equation}\label{cell:donor-comparison-error}
 C_{\mathrm{in}}\left(\int_0^H\mathfrak a(s)\idiff s
 +a_*(\varepsilon_{\mathrm{res}}+\varepsilon_{\mathrm{tail}})\right)
\end{equation}
of $\mathcal G(a_*\zeta e^{\lambda(s-T_{\mathrm{ev}})})$ in a fixed smaller analytic norm, where $\zeta=z(T_{\mathrm{ev}})/a_*$.
\end{proposition}

\begin{proof}
For $F\leq s\leq H$, \eqref{cell:relative-response}, \eqref{cell:certificate}, and the tail bound give
\begin{equation}\label{cell:post-certificate-amplitude}
 \left|e^{-\lambda s}z(s)-B_{\mathrm{in}}\right|
 \leq C_{\mathrm{in}}|B_{\mathrm{in}}|
 (\varepsilon_{\mathrm{res}}+\varepsilon_{\mathrm{tail}}).
\end{equation}
Before $F$, \eqref{cell:response-bound} and \eqref{cell:early-certificate} give $|z|\leq a_*/4$.  After $F$, \eqref{cell:post-certificate-amplitude} places $|z(s)|$ between $(1-\varepsilon)|B_{\mathrm{in}}|e^{\sigma s}$ and $(1+\varepsilon)|B_{\mathrm{in}}|e^{\sigma s}$, with $\varepsilon=C_{\mathrm{in}}(\varepsilon_{\mathrm{res}}+\varepsilon_{\mathrm{tail}})$.  The intermediate value theorem gives a first crossing, and taking logarithms gives \eqref{cell:exit-time-bound}.

Near that crossing, \eqref{cell:complement-bound} and \eqref{cell:exit-pointwise}, inserted into \eqref{cell:amplitude-remainder}, yield $\partial_s|z|\geq\sigma a_*/4$ for sufficiently small $c_*$.  This proves transversality.  Equation \eqref{cell:whole-ancient-envelope} follows from \eqref{cell:response-bound} and \eqref{cell:exit-time-bound}.  Finally, apply \eqref{cell:post-certificate-amplitude} at $s$ and $T_{\mathrm{ev}}$ and divide the two formulas.  On a fixed negative-time window this gives
\begin{equation}\label{cell:phase-matched-amplitude}
 \left|z(s)-a_*\zeta e^{\lambda(s-T_{\mathrm{ev}})}\right|
 \leq C_{\mathrm{in}}a_*(\varepsilon_{\mathrm{res}}+\varepsilon_{\mathrm{tail}}).
\end{equation}
The local Lipschitz bound for $\mathcal G$ and \eqref{cell:complement-bound} prove \eqref{cell:donor-comparison-error}.
\end{proof}

\subsection{Identification with transported cells and positive-time comparison}
\label{cell:material}

We finish by relating the contour construction to the actual material labels used in a finite Euler array.  For an analytic solution of \eqref{cell:normal}, let $V=K*W+E$.  The normal equation states that $V(Z_j)-\partial_sZ_j$ is tangent to the interface.  Solving
\begin{equation}\label{cell:tangential-reparametrization}
 \vartheta'=
 \frac{\bigl(V(s,Z_j(s,\vartheta))-\partial_sZ_j(s,\vartheta)\bigr)
 \cdot\partial_\theta Z_j(s,\vartheta)}{|\partial_\theta Z_j(s,\vartheta)|^2}
\end{equation}
therefore gives $\partial_sZ_j(s,\vartheta(s))=V(s,Z_j(s,\vartheta(s)))$.  The velocity has the bounded-vorticity logarithmic Lipschitz modulus plus an integrable Lipschitz modulus from $E$.  The unused analytic radius gives the spatial regularity required in \eqref{cell:tangential-reparametrization}.  Osgood uniqueness identifies the parametrized boundaries with the ambient flow images.  Their bounded complementary components are thus transported as well.

For comparison with another bounded-vorticity cell, extend the exterior streamfunction from a slightly larger unused tube, retaining its spatially constant translation.  This gives a global divergence-free exterior with integrable Lipschitz seminorm, agreeing with the original field wherever the cells are supported.  Two cells with the same exterior and the same initial labels have maximal flow displacement $D(s)$ satisfying
\begin{equation}\label{cell:forced-flow-uniqueness}
 \mathrm D^+D\leq CD\bigl(1+|\log D|\bigr)+\|\mathrm D E(s)\|_\infty D,
 \qquad D(0)=0.
\end{equation}
For completeness, the near region uses bounded vorticity and integrability of $|y|^{-1}$; the annulus $4D<|y|<1$ contributes $CD\int_{4D}^1r^{-1}\idiff r$; and the exterior region uses the conserved $L^1$ norm.  Osgood's lemma yields $D=0$.  Apply this comparison while holding the actual exterior history of a tagged finite-array cell fixed.  The analytic solution remains strictly inside the tube by \eqref{cell:center-closure}, so a first loss of agreement through the tube boundary is impossible.

We shall also use a fixed positive-time comparison after an exit.  Fix $U_*>0$.  Suppose a cell has the same two values and areas as $G(a_*\zeta)$, its initial support lies in one fixed ball after centering at the outer centroid, and its translated initial discrepancy satisfies $\|W(0,\xi_{\mathrm{out}}(0)+\cdot)-G(a_*\zeta)\|_1\leq\varepsilon$.  Remove an exterior reference trajectory $\chi'=E(s,\chi)$ starting at its outer centroid.  If $\int_0^{U_*}\|\mathrm D E(s)\|_\infty\idiff s\to0$, then we have
\begin{equation}\label{cell:positive-comparison}
 \sup_{0\leq s\leq U_*}
 \|W(s,\chi(s)+\cdot)-W_*^\zeta(s)\|_1\longrightarrow0
 \quad\text{as }\varepsilon+\int_0^{U_*}\|\mathrm D E(s)\|_\infty\idiff s\longrightarrow0,
\end{equation}
uniformly for $|\zeta|=1$.  To prove this, the self-induced speed is uniformly bounded by $C\|W\|_1^{1/2}\|W\|_\infty^{1/2}$, and Gronwall keeps both supports in a fixed ball.  Decompose the initial vorticities into the common reference datum and an $L^1$ difference.  The transported difference has velocity at most $C\sqrt\varepsilon$.  The common part obeys the same near/far estimate as \eqref{cell:forced-flow-uniqueness}, giving
\begin{equation}\label{cell:positive-flow-comparison}
 \mathrm D^+D\leq CD(1+|\log D|)+C\sqrt\varepsilon+C_{U_*}\|\mathrm D E(s)\|_\infty.
\end{equation}
Osgood's estimate implies uniform convergence of the flow maps.  Their measure preservation, together with the uniform modulus of the reference inverse flow, converts this into \eqref{cell:positive-comparison} by approximation of the reference scalar in $L^1$ by continuous compactly supported functions.  Rotations of the fixed reference scalar form a compact subset of $L^1$, so the approximation is uniform in $\zeta$.  Since the values $a$, $-a/2$, and $0$ are distinct, the two domain indicators converge in $L^1$ as well.  Their centroids and centered quadrupoles consequently converge.  This comparison uses bounded-vorticity transport after the exit and does not require analytic contours at positive times of arbitrary length.

\section{Finite arrays and transfer between scales}\label{cas:section}\label{sec:finite-array}

The local cell estimates give a response to a prescribed exterior field. We now apply them to the exterior field generated by the other material cells of one finite Euler solution. The main point is to control this field before the receiving cell has left a small neighborhood of the radial shield. Sources at smaller spatial scales evolve more rapidly. Their leading quadrupoles provide the input, whereas the response of the larger cells must be estimated from the available past. We first establish material separation on a common time interval, then estimate the different sources of exterior strain, and finally order the local exit arguments by deterministic times.

\subsection{Parameters, finite data, and the conclusion}

We retain the unit shield $\Omega$, the unstable eigenvalue $\lambda=\sigma-\mathrm{i}a/4$, and the exact family $G$ from Proposition~\ref{cell:donor}. For the section size $a_*>0$, let $W_*(s)=G(a_*e^{\lambda s})$ when $s\leq0$, and continue this vorticity by its bounded-vorticity Euler evolution when $s>0$. Its barred quadratic moment is
\begin{equation}\label{cas:donor-moment}
 Q_*(s)=\int_{\R^2}\overline y^{\,2}W_*(s,y)\idiff y,
\end{equation}
where vectors in the plane are identified with complex numbers in moment formulas. Choose $r=1-\delta\in(1/2,1)$ as in Corollary~\ref{cell:transform}, so that the absolutely convergent integral
\begin{equation}\label{cas:transform}
 L=\int_{\R}e^{-r\lambda s}Q_*(s)\idiff s
\end{equation}
is nonzero. The order in which $a_*$ and $\delta$ are fixed will be specified below.

Let $(s_n)_{n\geq0}$ be a positive sequence satisfying
\begin{equation}\label{cas:sequence}
 1<\nu\leq\frac{s_{n+1}}{s_n}\leq V<\infty,
 \qquad \sum_{n\geq0}r^ns_n<\infty.
\end{equation}
Define $m_0=\log(1/r)/\log\nu$, choose $M>m_0+1$, and, for a final large parameter $\Lambda$, set
\begin{equation}\label{cas:geometry}
 \begin{gathered}
 S_n=\Lambda s_n,\qquad R_n=e^{-(S_n-S_0)},\qquad
 h_n=S_n^{-M},\qquad d_n=h_nR_n,\\
 A_n=r^{-n},\qquad \tau_n=A_n^{-1}=r^n,\qquad
 c_n=R_ne_1,\qquad q_n=R_{n+1}/R_n.
 \end{gathered}
\end{equation}
Thus $R_n$ locates a cell, $d_n$ is its radius, and $A_n$ is its vorticity amplitude. The time variable in the $n$th unit cell is $s=A_nt$.

For later use, define the geometric factor, transfer coefficient, and clocks by
\begin{equation}\label{cas:clocks}
 \begin{gathered}
 \mathcal B(q)=(1-q)^{-4}+(1+q)^{-4},\qquad
 K_2=\frac{3e^{\mathrm{i}\pi/6}}{4\pi},\\
 D_n=-K_2\mathcal B(q_n)h_{n+1}^4q_n^4L,\qquad
 \ell_n=\frac1\sigma\log\frac{a_*}{|D_n|},\qquad
 \kappa_n=\sum_{j\geq0}r^j\ell_{n+j},\qquad
 B_{0,n}=a_*e^{-\sigma\kappa_n}.
 \end{gathered}
\end{equation}
For sufficiently large $\Lambda$, every $\ell_n$ is positive and the series defining $\kappa_n$ converges. The notation $B_{0,n}$ denotes an amplitude, not a spatial ball.

Put $O=\operatorname{diag}(1,-1)$. At a finite cutoff $N$, let
\begin{equation}\label{cas:finite-data}
 \begin{gathered}
 W_{n,0}^N=\Omega\quad(n<N),\qquad
 W_{N,0}^N=G\bigl(a_*e^{-\lambda\kappa_N}\bigr),\\
 \omega_{n,+}^N(0,x)=(-1)^nA_n
 W_{n,0}^N\bigl(d_n^{-1}O^n(x-c_n)\bigr),\qquad
 \omega_{n,-}^N(0,x)=\omega_{n,+}^N(0,-x),\\
 \omega^N(0,x)=\sum_{n=0}^N
 \bigl(\omega_{n,+}^N(0,x)+\omega_{n,-}^N(0,x)\bigr).
 \end{gathered}
\end{equation}
The canonical bounded-vorticity solution from this datum is denoted by $\omega^N$. Each summand is transported by its flow and retains its label. The reflection and the factor $(-1)^n$ occur together because vorticity changes sign under an orientation-reversing Euler symmetry. Consequently every cell has the same forward unit dynamics. Let $c_n^{\mathrm{out},N}(t)$ denote the centroid of the transported positive outer domain. Its centered unit vorticity and barred quadratic moment are
\begin{equation}\label{cas:centered-unit-cell}
 \begin{gathered}
 W_n^N(s,y)=(-1)^nA_n^{-1}\omega_{n,+}^N
 \bigl(\tau_ns,c_n^{\mathrm{out},N}(\tau_ns)+d_nO^ny\bigr),\\
 Q_n^N(s)=\int_{\R^2}\overline y^{\,2}W_n^N(s,y)\idiff y.
 \end{gathered}
\end{equation}
The contour size $\mathfrak X_n^N$ below is the spectrally adapted norm from \eqref{cell:spectral-norm}, evaluated at the decreasing radius supplied by Proposition~\ref{cell:response}.

\begin{theorem}[Uniform finite arrays]\label{cas:finite}
There are choices of $a_*>0$ and $r\in(0,1)$ for which the following statement holds. Fix a sequence satisfying \eqref{cas:sequence}, and let the geometric quantities be given by \eqref{cas:geometry}--\eqref{cas:clocks}. Fix $B=1$, a clock tolerance $0<\varepsilon_t<\delta B/10$, and a moment tolerance $\varepsilon_Q>0$. There are finite windows $T_-,T_+>0$ and a sufficiently large $\Lambda$ such that, for every $N$, the solution with datum \eqref{cas:finite-data} has the following properties on
\begin{equation}\label{cas:common-time}
 0\leq t\leq T:=\kappa_0+B+T_++3.
\end{equation}
Every label remains in $B(\pm c_n,R_n/8)$. Every positive cell has a first transverse crossing of the local section $|z|=a_*$ at a time $t_n^N=\tau_n\widetilde\kappa_n^N$, and the corresponding unit times satisfy
\begin{equation}\label{cas:event-recursion}
 \begin{gathered}
 \widetilde\kappa_n^N=\ell_n+r\widetilde\kappa_{n+1}^N+e_n^N
       \quad(n<N),\qquad |e_n^N|\leq\varepsilon_t,\\
 |\widetilde\kappa_N^N-\kappa_N|\leq\varepsilon_t,
 \qquad
 |\widetilde\kappa_n^N-\kappa_n|\leq\varepsilon_t/\delta.
 \end{gathered}
\end{equation}
Before that crossing, the centered analytic contour norm satisfies
\begin{equation}\label{cas:ancient-envelope}
 \mathfrak X_n^N(s)\leq Ca_*e^{\sigma(s-\widetilde\kappa_n^N)}
       \qquad(0\leq s\leq\widetilde\kappa_n^N).
\end{equation}
Let $Q_n^N(s)$ be the barred quadratic moment in the positive cell's unit orientation, centered at its actual outer centroid. For some $\zeta_n^N\in\mathbb C$ with $|\zeta_n^N|=1$, it satisfies
\begin{equation}\label{cas:window}
 \sup_{-T_-\leq v\leq T_+}
 \bigl|Q_n^N(\widetilde\kappa_n^N+v)-\zeta_n^NQ_*(v)\bigr|
 \leq\varepsilon_Q.
\end{equation}
The support radii, impulses, and moments after the crossing satisfy the polynomial bounds established in \eqref{cas:post-support} and \eqref{cas:post-impulse} below. All constants are independent of $n$ and $N$.

Finally, there are four fixed real functions in $C_c^\infty(\R^2)$ such that one of them, denoted by $\varphi$, satisfies along an infinite subsequence
\begin{equation}\label{cas:witness}
 \int_{\R^2}\varphi(x)\omega^N(t_0^N,x)\idiff x
 \geq c d_0^4a_*>0.
\end{equation}
Its integral against the complete radial shield datum is zero.
\end{theorem}

\subsection{Choice of constants and elementary scale estimates}

\begin{proof}[Proof of Theorem~\ref{cas:finite}]
We first choose the local constants without referring to an output bound that will depend on them. The donor expansion from Proposition~\ref{cell:donor} gives
\begin{equation}\label{cas:pole}
 Q_*(s)=(1-\mathrm{i}\sqrt3)a_*e^{\lambda s}
             +O(a_*^2e^{2\sigma s})\quad(s\leq0),\qquad
 L=\frac{a_*}{\sigma\delta}+O(1+a_*^2).
\end{equation}
The positive-time error in the last formula is bounded uniformly in $a_*$ because the donor has support radius at most $C(1+s)$. Let $C_{\rm an}\geq1$ be the fixed multipole-to-analytic-norm constant in \eqref{cas:analytic-source} below. It depends only on the unit cell and the fixed norm radii. Corollary~\ref{cell:transform} gives a constant $K_{\rm iso}=5C_{\rm an}$ such that
\begin{equation}\label{cas:absolute-transform}
 |L|\geq\frac{a_*}{2\sigma\delta},\qquad
 C_{\rm an}\int_{\R}e^{-r\sigma s}|Q_*(s)|\idiff s
 \leq K_{\rm iso}|L|.
\end{equation}
Thus $K_{\rm iso}$ can be fixed before $a_*$ and $\delta$, by requiring $\delta\leq c_{\rm tr}a_*$ afterwards. Only the isolated quadrupole enters this constant. The other moments of an actual cell will occur with small geometric factors, and their bounds will be absorbed by the final choice of $\Lambda$.

Choose a local absolute-input budget
\begin{equation}\label{cas:input-budget}
 K_{\rm in}\geq32e^{2\sigma B}(K_{\rm iso}+1),
\end{equation}
with a larger numerical factor if required by the fixed projections. Let $C_{\rm in}$ and $\eta_0(K_{\rm in})$ be the constants in Proposition~\ref{cell:response}. Choose $a_*$ sufficiently small that
\begin{equation}\label{cas:section-choice}
 2a_*e^{2\sigma(B+1)}\leq\eta_0(K_{\rm in}),
\end{equation}
and then choose $\delta$ sufficiently small for \eqref{cas:absolute-transform}. Fix a contour-cone constant $C_{\rm cone}\geq4C_{\rm in}$, enlarged for the initial manifold chart. The negative window $T_-$ will be chosen large enough for the early-crossing condition of Proposition~\ref{cell:first-exit}, for the ancient transform tail, and for \eqref{cas:uncertified-input} below. Set $W=B+2T_-$. Choose $T_+$ large enough for the future transform tail. The sequence, the exponent $M$, the windows, and all requested errors are fixed before increasing $\Lambda$.

The logarithm in \eqref{cas:clocks} has the explicit form
\begin{equation}\label{cas:ell-size}
 \sigma\ell_n=4(S_{n+1}-S_n)+4M\log S_{n+1}+\gamma_n,
 \qquad |\gamma_n|\leq C.
\end{equation}
Indeed, the remaining factors in $D_n$ are fixed, and $\mathcal B(q_n)$ stays in a compact subset of $(0,\infty)$. Since $S_{n+1}\geq\nu S_n$, formula \eqref{cas:ell-size} gives $\ell_n\geq cS_n$ for large $\Lambda$. It also gives
\begin{equation}\label{cas:clock-size}
 \kappa_0\leq C\Lambda,\qquad
 t_n:=r^n\kappa_n=\sum_{k\geq n}r^k\ell_k\downarrow0,\qquad
 \kappa_n\leq A_n\kappa_0\leq C\Lambda A_n.
\end{equation}
To obtain the first inequality, use \eqref{cas:sequence} for the differences $s_{n+1}-s_n$ and the bounds $\log s_n\leq\log s_0+n\log V$. The second follows by multiplying the definition of $\kappa_n$ by $r^n$. We will use the exact identities
\begin{equation}\label{cas:clock-identities}
 \kappa_n=\ell_n+r\kappa_{n+1},\qquad
 B_{0,n}=|D_n|e^{-\sigma r\kappa_{n+1}}.
\end{equation}

Because $s_n\geq s_0\nu^n$, we have $A_n\leq(s_n/s_0)^{m_0}$. Consequently, we obtain
\begin{equation}\label{cas:mass-small}
 \sup_n A_nh_n\leq C\Lambda^{-M},\qquad
 B_{\rm mass}:=\sum_{n\geq0}A_nh_n^2\leq C\Lambda^{-2M}.
\end{equation}
Combining this with \eqref{cas:clock-size} yields
\begin{equation}\label{cas:waiting-small}
 \begin{gathered}
 \sup_n\kappa_nh_n\leq C\Lambda^{1-M},\qquad
 \sup_n\kappa_n^2h_n^2\leq C\Lambda^{2-2M},\\
 \sup_n\kappa_nh_n^2\leq C\Lambda^{1-2M},\qquad
 \sup_n\kappa_nA_n^{-1}B_{\rm mass}\leq C\Lambda^{1-2M}.
 \end{gathered}
\end{equation}
Moreover, $\kappa_n$ is bounded by a fixed power of $S_n$. Thus every fixed polynomial in $S_n$ and $\kappa_n$, multiplied by $e^{-cS_n}$, tends to zero uniformly in $n$ as $\Lambda\to\infty$.

\subsection{Material separation and bounds valid after a crossing}

The finite Euler solution preserves the two values, areas, and circulation of every label. A near-far split of the Biot--Savart integral therefore gives
\begin{equation}\label{cas:label-bounds}
 \|\omega_{j,\pm}^N(t)\|_1\leq CA_jd_j^2,\qquad
 \|\omega_{j,\pm}^N(t)\|_\infty\leq CA_j,\qquad
 \|K*\omega_{j,\pm}^N(t)\|_\infty\leq CA_jd_j.
\end{equation}
For the last inequality, split at radius $d_j$: the near integral is at most $C\|\omega\|_\infty d_j$, and the far integral is at most $C\|\omega\|_1/d_j$.

Bootstrap all supports inside $B(\pm c_j,R_j/8)$. The datum is even, so uniqueness of the finite bounded-vorticity evolution makes its velocity odd and preserves the pairing of each label. At a point $x$ in the positive $n$th tube, every coarser paired source vanishes at the origin and is separated from the segment joining $0$ to $x$. Its gradient is bounded by $CA_jd_j^2/R_j^2=CA_jh_j^2$. Hence, we obtain
\begin{equation}\label{cas:coarse-velocity}
 |u_{<n}(t,x)|\leq C|x|\sum_{j<n}A_jh_j^2.
\end{equation}
Use \eqref{cas:label-bounds} for the cell itself, its partner, and all finer labels. This gives
\begin{equation}\label{cas:tube-velocity}
 \frac{|u^N(t,x)|}{R_n}
 \leq C\left(A_nh_n+B_{\rm mass}
       +\sum_{j>n}A_jh_j\frac{R_j}{R_n}\right)
 \leq C\Lambda^{-M}.
\end{equation}
The final series is bounded by its first term times a uniform constant because $R_{j+1}/R_j\leq e^{-(\nu-1)\Lambda s_j}$. Since $T\leq C\Lambda$, the displacement is at most $C\Lambda^{1-M}R_n$. Together with the initial radius $Cd_n$, this is strictly smaller than the bootstrap margin for sufficiently large $\Lambda$. A first-exit argument proves confinement on \eqref{cas:common-time}. It also proves the sharper displacement bound just obtained for every material point and every domain centroid.

Remove the positive $n$th cell from the velocity and denote the remaining field by $u_{\rm ext,n}^N$. Its reference trajectory starts at $c_n$ and satisfies
\begin{equation}\label{cas:reference-center}
 \frac{\idiff}{\idiff t}c_n^{\rm ref}(t)
   =u_{\rm ext,n}^N(t,c_n^{\rm ref}(t)),\qquad
 c_n^{\rm ref}(0)=c_n.
\end{equation}
The preceding estimates keep it inside the same tube. The negative reference trajectory and negative outer centroid are the opposites of their positive partners.

Normalize the exterior to unit length $d_n$ and unit time $A_n^{-1}$. Its analytic coefficient norm, denoted by $\mathfrak a_{\rm raw,n}$, and its second spatial derivative satisfy on the whole physical tube
\begin{equation}\label{cas:raw-exterior}
 \mathfrak a_{\rm raw,n}\leq C(A_n^{-1}B_{\rm mass}+h_n^2),
 \qquad \|D_y^2E_n\|_\infty\leq Ch_n^3.
\end{equation}
For the first estimate, sum $A_jd_j^2/\operatorname{dist}^2$ and divide by $A_n$. For the second, sum $A_jd_j^2/\operatorname{dist}^3$ and multiply by $d_n/A_n$. The own partner contributes $Ch_n^3$. A coarser term with $j=n-k$ divided by $h_n^3$ is at most
\begin{equation}\label{cas:curvature-series}
 r^k\left(\frac{h_{n-k}}{h_n}\right)^2\frac{R_n}{R_{n-k}}
 \leq\left(rV^{2M}e^{-(\nu-1)\Lambda s_0}\right)^k,
\end{equation}
whose sum is small. The finer sum is bounded by the same method. Harmonic Cauchy estimates on a slightly larger tube give the first inequality of \eqref{cas:raw-exterior} also for the full analytic coefficient norm on every radius-four unit ball under consideration.

The accumulated raw strain on the common physical interval satisfies
\begin{equation}\label{cas:raw-integral}
 A_nT\mathfrak a_{\rm raw,n}
 \leq CT(B_{\rm mass}+A_nh_n^2)
 \leq C\Lambda^{1-2M}.
\end{equation}
After a crossing, restart the reference trajectory at the actual outer centroid, and let $v\geq0$ denote unit time since that crossing. The self-induced speed is bounded by a fixed constant by \eqref{cas:label-bounds}, whereas the relative exterior velocity is bounded by $\mathfrak a_{\rm raw,n}$ times distance. Gronwall's inequality and \eqref{cas:raw-integral} therefore give the support estimate
\begin{equation}\label{cas:post-support}
 \rho_n(v)\leq C(1+v).
\end{equation}
These unit balls fit inside the physical tube because $h_n(1+A_nT)\leq C\Lambda^{1-M}$.

Let $I_n$ be the signed complex impulse of the unit cell. Its value is independent of the center because the circulation is zero. Before the crossing, the exact harmonic cancellation in the local impulse identity gives
\begin{equation}\label{cas:ancient-impulse-equation}
 I_n'=\mathrm D E_n(\xi_n)I_n+O(h_n^3\mathfrak X_n),\qquad I_n(0)=0,
\end{equation}
where $\xi_n$ denotes the actual outer centroid in the chosen unit coordinates. If the contour has an exponential envelope relative to an event time, variation of constants in \eqref{cas:ancient-impulse-equation} gives the same envelope for $I_n$, multiplied by $Ch_n^3$. In particular, after a crossing has been certified, its impulse satisfies
\begin{equation}\label{cas:ancient-impulse}
 |I_n(v)|\leq Ch_n^3a_*e^{\sigma v}\qquad(v\leq0).
\end{equation}
For $v\geq0$, Taylor expansion of the exterior about the reference center and \eqref{cas:post-support} give
\begin{equation}\label{cas:post-impulse}
 |I_n'|\leq C\mathfrak a_{\rm raw,n}|I_n|+Ch_n^3(1+v)^2,
 \qquad |I_n(v)|\leq Ch_n^3(1+v)^3.
\end{equation}
Centered absolute quadratic and cubic moments are bounded by $C(1+v)^2$ and $C(1+v)^3$. These estimates do not require contour regularity after the crossing.

\subsection{The adjacent source and the transfer coefficient}

The leading source for cell $n$ is the paired level $n+1$. Work in the parent's oriented physical coordinates. Write $c,-c$ for the current child outer centroids, $I_c$ for the positive child's complex impulse, and $M_c$ for its centered unbarred quadratic moment. The negative child has impulse $-I_c$ and quadratic moment $M_c$. Expanding each kernel about its own centroid gives
\begin{equation}\label{cas:paired-expansion}
 \begin{aligned}
 \int\frac{\omega_{\rm pair}(y)}{(x-y)^2}\idiff y
 ={}&2I_c\bigl((x-c)^{-3}-(x+c)^{-3}\bigr)
 +3M_c\bigl((x-c)^{-4}+(x+c)^{-4}\bigr)+\mathcal R_3.
 \end{aligned}
\end{equation}
Since $|x|\asymp R_n$ and $|c|\asymp R_{n+1}$, the difference multiplying $I_c$ is bounded by $CR_{n+1}/R_n^4$. In physical units $I_c=A_{n+1}d_{n+1}^3I_{n+1}$. Thus \eqref{cas:ancient-impulse} and \eqref{cas:post-impulse} show that the impulse contribution, relative to $A_{n+1}d_{n+1}^4/R_n^4$, costs only $Ch_{n+1}^2$ times its age envelope. Estimating the two impulses separately would give an incorrect loss of a negative power of $q_n$.

The cubic remainder in \eqref{cas:paired-expansion} is bounded by the absolute cubic moment divided by $R_n^5$. Relative to the same quadratic scale, its coefficient is at most $Ch_{n+1}q_n$. Before a crossing, subtract the radial shield centered at the actual outer centroid. Its exterior field is exactly zero, and the absolute moments of the residual have the contour envelope. This subtraction is made at each time and is never differentiated. The displacement estimate following \eqref{cas:tube-velocity} changes the quadratic kernel factors by a relative $O(\Lambda^{1-M})$. Differentiating the exterior kernel $l-1$ additional times gives, in the parent's unit Taylor coefficient of order $l\geq1$, an extra factor at most $C^l(1+l)^3h_n^{l-1}$. Its weighted sum is uniformly finite. Consequently the expansion controls the full analytic exterior norm used in Proposition~\ref{cell:response}.

In the parent's orientation the child is a negative reflected unit cell. Therefore its unbarred physical moment is
\begin{equation}\label{cas:reflection-moment}
 M_c=-A_{n+1}d_{n+1}^4Q_{n+1}^N.
\end{equation}
Combining the local growing projection with \eqref{cas:paired-expansion} gives, in parent unit time,
\begin{equation}\label{cas:adjacent-forcing}
 f_{{\rm gr,adj},n}(s)
 =-K_2\mathcal B(q_n)r^{-1}h_{n+1}^4q_n^4
 \widehat Q_{n+1}^N\bigl(s/r-\widetilde\kappa_{n+1}^N\bigr)
 +\mathcal R_{{\rm adj},n}(s),
\end{equation}
where $\widehat Q_{n+1}^N(v)=Q_{n+1}^N(\widetilde\kappa_{n+1}^N+v)$ and the error has the bounds just proved. The change of variable $s=r(v+\widetilde\kappa_{n+1}^N)$ cancels the factor $r^{-1}$. For the isolated donor, the resulting growing coefficient is exactly
\begin{equation}\label{cas:adjacent-transform}
 D_n\zeta_{n+1}^Ne^{-\lambda r\widetilde\kappa_{n+1}^N}.
\end{equation}
Here the donor phase is $\zeta_{n+1}^N$, rather than its complex conjugate. Formula \eqref{cas:reflection-moment} is the reason for this choice of alternating orientations.

We make the analytic-norm estimate explicit, including the part of a child trajectory after it leaves the analytic chart. For a certified child, choose constants $C_{\rm cert}$ and $C_{\rm post}$ so that the radial-subtracted absolute cubic moment before exit, the whole absolute cubic moment after exit, and the rescaled impulse $h_{n+1}^{-3}|I_{n+1}|$ are bounded by the age majorant
\begin{equation}\label{cas:age-majorant}
 \mathcal P(v)=
 \begin{cases}
 C_{\rm cert}a_*e^{\sigma v},&v\leq0,\\
 C_{\rm post}(1+v)^3,&v>0.
 \end{cases}
\end{equation}
Here $C_{\rm cert}$ may depend on $K_{\rm in}$; it is fixed before $T_-$. The constant $C_{\rm post}$ follows from \eqref{cas:post-support}--\eqref{cas:post-impulse}. Both are independent of $n,N,\Lambda$. Define the small geometric error by
\begin{equation}\label{cas:geometric-error}
 \varepsilon_{{\rm geom},n}
 =C_{\rm geom}\bigl(\Lambda^{1-M}+h_{n+1}^2+h_{n+1}q_n\bigr).
\end{equation}
Here $C_{\rm geom}$ is fixed after the sequence and $M$; it absorbs the constant in the material displacement estimate, and is independent of $n,N,\Lambda$. In particular, it is not included in the earlier constant $C_{\rm an}$. For $v=s/r-\widetilde\kappa_{n+1}^N$ in the available history, the full analytic coefficient norm and the scalar projection error satisfy
\begin{equation}\label{cas:analytic-source}
 \begin{aligned}
 \mathfrak a_{{\rm adj},n}(s)
 &\leq C_{\rm an}|K_2\mathcal B(q_n)|r^{-1}h_{n+1}^4q_n^4
 \bigl(|\widehat Q_{n+1}^N(v)|+
             \varepsilon_{{\rm geom},n}\mathcal P(v)\bigr),\\
 |\mathcal R_{{\rm adj},n}(s)|
 &\leq C_{\rm an}|K_2\mathcal B(q_n)|r^{-1}h_{n+1}^4q_n^4
       \varepsilon_{{\rm geom},n}\mathcal P(v).
 \end{aligned}
\end{equation}
Indeed, the paired impulse term is bounded by $h_{n+1}^2\mathcal P(v)$, the cubic remainder by $h_{n+1}q_n\mathcal P(v)$, and the change of the leading kernel by $\Lambda^{1-M}\mathcal P(v)$. The Taylor coefficients of order $l\geq1$ have the additional factors $C^l(1+l)^3h_n^{l-1}$ already derived from \eqref{cas:paired-expansion}. Their weighted sum is bounded by an absolute constant for all sufficiently large $\Lambda$. This proves \eqref{cas:analytic-source} without using a future contour norm.

The error majorant is integrable with the donor weight because its weighted integral satisfies
\begin{equation}\label{cas:age-integral}
 \begin{aligned}
 \frac1{|L|}\int_{\R}e^{-r\sigma v}\mathcal P(v)\idiff v
 &=\frac{C_{\rm cert}a_*}{\sigma\delta|L|}
   +\frac{C_{\rm post}}{|L|}
       \int_0^\infty e^{-r\sigma v}(1+v)^3\idiff v\\
 &\leq C(C_{\rm cert}+1).
 \end{aligned}
\end{equation}
The last inequality uses $r\geq1/2$, \eqref{cas:absolute-transform}, and $\delta\leq c_{\rm tr}a_*$. Thus the weighted contribution of the second term in \eqref{cas:analytic-source}, divided by $B_{0,n}$, is $o(1)$ for fixed $C_{\rm cert}$. It does not enlarge $K_{\rm iso}$ or $K_{\rm in}$.

\subsection{A causal estimate for the coarser cells}

We next estimate the coarser response using only a stopped history. Fix a receiver $n$ and a time before its deterministic certificate. For $k<n$, suppose the available part of the history satisfies the cone bounds
\begin{equation}\label{cas:stopped-cones}
 \mathfrak X_k(t)\leq C_{\rm cone}B_{0,k}e^{\sigma A_kt}.
\end{equation}
Physical time is used in this subsection. If $t\leq\tau_n\kappa_n$, the clock identities give
\begin{equation}\label{cas:coarse-smallness}
 \mathfrak X_k(t)\leq C_{\rm cone}a_*e^{-\sigma(\kappa_k-r^{n-k}\kappa_n)}
 \leq C_{\rm cone}a_*e^{-\sigma\ell_k}.
\end{equation}
This is uniformly smaller than any fixed multiple of $\delta$ for large $\Lambda$.

For each such coarser cell subtract its radial shield at the actual outer centroid. The residual has $L^1$ norm at most $CA_kd_k^2\mathfrak X_k$. The pair of residuals has zero circulation and zero impulse, and its absolute second moment about the origin is at most $CA_kh_k^2R_k^4\mathfrak X_k$. Thus the physical strain coupling from paired level $j$ to level $k$ is bounded by
\begin{equation}\label{cas:couplings}
 \beta_{kj}\leq CA_jh_j^2
 \left(\frac{R_j}{\max(R_k,R_j)}\right)^4.
\end{equation}
The diagonal term represents the cell's own partner. For a finer source, \eqref{cas:couplings} uses its paired cancellations. For a coarser source, the direct residual $L^1$ estimate suffices.

Let $Q_f(t)$ be the sum of the absolute second moments of the paired levels $j\geq n$, using radial-subtracted residuals before a crossing and the whole paired labels afterwards. The contour norm inequalities of the local construction imply
\begin{equation}\label{cas:head-system}
 \mathrm D^+\mathfrak X_k\leq\sigma A_k\mathfrak X_k+CA_k\mathfrak X_k^2
          +\sum_{j<n}\beta_{kj}\mathfrak X_j+C Q_f/R_k^4\qquad(k<n).
\end{equation}
Define $\mathfrak H_n(t)=\max_{k<n}R_k^4\mathfrak X_k(t)$, with $\mathfrak H_n=0$ when $n=0$. Multiplying \eqref{cas:couplings} by $R_k^4/R_j^4$ gives a row sum at most $CB_{\rm mass}$: if $j\geq k$, the factors of $R_k/R_j$ cancel, and if $j<k$, the remaining factor is at most one. Since $A_k\leq rA_n$, equations \eqref{cas:coarse-smallness} and \eqref{cas:mass-small} yield
\begin{equation}\label{cas:causal-ode}
 \mathrm D^+\mathfrak H_n\leq(\sigma A_n-\sigma_{\rm gap}A_n)\mathfrak H_n+CQ_f,
 \qquad \mathfrak H_n(0)=0,\qquad \sigma_{\rm gap}=\sigma\delta/2.
\end{equation}
Variation of constants gives, on every available prefix,
\begin{equation}\label{cas:causal-pointwise}
 \mathfrak H_n(t)\leq Ce^{\sigma A_nt}\int_0^te^{-\sigma A_nv}Q_f(v)\idiff v.
\end{equation}
Keeping the factor $e^{-\sigma_{\rm gap}A_n(t-v)}$ before integrating once more gives
\begin{equation}\label{cas:causal-integrated}
 \int_0^{H/A_n}e^{-\sigma A_nt}\mathfrak H_n(t)\idiff t
 \leq\frac{C}{\sigma_{\rm gap}A_n}
       \int_0^{H/A_n}e^{-\sigma A_nt}Q_f(t)\idiff t.
\end{equation}
In particular, neither estimate uses a future bound for $Q_f$.

To keep the time normalization explicit, put
\begin{equation}\label{cas:head-parameters}
 \begin{gathered}
 J_Q(H)=\int_0^{H/A_n}e^{-\sigma A_nt}Q_f(t)\idiff t,
 \qquad J_f(H)=J_Q(H)/R_n^4,\\
 Z_n=\sum_{k<n}\frac{A_kh_k^2}{R_k^4},\qquad
 \mu_n=\frac{CR_n^4Z_n}{\sigma_{\rm gap}A_n}.
 \end{gathered}
\end{equation}
The physical coarser strain is at most $Z_n\mathfrak H_n$, so its unit value is at most $A_n^{-1}Z_n\mathfrak H_n$. The identity $\idiff s=A_n\idiff t$, together with \eqref{cas:causal-pointwise}--\eqref{cas:causal-integrated}, gives
\begin{equation}\label{cas:head-input}
 \begin{aligned}
 J_{\rm head}(H)&:=\int_0^He^{-\sigma s}\mathfrak a_{\rm head}(s)\idiff s
                \leq\mu_nJ_f(H),\\
 \mathfrak a_{\rm head}(s)&\leq C\sigma_{\rm gap}\mu_ne^{\sigma s}J_f(H)
                \qquad(0\leq s\leq H).
 \end{aligned}
\end{equation}
The rapidly separated radii imply
\begin{equation}\label{cas:head-small-factor}
 \mu_n\leq Crh_0^2e^{-4(S_n-S_{n-1})}
 \leq C\Lambda^{-2M}e^{-cS_n}\quad(n\geq1),\qquad \mu_0=0.
\end{equation}
For example, the contribution of $k<n$ to $R_n^4Z_n/A_n$ is $r^{n-k}h_k^2e^{-4(S_n-S_k)}$, and the sum is bounded by a constant times its last term.

\subsection{The remaining weighted inputs}

We now bound $J_f$ and the direct finer inputs, still on available prefixes. The material tubes give a useful bound valid at every positive age: the absolute second moment of a departed paired level $j$ is at most
\begin{equation}\label{cas:whole-pair-moment}
 C A_jd_j^2R_j^2=C A_jh_j^2R_j^4.
\end{equation}
Before departure the same bound is multiplied by its contour envelope.

For a depth $m\geq1$, split the weighted integral at the nominal child event, allowing a fixed clock-box displacement. Integrating the ancient exponential costs a factor at most $(1-r^m)^{-1}\leq\delta^{-1}$, and integrating the whole-label bound \eqref{cas:whole-pair-moment} costs at most $Cr^{-m}$. The direct input therefore satisfies
\begin{equation}\label{cas:depth-input}
 \frac{J_{n,m}}{B_{0,n}}
 \leq C(1+C_{\rm cone})r^{-m}h_{n+m}^2
 \exp\bigl(-4(S_{n+m}-S_n)+\sigma\kappa_n-\sigma r^m\kappa_{n+m}\bigr).
\end{equation}
For an uncertified child, extend its available cone up to the nominal event and then use \eqref{cas:whole-pair-moment}; this only enlarges the estimate. By \eqref{cas:ell-size} and \eqref{cas:clock-identities}, the leading exponent for $m\geq2$ is
\begin{equation}\label{cas:depth-exponent}
 \begin{aligned}
 &4\sum_{j=0}^{m-1}r^j(S_{n+j+1}-S_{n+j})-4(S_{n+m}-S_n)\\
 &\hspace{1cm}=-4\sum_{j=1}^{m-1}(1-r^j)(S_{n+j+1}-S_{n+j})
 \leq-4\delta(S_{n+m}-S_{n+1}).
 \end{aligned}
\end{equation}
The logarithmic correction is bounded by
\begin{equation}\label{cas:depth-log}
 4M\sum_{j=0}^{m-1}r^j\log S_{n+j+1}
 \leq\frac{4M}{\delta}\log S_n+\frac{4M\log V}{\delta^2}.
\end{equation}
Also $r^{-m}h_{n+m}^2\leq h_n^2(r^{-1}\nu^{-2M})^m$, and the ratio in parentheses is less than one. Hence, we obtain
\begin{equation}\label{cas:deep-small}
 \sum_{m\geq2}\frac{J_{n,m}}{B_{0,n}}
 \leq C(1+C_{\rm cone})S_n^Ce^{-cS_n}
 =o(1)(1+C_{\rm cone}),
\end{equation}
uniformly in $n$ and $N$. All occurrences of $o(1)$ in this proof refer to $\Lambda\to\infty$ with the previously selected constants fixed.

At depth one, the clock identity gives the explicit cancellation
\begin{equation}\label{cas:adjacent-coarse-loss}
 h_{n+1}^2q_n^4e^{\sigma\kappa_n-\sigma r\kappa_{n+1}}
 =\frac{a_*}{|K_2\mathcal B(q_n)L|}\,h_{n+1}^{-2}
 \leq C S_n^{2M}.
\end{equation}
The receiver's own paired residual contributes at most $CC_{\rm cone}h_n^2HB_{0,n}$ to $J_f(H)$. Equations \eqref{cas:deep-small}--\eqref{cas:adjacent-coarse-loss} thus give
\begin{equation}\label{cas:Qf-bound}
 \frac{J_f(H)}{B_{0,n}}
 \leq C(1+C_{\rm cone})
       \bigl(S_n^{2M}+h_n^2H+S_n^Ce^{-cS_n}\bigr)=:P_n.
\end{equation}
For the terminal receiver, the adjacent and deeper contributions are absent. The right side is bounded by a fixed polynomial in $S_n$ on the intervals under consideration. Substitution into \eqref{cas:head-input} yields
\begin{equation}\label{cas:head-final}
 J_{\rm head}(H)\leq\mu_nP_nB_{0,n},\qquad
 \mathfrak a_{\rm head}(s)\leq C\mu_nP_na_*e^{\sigma(s-\kappa_n)}.
\end{equation}

The own partner produces unit strain at most $Ch_n^2\mathfrak X_n$. On a stopped cone, its weighted input is consequently bounded by
\begin{equation}\label{cas:partner-input}
 J_{\rm partner}\leq CC_{\rm cone}h_n^2\kappa_n B_{0,n}
                  =o(1)B_{0,n}.
\end{equation}

For the adjacent child, the paired expansion gives a sharper estimate than \eqref{cas:depth-input}. Before that child is certified, its unit strain at the receiver is at most
\begin{equation}\label{cas:uncertified-strain}
 \mathfrak a_{{\rm adj},n}(s)
 \leq Cr^{-1}h_{n+1}^4q_n^4C_{\rm cone}B_{0,n+1}e^{\sigma s/r}.
\end{equation}
Integrating and using \eqref{cas:absolute-transform} and \eqref{cas:clock-identities} gives
\begin{equation}\label{cas:uncertified-input}
 B_{0,n}^{-1}\int_0^se^{-\sigma v}\mathfrak a_{{\rm adj},n}(v)\idiff v
 \leq CC_{\rm cone}e^{\sigma\delta(s/r-\kappa_{n+1})}.
\end{equation}
At an available time before child certification, $s/r\leq\kappa_{n+1}-W$. Increasing $T_-$, and hence $W$, makes the coefficient of $C_{\rm cone}$ in \eqref{cas:uncertified-input} as small as prescribed.

For a certified child, split every prefix into the ancient part before $-T_-$, the traversed part of its comparison window $[-T_-,T_+]$, and the future part after $T_+$. The actual quadrupole before $-T_-$ is bounded by $C_{\rm cert}a_*e^{\sigma v}$, so its weighted tail is at most $C_{\rm cert}a_*e^{-\sigma\delta T_-}/(\sigma\delta)$. Choose $T_-$ after $C_{\rm cert}$ so that this is as small a multiple of $|L|$ as required. The future quadrupole tail is at most $C\int_{T_+}^\infty e^{-r\sigma v}(1+v)^2\idiff v$, which is made small by $T_+$. On the middle part, only the quadrupole comparison is needed: its error contributes at most $\varepsilon_Q\int_{-T_-}^{T_+}e^{-r\sigma v}\idiff v$. The impulse and cubic terms are already $o(1)$ by \eqref{cas:age-integral}; no analytic contour estimate after exit enters this argument. Changing variables from $s$ to $v$ in \eqref{cas:analytic-source} and dividing by $B_{0,n}$ gives the clock factor $e^{-\sigma r(\widetilde\kappa_{n+1}^N-\kappa_{n+1})}\leq e^{\sigma B}$. Using \eqref{cas:absolute-transform}, choosing the tails and the middle error small, and finally increasing $\Lambda$, we obtain
\begin{equation}\label{cas:certified-prefix}
 J_{{\rm adj},n}(H)\leq K_{\rm pref}B_{0,n},\qquad
 K_{\rm pref}\leq e^{\sigma B}(K_{\rm iso}+1)
\end{equation}
for every prefix. A prefix ending inside the window only shortens a nonnegative integral. Thus $K_{\rm pref}$ does not depend on a later contour-cone bound or on a new output constant from the forced-cell estimate.

\subsection{Certification, exit times, and chronological induction}

The deterministic unit certificate for level $n$ is
\begin{equation}\label{cas:certificate}
 F_{c,n}=\kappa_n-W>0,
\end{equation}
with physical time $\tau_nF_{c,n}$. For large $\Lambda$, the child's entire comparison window lies before the parent's certificate, because the clocks satisfy
\begin{equation}\label{cas:certificate-order}
 r(\kappa_{n+1}+B+T_++2)<\kappa_n-W
 \quad\Longleftrightarrow\quad
 \ell_n>W+r(B+T_++2).
\end{equation}
The last inequality follows uniformly from \eqref{cas:ell-size}. Hence certification proceeds from the smallest spatial scale towards the largest.

Before the certificate, the partner and all finer sources have the unconditional unit strain bound $Ch_n^2$. Their contribution to the relative-response error of Proposition~\ref{cell:response} is at most
\begin{equation}\label{cas:wait-fine}
 C(1+\kappa_n)\kappa_nh_n^2=o(1).
\end{equation}
For the coarser sources, \eqref{cas:head-final} gives instead
\begin{equation}\label{cas:wait-head}
 \int_0^{F_{c,n}}\mathfrak a_{\rm head}(s)\idiff s
 \leq C\mu_nP_na_*e^{-\sigma W}.
\end{equation}
Multiplication by $1+\kappa_n$ still gives $o(1)$ by \eqref{cas:head-small-factor}. This use of the causal estimate is necessary: the raw coarser bound over the whole waiting interval would produce $\kappa_n^2A_n^{-1}B_{\rm mass}$, which has not been bounded.

After the certificate, use the fixed deterministic horizon $H_n=\kappa_n+B+2$. The length of this remaining interval is $W+B+2$. No coarser cone is needed there. The raw bound \eqref{cas:raw-exterior} gives
\begin{equation}\label{cas:late-input}
 \begin{gathered}
 B_{0,n}^{-1}\int_{F_{c,n}}^{H_n}e^{-\sigma s}\mathfrak a_n(s)\idiff s
 \leq C\mathfrak a_{\rm raw,n}e^{\sigma W}/a_*=o(1),\\
 (1+\kappa_n)\int_{F_{c,n}}^{H_n}\mathfrak a_n(s)\idiff s
 \leq C(1+\kappa_n)(W+B+2)(A_n^{-1}B_{\rm mass}+h_n^2)=o(1).
 \end{gathered}
\end{equation}
The second estimate follows from \eqref{cas:waiting-small}. The raw pointwise strain near the section also tends to zero uniformly and is smaller than any prescribed fixed multiple of $a_*$.

We now perform the finite induction. Initially every cell is in its contour chart, and every uncertified cone has a strict margin. Stop at a first failure of a cone, of the analytic-radius margin, or of a local reference-tube margin before the next deterministic certificate. The estimates on this stopped history remain valid. The normalized absolute weighted input to an uncertified receiver is at most
\begin{equation}\label{cas:bootstrap-budget}
 K_{\rm pref}+1+\varepsilon_{\rm cone}C_{\rm cone}
                    +o(1)(1+C_{\rm cone}),
\end{equation}
where the term $1$ allows the terminal seed and $\varepsilon_{\rm cone}$ is made small by \eqref{cas:uncertified-input}. Choose $\varepsilon_{\rm cone}C_{\rm cone}\leq1$ and then make the last term at most one. The reserve in \eqref{cas:input-budget} yields the more precise bound
\begin{equation}\label{cas:reserved-budget}
 K_{\rm pref}+1+\varepsilon_{\rm cone}C_{\rm cone}
                    +o(1)(1+C_{\rm cone})
 \leq \frac{K_{\rm in}}{4e^{\sigma B}}.
\end{equation}
The same bound, after a further harmless increase of $\Lambda$, includes the late input in \eqref{cas:late-input}. Apply Proposition~\ref{cell:response} with amplitude $B_{0,n}$. It improves the contour cone by a factor at least two. Its contour contribution to the analytic-radius loss is bounded by $CK_{\rm in}a_*e^{-\sigma W}$, and its exterior contribution is small by \eqref{cas:wait-fine}--\eqref{cas:wait-head}. Its reference-center estimate has no waiting-time factor. Thus all three possible first failures are excluded.

At the next certificate define the actual growing coefficient from the past by
\begin{equation}\label{cas:input-coefficient}
 B_n^{\rm in}=z_{{\rm seed},n}
       +\int_0^{F_{c,n}}e^{-\lambda s}f_{{\rm gr},n}(s)\idiff s,
\end{equation}
where $z_{{\rm seed},n}=0$ for $n<N$ and $z_{{\rm seed},N}=a_*e^{-\lambda\kappa_N}$. At the terminal level, the partner and coarser estimates give
\begin{equation}\label{cas:terminal-input}
 B_N^{\rm in}=a_*e^{-\lambda\kappa_N}+o(B_{0,N}).
\end{equation}
For a nonterminal receiver, \eqref{cas:certificate-order} places the child's comparison window in the past. Using \eqref{cas:adjacent-forcing}, splitting the donor transform at the fixed windows, and applying \eqref{cas:head-final}, \eqref{cas:deep-small}, and \eqref{cas:partner-input}, we obtain
\begin{equation}\label{cas:matched-input}
 B_n^{\rm in}
 =D_n\zeta_{n+1}^Ne^{-\lambda r\widetilde\kappa_{n+1}^N}
          (1+\varepsilon_n),\qquad |\varepsilon_n|\leq\varepsilon_D.
\end{equation}
Here $\varepsilon_D>0$ is prescribed sufficiently small before the final choice of $\Lambda$. The transform tails are made smaller than a fixed fraction of $\varepsilon_D|L|$ by the windows. On the comparison window, the moment error contributes at most
$\varepsilon_Q\int_{-T_-}^{T_+}e^{-r\sigma v}\idiff v$ times the transfer prefactor; use a smaller moment tolerance if necessary. The remaining normalized errors are bounded by
\begin{equation}\label{cas:matching-errors}
 C\bigl(\Lambda^{1-M}+h_{n+1}^2+h_{n+1}q_n\bigr)
       +o(1),
\end{equation}
with the uniform moment envelopes already established. The lower bound for $|L|$ in \eqref{cas:absolute-transform} justifies every normalization.

Assume inductively that the child event belongs to its clock box $|\widetilde\kappa_{n+1}^N-\kappa_{n+1}|<B$. Then \eqref{cas:matched-input} implies
\begin{equation}\label{cas:coefficient-comparison}
 e^{-\sigma rB}(1-\varepsilon_D)
 \leq\frac{|B_n^{\rm in}|}{B_{0,n}}
 \leq e^{\sigma rB}(1+\varepsilon_D).
\end{equation}
For $\varepsilon_D\leq1/2$, the lower bound in \eqref{cas:coefficient-comparison} is at least $\tfrac12e^{-\sigma B}$. Therefore \eqref{cas:reserved-budget}, divided by $|B_n^{\rm in}|/B_{0,n}$, is at most $K_{\rm in}/2$. This verifies the same absolute-input budget in Proposition~\ref{cell:response} after changing its normalization from $B_{0,n}$ to $|B_n^{\rm in}|$. The predicted crossing time satisfies
\begin{equation}\label{cas:predicted-event}
 T_n^{\rm pred}:=\frac1\sigma\log\frac{a_*}{|B_n^{\rm in}|}
 =\ell_n+r\widetilde\kappa_{n+1}^N+O(\varepsilon_D).
\end{equation}
Choose $\varepsilon_D$ small compared with $\delta B$. Then $T_n^{\rm pred}$ lies strictly inside $|s-\kappa_n|<B$, and the input amplitude satisfies
\begin{equation}\label{cas:local-admissibility}
 |B_n^{\rm in}|e^{\sigma H_n}\leq2a_*e^{2\sigma(B+1)}.
\end{equation}
This is admissible by \eqref{cas:section-choice}. The certificate precedes $T_n^{\rm pred}$ by the fixed interval required in Proposition~\ref{cell:first-exit}. Equations \eqref{cas:wait-fine}--\eqref{cas:late-input} verify its relative-response and late-tail assumptions; the pointwise raw bound verifies transversality. That proposition gives the first event, with its clock error arbitrarily small after increasing $\Lambda$. At the terminal level use \eqref{cas:terminal-input} instead. Taking the total error smaller than $\varepsilon_t$ proves the recursion in \eqref{cas:event-recursion}.

The clock boxes close by backward contraction. Iterating the recursion gives
\begin{equation}\label{cas:clock-contraction}
 |\widetilde\kappa_n^N-\kappa_n|
 \leq\varepsilon_t\sum_{j\geq0}r^j
 =\varepsilon_t/\delta<B/10.
\end{equation}
The same local exit result gives the entire envelope \eqref{cas:ancient-envelope}, not only a bound near the section. Its negative-window comparison is valid because $W=B+2T_-$. More precisely, the comparison error in \eqref{cell:donor-comparison-error} is bounded by the unweighted exterior integral plus $a_*$ times the relative-response and late-tail errors in \eqref{cas:wait-fine}--\eqref{cas:late-input}. It tends to zero with $\Lambda$ even after the nonzero matching tolerance $\varepsilon_D$ has been fixed, because the comparison uses the actual exit time and phase. Thus the required moment tolerance does not feed recursively into itself. For the positive window, restart the exterior reference at the actual outer centroid at exit. On the fixed ball containing the cell through unit time $T_+$, the exterior Lipschitz norm tends to zero by \eqref{cas:raw-exterior}. The comparison \eqref{cell:positive-comparison} then gives convergence of the material domain indicators to the isolated orbit, uniformly in the exit phase. The supports remain in a fixed unit ball by \eqref{cas:post-support}; their centroids and centered quadratic moments therefore converge as well. This proves \eqref{cas:window}. It does not require a future coarser cone. The estimates \eqref{cas:post-support} and \eqref{cas:post-impulse} remain valid at all larger positive ages, so the newly certified cell supplies every later prefix estimate. This completes the chronological induction for all finite $N$.

\subsection{A fixed observable at the outer scale}

At the outer first event, the growing coordinate has modulus $a_*$ and the complementary local coordinate tends to zero as $\Lambda\to\infty$. The donor expansion and the small fixed choice of $a_*$ imply that the centered unit quadratic moment has modulus at least $ca_*$. Its unit impulse has modulus at most $Ch_0^3a_*$ by \eqref{cas:ancient-impulse}.

There are no coarser cells at level zero. The partner and finer labels give exterior speed at most $Ch_0^2R_0$ near $c_0$, so the reference center moves by at most $CT h_0d_0=O(\Lambda^{1-M})d_0$. The actual outer centroid differs from that reference by a bounded multiple of $d_0$ by Proposition~\ref{cell:response}. For any zero-circulation scalar, translation of the unbarred quadratic moment obeys the exact identity
\begin{equation}\label{cas:translation-identity}
 \int(x-c_0)^2\omega(x)\idiff x
 =\int(x-c^{\rm out})^2\omega(x)\idiff x
   +2(c^{\rm out}-c_0)\int(x-c^{\rm out})\omega(x)\idiff x.
\end{equation}
The last term is smaller than half the main moment for large $\Lambda$. Hence the quadratic moment of the outer positive label about the original center has modulus at least $cA_0d_0^4a_*=cd_0^4a_*$.

Choose a smooth cutoff $\chi$ equal to one on $B(c_0,R_0/8)$ and supported in $B(c_0,R_0/4)$. All other tubes are outside its support. The four real tests
\begin{equation}\label{cas:four-tests}
 \chi\operatorname{Re}\bigl((x-c_0)^2\bigr),\quad
 -\chi\operatorname{Re}\bigl((x-c_0)^2\bigr),\quad
 \chi\operatorname{Im}\bigl((x-c_0)^2\bigr),\quad
 -\chi\operatorname{Im}\bigl((x-c_0)^2\bigr)
\end{equation}
therefore isolate that label. At each cutoff one of their integrals is at least a fixed positive multiple of $d_0^4a_*$. One test has this property along an infinite subsequence. Every test in \eqref{cas:four-tests} has zero integral against the complete quiet datum by radiality. This proves \eqref{cas:witness} and the theorem.
\end{proof}

\section{Selecting the scales from the vorticity envelope}
\label{sec:envelopes}

The finite-array theorem requires two properties of the logarithmic
scales: successive scales have controlled ratios, and their
amplitude-weighted sum is finite.  We now obtain both properties
from \(\mathcal I_\Theta<\infty\).  This step also identifies the
Lebesgue growth of the complete stationary datum.  Its proof uses
only monotonicity of the prescribed envelope.

\subsection{A regular minorant}

\begin{lemma}\label{env:minorant}
Suppose that \(\Theta:[p_0,\infty)\to[1,\infty)\) is finite and
nondecreasing and that \(\mathcal I_\Theta<\infty\).
For every \(c_{\mathrm g}\in(0,1)\) there are \(p_1\geq p_0\),
a positive \(C^\infty\) function \(g\) on \([p_1,\infty)\), and
a finite constant \(D_{\mathrm g}\) such that
\begin{equation}\label{env:regular}
 g(p)\leq\Theta(p),\qquad
 c_{\mathrm g}\leq\frac{pg'(p)}{g(p)}\leq D_{\mathrm g},
 \qquad
 \int_{p_1}^\infty\frac{\idiff p}{g(p)}<\infty.
\end{equation}
\end{lemma}

\begin{proof}
Choose \(K_{\mathrm g}>2\), set
\(a_j=\Theta(p_0 2^j)\), and define
\(b_0=a_0\) and \(b_{j+1}=\min(a_{j+1},K_{\mathrm g}b_j)\).
Then \(b_j\) is nondecreasing,
\(b_j\leq a_j\), and \(b_{j+1}\leq K_{\mathrm g}b_j\).
The bound
\begin{equation}\label{env:dyadic-summability}
 \int_{p_0 2^{j-1}}^{p_0 2^j}
       \frac{\idiff p}{\Theta(p)}
 \geq \frac{p_0 2^{j-1}}{a_j}\qquad(j\geq1)
\end{equation}
shows that \(\sum_j2^j/a_j<\infty\).
With \(e_j=2^j/b_j\), the defining recurrence gives
\begin{equation}\label{env:recurrence}
 e_{j+1}
 =\max\left(\frac{2^{j+1}}{a_{j+1}},
                      \frac{2}{K_{\mathrm g}}e_j\right)
 \leq \frac{2^{j+1}}{a_{j+1}}
                      +\frac{2}{K_{\mathrm g}}e_j.
\end{equation}
Summing \eqref{env:recurrence} and using \(2/K_{\mathrm g}<1\)
proves \(\sum_je_j<\infty\).
Consequently the step function
\(B(p)=b_j\) for \(p_0 2^j\leq p<p_0 2^{j+1}\) satisfies
\begin{equation}\label{env:step}
 B\leq\Theta,\qquad B(2p)\leq K_{\mathrm g}B(p),\qquad
 \int_{p_0}^\infty B(p)^{-1}\idiff p<\infty.
\end{equation}

To impose a positive lower logarithmic slope, define
\begin{equation}\label{env:lower-slope}
 a_{\mathrm g}(p)=\inf_{q\geq p}
                      \left(\frac pq\right)^{c_{\mathrm g}}B(q).
\end{equation}
This infimum is positive.  Indeed, monotonicity and
\eqref{env:step} imply
\begin{equation}\label{env:positive-infimum}
 \frac{q}{2B(q)}
 \leq \int_{q/2}^q\frac{\idiff s}{B(s)}
 \longrightarrow0,
\end{equation}
so \(B(q)/q^{c_{\mathrm g}}\to\infty\).
Changing variables in the infimum gives, for \(\vartheta\geq1\),
\begin{equation}\label{env:slope-doubling}
 a_{\mathrm g}(\vartheta p)\geq
                 \vartheta^{c_{\mathrm g}}a_{\mathrm g}(p),
 \qquad
 a_{\mathrm g}(2p)\leq K_{\mathrm g}a_{\mathrm g}(p),
 \qquad a_{\mathrm g}(p)\leq B(p).
\end{equation}
Splitting the variable \(q\) in \eqref{env:lower-slope}
into dyadic intervals also yields
\begin{equation}\label{env:reciprocal-majorant}
 \frac1{a_{\mathrm g}(p)}
 \leq 2^{c_{\mathrm g}}\sum_{j\geq0}
                   \frac{2^{c_{\mathrm g}j}}{B(2^jp)}.
\end{equation}
Tonelli's theorem and \(s=2^jp\) therefore give
\begin{equation}\label{env:reciprocal-integral}
 \int_{p_0}^\infty\frac{\idiff p}{a_{\mathrm g}(p)}
 \leq \frac{2^{c_{\mathrm g}}}{1-2^{c_{\mathrm g}-1}}
                      \int_{p_0}^\infty\frac{\idiff s}{B(s)}
 <\infty.
\end{equation}

Choose a nonnegative \(\varphi_{\mathrm g}\in C_c^\infty(0,\log2)\)
with integral one, and set
\begin{equation}\label{env:smoothing}
 g(p)=\int_0^{\log2}a_{\mathrm g}(pe^{-v})
                     \varphi_{\mathrm g}(v)\idiff v,
 \qquad p\geq2p_0.
\end{equation}
The inequalities \eqref{env:slope-doubling} imply
\(a_{\mathrm g}(p)/K_{\mathrm g}\leq g(p)\leq a_{\mathrm g}(p)\)
and \(g(\vartheta p)\geq\vartheta^{c_{\mathrm g}}g(p)\).
Thus \(pg'(p)\geq c_{\mathrm g}g(p)\).
Differentiation in logarithmic coordinates gives
\begin{equation}\label{env:smoothing-derivative}
 |pg'(p)|
 \leq \norm{\varphi_{\mathrm g}'}_1a_{\mathrm g}(p)
 \leq K_{\mathrm g}\norm{\varphi_{\mathrm g}'}_1g(p).
\end{equation}
Equations \eqref{env:reciprocal-integral} and
\eqref{env:smoothing-derivative} prove \eqref{env:regular}
with \(p_1=2p_0\).
\end{proof}

\subsection{Logarithmic scales and the complete datum}

Fix the local parameters \(a_*\) and \(r\) in the order specified
in Theorem~\ref{cas:finite}.  Apply Lemma~\ref{env:minorant},
choose \(s_0\geq p_1\), and define \(s_n\) by
\begin{equation}\label{env:sequence}
 g(s_n)=g(s_0)r^{-n}.
\end{equation}
The lower logarithmic slope in \eqref{env:regular} makes this
definition unique.  Integrating the two derivative bounds gives
\begin{equation}\label{env:ratios}
 1<\nu:=r^{-1/D_{\mathrm g}}
 \leq\frac{s_{n+1}}{s_n}
 \leq r^{-1/c_{\mathrm g}}=:V<\infty.
\end{equation}
Moreover, monotonicity gives
\begin{equation}\label{env:clock-sum}
 \int_{s_n}^{s_{n+1}}\frac{\idiff s}{g(s)}
 \geq\frac{r(s_{n+1}-s_n)}{g(s_n)}
 \geq\frac{r(\nu-1)}{g(s_0)}r^ns_n.
\end{equation}
The last condition in \eqref{env:regular} therefore implies
\(\sum_{n\geq0}r^ns_n<\infty\), as required by
Theorem~\ref{cas:finite}.

Choose \(M>D_{\mathrm g}+1\) and then choose the final parameter
\(\Lambda\) sufficiently large for that theorem.  In its notation,
the scales and the complete datum are
\begin{align}
 S_n&=\Lambda s_n,\quad
 R_n=e^{-(S_n-S_0)},\quad h_n=S_n^{-M},\quad
 d_n=h_nR_n,\quad A_n=r^{-n},\quad c_n=R_ne_1,
                                                   \label{env:scales}\\
 \omega_0(x)
 &=\sum_{n\geq0}(-1)^nA_n
       \left(\Omega\left(\frac{x-c_n}{d_n}\right)
             +\Omega\left(\frac{x+c_n}{d_n}\right)\right).
                                                   \label{env:datum}
\end{align}
The reflection used in the finite evolutions does not appear in
\eqref{env:datum}, because \(\Omega\) is radial.
The disks in \eqref{env:datum} are pairwise disjoint, lie in a fixed
ball, and accumulate only at the origin.

\begin{proposition}[Vorticity growth]\label{env:norms}
The datum \eqref{env:datum} belongs to
\(\mathscr S\cap Y^\Theta\).  There are constants
\(c,C>0\), independent of \(p\), such that
\begin{equation}\label{env:two-sided}
 c\,g(p)\leq\norm{\omega_0}_{L^p}\leq Cg(p)
 \quad(p\geq S_0),\qquad
 \norm{\omega_0}_{L^p}\leq C\Theta(p)\quad(p\geq p_0).
\end{equation}
Its canonical velocity is stationary, and no distributional source
occurs at the accumulation point.
\end{proposition}

\begin{proof}
Disjointness gives the exact formula
\begin{equation}\label{env:norm-formula}
 \norm{\omega_0}_{L^p}^p
 =2\norm{\Omega}_{L^p}^p\sum_{n\geq0}A_n^pd_n^2,
 \qquad
 \norm{\Omega}_{L^p}^p
 =\frac\pi3a^p+\frac{2\pi}3\left(\frac a2\right)^p.
\end{equation}
Since \(A_n=g(s_n)/g(s_0)\) and \(s_n=S_n/\Lambda\),
\eqref{env:regular} implies
\begin{equation}\label{env:amplitude-bound}
 A_n\leq Cg(p)\max\left(1,\left(\frac{S_n}{p}\right)^{D_{\mathrm g}}\right)
 \qquad(p\geq p_1).
\end{equation}
For \(S>0\), maximization of \(S^{D_{\mathrm g}p}e^{-S}\)
at \(S=D_{\mathrm g}p\) yields
\begin{equation}\label{env:maximum}
 e^{-S}\max\left(1,\left(\frac Sp\right)^{D_{\mathrm g}p}\right)
 \leq\max\left(1,\left(\frac{D_{\mathrm g}}e\right)^{D_{\mathrm g}p}\right).
\end{equation}
Combining \eqref{env:amplitude-bound} and \eqref{env:maximum},
we obtain
\begin{equation}\label{env:norm-upper}
 \sum_{n\geq0}A_n^pd_n^2
 \leq C^pg(p)^p e^{2S_0}
       \max\left(1,\left(\frac{D_{\mathrm g}}e\right)^{D_{\mathrm g}p}\right)
       \sum_{n\geq0}S_n^{-2M}e^{-S_n}.
\end{equation}
The last series is finite by \eqref{env:ratios}.
Its \(p\)-th root and \(e^{2S_0/p}\) are uniformly bounded
for \(p\geq p_1\).  Equations \eqref{env:norm-formula} and
\eqref{env:norm-upper} prove the upper bound.
They also prove membership in every finite \(L^p\);
the remaining compact interval of exponents is covered by
interpolation and \(\Theta\geq1\).

For the lower bound, fix \(p\geq S_0\) and choose \(n\) with
\(S_n\leq p<S_{n+1}\leq VS_n\).
The upper logarithmic slope gives
\begin{equation}\label{env:amplitude-lower}
 A_n=\frac{g(S_n/\Lambda)}{g(s_0)}
 \geq\frac{g(p)}{g(s_0)(\Lambda V)^{D_{\mathrm g}}}.
\end{equation}
Also, the physical radius satisfies
\begin{equation}\label{env:area-lower}
 d_n^{2/p}
 =\exp\left(-\frac{2M\log S_n}{p}
            -\frac{2(S_n-S_0)}p\right)
 \geq e^{-2M/e-2}.
\end{equation}
The contribution of this one cell in \eqref{env:norm-formula}
and \eqref{env:amplitude-lower}--\eqref{env:area-lower}
prove the lower bound in \eqref{env:two-sided}.
Stationarity, including at the accumulation point, follows from
Lemma~\ref{uniq:stationary-lemma}; all its hypotheses have now
been verified.
\end{proof}

\subsection{The finite initial data}

Let \(\omega_0^N\) be the actual finite initial datum in
Theorem~\ref{cas:finite}.  At levels \(n<N\) its cells are exactly
quiet, while at level \(N\) they are the exact manifold seed.
No levels beyond \(N\) are present.
The two interfaces of every manifold cell enclose exactly the
same areas as those of \(\Omega\).  Consequently, we obtain
\begin{equation}\label{env:finite-norms}
 \norm{\omega_0^N}_{L^p}^p
       =2\norm{\Omega}_{L^p}^p\sum_{n=0}^NA_n^pd_n^2
       \leq\norm{\omega_0}_{L^p}^p.
\end{equation}
The symmetric-difference estimate in Proposition~\ref{cell:donor}
also gives
\begin{equation}\label{env:initial-convergence}
 \norm{\omega_0^N-\omega_0}_{L^p}^p
 \leq C_pA_N^pd_N^2a_*e^{-\sigma\kappa_N}
       +2\norm{\Omega}_{L^p}^p
                    \sum_{n>N}A_n^pd_n^2
 \longrightarrow0
 \qquad(1\leq p<\infty).
\end{equation}
The constants in this limit may depend on \(p\), which is fixed
before \(N\to\infty\).
Every finite Euler evolution conserves
\eqref{env:finite-norms}.  These are the uniform estimates used
in the passage to the limit.

\section{The limiting Lagrangian evolution}\label{sec:limit}

We now pass from the actual finite Euler solutions to one solution
with the complete datum \eqref{env:datum}.  There are two separate
points to prove.  First, compactness must preserve the nonlinear
equation and its strong initial trace.  Second, the outer-cell
moment from Theorem~\ref{cas:finite} must remain nonzero.
The transport theorem is applied only after the limiting velocity
has been identified.

\subsection{Velocity compactness}

Write \((\omega^N,u^N)\) for the finite solutions on their common
interval \([0,T]\).  By Theorem~\ref{cas:finite}, all vorticities
are supported in one ball \(B_{R_*}\).
Conservation of \eqref{env:finite-norms} implies
\begin{equation}\label{lim:uniform-norms}
 \sup_N\sup_{0\leq t\leq T}\norm{\omega^N(t)}_{L^p}
       \leq\norm{\omega_0}_{L^p},
 \qquad 1\leq p<\infty.
\end{equation}
Splitting the Biot--Savart integral into \(|x-y|<1\) and its
complement, and using \(p=4\) in the first part, gives
\begin{equation}\label{lim:bounded-speed}
 \norm{u^N(t)}_\infty
 \leq C\left(\norm{\omega^N(t)}_4+\norm{\omega^N(t)}_1\right)
 \leq C.
\end{equation}
Zero circulation yields the uniform far-field bound
\begin{equation}\label{lim:far-field}
 |u^N(t,x)|
 =\left|\int_{B_{R_*}}
       \left(K(x-y)-K(x)\right)\omega^N(t,y)\idiff y\right|
 \leq CR_*\norm{\omega_0}_1|x|^{-2}
 \quad(|x|\geq2R_*).
\end{equation}
In particular, the velocities have uniformly small \(L^2\)
tails at spatial infinity.

The fixed-exponent Fourier identity for divergence-free fields
gives \(\norm{\nabla u^N}_2=\norm{\omega^N}_2\).
Differentiating the canonical convolution distributionally gives
\begin{equation}\label{lim:time-derivative}
 \partial_tu^N
   =-\sum_{j=1}^2\partial_jK*
                      \left(u_j^N\omega^N\right).
\end{equation}
Each \(\partial_jK\) is an order-zero Fourier multiplier
bounded on \(L^2\).  Thus \eqref{lim:uniform-norms} and
\eqref{lim:bounded-speed} imply
\begin{equation}\label{lim:H1-time}
 \sup_N\left(
       \norm{u^N}_{L^\infty(0,T;H^1)}
       +\norm{\partial_tu^N}_{L^\infty(0,T;L^2)}
          \right)<\infty.
\end{equation}
At each rational time, local Rellich compactness and
\eqref{lim:far-field} give compactness in \(L^2(\R^2)\).
A diagonal extraction and the common time modulus from
\eqref{lim:H1-time} therefore give
\begin{equation}\label{lim:velocity-strong}
 u^N\longrightarrow u
       \quad\text{in }C\left([0,T];L^2(\R^2)\right).
\end{equation}

After a further diagonal extraction at integer exponents,
the vorticities converge weakly star in
\(L^\infty(0,T;L^p)\), \(1<p<\infty\), to a scalar \(\omega\)
with the same common support.  Their \(L^1\) bounds follow
from that support and the \(L^2\) bound.
For a smooth compactly supported spacetime test, the nonlinear
flux passes to the limit: the factor \(u^N-u\) converges
strongly in \(L^2\), while \(\omega^N\) is bounded in
\(L^\infty_tL^2_x\); the remaining term uses the weak
convergence of \(\omega^N\).
Together with \eqref{env:initial-convergence}, this proves
\begin{equation}\label{lim:vorticity-equation}
 \partial_t\omega+\Div(u\omega)=0,\qquad
 \omega|_{t=0}=\omega_0
\end{equation}
in the distributional formulation with initial trace.

Taking the curl and divergence in \eqref{lim:velocity-strong}
shows that \(\curl u=\omega\) and \(\Div u=0\).
The zero-circulation, compactly supported scalar \(\omega\)
has \(K*\omega\in L^2\).  Hence \(u-K*\omega\) is a
divergence-free and curl-free \(L^2\) vector field, and is
zero.  This identifies the canonical velocity.

\subsection{Exact transport and strong continuity}

We state explicitly the standard transport input used at this point.
If \(v\) is divergence free,
\(v\in L^1(0,T;W_{\mathrm{loc}}^{1,q})\) for some \(q>1\),
and \(v/(1+|x|)\in L^1_t(L^1_x+L^\infty_x)\),
the Sobolev transport theorem gives its unique regular Lagrangian
flow.  For a scalar in a dual integrability class,
distributional transport agrees with transport by that flow.
When the divergence is zero the flow preserves measure and has
an almost-everywhere inverse.  We use this consequence of
DiPerna--Lions theory \cite{DiPernaLions1989}.

Here \(u\in L^\infty_tW_{\mathrm{loc}}^{1,4}\) by
\eqref{lim:uniform-norms} and the fixed-exponent
Calder\'on--Zygmund estimate.  It is bounded by
\eqref{lim:bounded-speed}, and \(\omega\in L^\infty_tL^2_x\).
Thus the spatial growth and scalar integrability hypotheses
hold.  If \(X(t,s,x)\) denotes the resulting two-time flow,
uniqueness for the linear equation with this fixed velocity
identifies \eqref{lim:vorticity-equation} with
\begin{equation}\label{lim:exact-transport}
 \omega(t,x)=\omega_0\left(X(0,t,x)\right),\qquad
 \omega\left(t,X(t,0,x_0)\right)=\omega_0(x_0)
 \quad\text{almost everywhere}.
\end{equation}
The representative in \eqref{lim:exact-transport} is taken
at every time.  Measure preservation gives
\begin{equation}\label{lim:conservation}
 \norm{\omega(t)}_{L^p}=\norm{\omega_0}_{L^p}
 \quad(1\leq p<\infty),\qquad
 \int_{\R^2}\beta(\omega(t,x))\idiff x
 =\int_{\R^2}\beta(\omega_0(x))\idiff x
\end{equation}
for every Borel function \(\beta\) with \(\beta(0)=0\)
and \(\beta(\omega_0)\in L^1\).
The corresponding distributional renormalized transport
identities follow from the same theorem.

For completeness, strong continuity does not require continuity
of the transported scalar in space.  Fix \(s\in[0,T]\)
and approximate \(\omega(s)\) in \(L^p\) by
\(f\in C_c^\infty\).  The group property and measure preservation
give
\begin{align}
 \norm{\omega(t)-\omega(s)}_p
 &\leq 2\norm{\omega(s)-f}_p
       +\norm{f\circ X(s,t,\cdot)-f}_p,
                                                    \label{lim:continuity}\\
 |X(s,t,x)-x|&\leq\norm{u}_\infty|t-s|
       \quad\text{for almost every }x.               \label{lim:displacement}
\end{align}
Uniform continuity of \(f\), together with its bounded support
and \eqref{lim:displacement}, makes the last term in
\eqref{lim:continuity} tend to zero.
Letting the approximation error tend to zero proves
\(\omega\in C([0,T];L^p)\) for every finite \(p\).

In fact, the finite vorticities converge strongly in the same
spaces:
\begin{equation}\label{lim:vorticity-strong}
 \omega^N\longrightarrow\omega
       \quad\text{in }C\left([0,T];L^p(\R^2)\right),
 \qquad 1\leq p<\infty.
\end{equation}
To prove this for \(p>1\), consider any sequence
\(t_N\to t\).
The convergence \eqref{lim:velocity-strong} identifies the
distributional limit of \(\omega^N(t_N)=\curl u^N(t_N)\)
as \(\omega(t)\).  The uniform \(L^p\) bound then gives weak
\(L^p\) convergence.  By \eqref{env:initial-convergence} and
\eqref{lim:conservation}, its norms satisfy
\begin{equation}\label{lim:norm-convergence}
 \norm{\omega^N(t_N)}_p
       =\norm{\omega_0^N}_p
       \longrightarrow\norm{\omega_0}_p
       =\norm{\omega(t)}_p.
\end{equation}
Uniform convexity of \(L^p\) gives strong convergence at
these moving times.  The already proved continuity of
\(\omega\) then proves \eqref{lim:vorticity-strong} by
contradiction.  The case \(p=1\) follows from the case \(p=2\)
and the common compact support.

\subsection{Momentum and energy}

Passing \eqref{lim:time-derivative} to the limit gives
\(\partial_tu\in L^\infty_tL^2_x\).
Also \(u\cdot\nabla u\in L^\infty_tL^2_x\), since
\(u\) is bounded and \(\nabla u\in L^\infty_tL^2_x\).
The vorticity equation and the divergence constraint imply
\begin{equation}\label{lim:projected-momentum}
 \partial_tu+\mathbb P(u\cdot\nabla u)=0,
\end{equation}
where \(\mathbb P\) is the orthogonal projection in
\(L^2(\R^2;\R^2)\) onto divergence-free vector fields.
Indeed, the left side of \eqref{lim:projected-momentum}
is in \(L^2\), has zero divergence and zero curl, and
therefore vanishes.
The complementary gradient defines a pressure \(P\),
so \eqref{eq:momentum} holds.

The regularity just established permits \(L^2\) differentiation
of the energy.  Since \(\mathbb Pu=u\), the energy derivative satisfies
\begin{equation}\label{lim:energy}
 \frac{\diff}{\diff t}\frac12\norm{u(t)}_2^2
 =-\int_{\R^2}u\cdot(u\cdot\nabla u)\idiff x=0
 \quad\text{for almost every }t.
\end{equation}
For the last equality, first insert a smooth spatial cutoff
equal to one on \(B_R\), integrate by parts, and let
\(R\to\infty\).
The boundary error is bounded by
\(CR^{-1}\norm{u}_\infty
 \int_{R\leq|x|\leq2R}|u|^2\idiff x\), which tends to zero.
Thus no separate energy-admissibility assumption is needed
for the constructed solution.

\subsection{The fixed observable and the main theorem}

Theorem~\ref{cas:finite} provides four fixed smooth real tests
near \(c_0\), one of which detects the outer-cell quadratic
moment along an infinite subsequence.
Fix that test and call it \(\phi_*\).
Its quiet moment is zero.  The corresponding first-exit times
satisfy
\begin{equation}\label{lim:event-times}
 t_0^N\in
 \left[\kappa_0-\frac{\varepsilon_t}{1-r},
       \kappa_0+\frac{\varepsilon_t}{1-r}\right]
 \Subset(0,T).
\end{equation}
After another subsequence, \(t_0^N\to t_*\in(0,T)\).
The uniform strong convergence in \eqref{lim:vorticity-strong}
passes the separation estimate to the limit:
\begin{equation}\label{lim:gap}
 \int_{\R^2}\phi_*(x)
            \left(\omega(t_*,x)-\omega_0(x)\right)\idiff x
 =\lim_{N\to\infty}
        \int_{\R^2}\phi_*(x)\omega^N(t_0^N,x)\idiff x
 \geq c\,d_0^4a_*>0.
\end{equation}
Alternatively, \eqref{lim:velocity-strong} suffices because
\(\int\phi_*\omega^N=-\int\nabla^\perp\phi_*\cdot u^N\).
The positivity in \eqref{lim:gap} is at a fixed physical
scale and is independent of the cutoff.

\begin{proof}[Proof of Theorem~\ref{thm:main}]
For part~\ref{main:unique},
Lemma~\ref{uniq:stationary-lemma} and its explicit flow
\eqref{uniq:quiet-flow} construct the quiet evolution
from every \(\omega_0\in\mathscr S\).
Theorem~\ref{uniq:admissible-theorem}, together with the energy
identity for the class \eqref{eq:solution-class}, identifies
every other Lagrangian evolution with it.
The regular Lagrangian flows then coincide almost everywhere.

For part~\ref{main:nonunique}, choose the envelope scales in
\eqref{env:sequence}--\eqref{env:scales} and the datum
\eqref{env:datum}.  Proposition~\ref{env:norms} verifies its
membership in \(\mathscr S\cap Y^\Theta\).
Theorem~\ref{cas:finite} and
\eqref{lim:velocity-strong}--\eqref{lim:exact-transport}
produce a Lagrangian evolution with that datum.
Call it the \(\mathrm{trig}\) branch and call the stationary
evolution the \(\mathrm{quiet}\) branch.
Equations \eqref{lim:conservation}, \eqref{lim:energy},
and \eqref{lim:gap} give all the assertions on \([0,T]\).
Finally the Euler rescaling
\begin{equation}\label{lim:rescaling}
 \widetilde\omega^j(t,x)=T\omega^j(Tt,x),\qquad
 \widetilde u^j(t,x)=Tu^j(Tt,x),\qquad
 \widetilde X^j(t,x)=X^j(Tt,x)
\end{equation}
gives the result on \([0,1]\).
The datum is multiplied by \(T\), which preserves
\(\mathscr S\cap Y^\Theta\), and the separation time becomes
\(t_*/T\in(0,1)\).
\end{proof}

\vspace{0.2in}
\noindent{\bf Acknowledgments.}
This work was supported in part by NSF Grant DMS-2508463,
the Choi Family Gift Fund, and the Mike Yan Gift Fund.

\vspace{0.2in}
\noindent{\bf Disclosure of AI assistance.}
The authors' investigation of nonuniqueness for the
two-dimensional Euler equations began in early 2026 with
the exploration of several possible scenarios and
accompanying numerical experiments. Following OpenAI's
announcement of its Navier--Stokes result in September
2026, the authors used GPT-6 to organize their research
notes and to assist in developing, checking, and refining
the mathematical arguments in
Sections~\ref{sec:uniqueness}--\ref{sec:limit}.
This assistance included the dyadic comparison for the
optimized modulus, the long-time response and chronological
transfer estimates, and the two-sided norm estimate for
the selected envelope. An initial draft of the present
manuscript was completed in mid-September 2026 and has
since undergone further checking and revision.

Throughout this process, the authors directed the
AI-assisted mathematical exploration, critically examined
the resulting arguments, and revised the proofs and
exposition. The proofs are presented analytically in the
manuscript; neither numerical experiments nor AI-generated
output is used as a substitute for mathematical
justification. The authors take full responsibility for
the mathematical content, correctness, and presentation
of the results.

\bibliographystyle{myalpha}
\bibliography{references}
\end{document}

%% file: macro-3.tex
\usepackage[usenames,dvipsnames]{xcolor}
\usepackage{fancyhdr}
\usepackage{amsmath,amsfonts,amsbsy,amsgen,amscd,amssymb,amsthm}
\usepackage{url}

\usepackage{mathtools}

\mathtoolsset{showonlyrefs}

\usepackage[font=small,margin=0.25in,labelfont={bf},labelsep={colon}]{caption}
\usepackage{subcaption}

\usepackage{tikz}
\usepackage{microtype}
\usepackage{enumitem}

\definecolor{dark-gray}{gray}{0.3}
\definecolor{dkgray}{rgb}{.4,.4,.4}
\definecolor{dkblue}{rgb}{0,0,.5}
\definecolor{medblue}{rgb}{0,0,.75}
\definecolor{rust}{rgb}{0.5,0.1,0.1}

\usepackage[colorlinks=true]{hyperref}

\hypersetup{urlcolor=rust}
\hypersetup{citecolor=dkblue}
\hypersetup{linkcolor=dkblue}

\usepackage{setspace}
\usepackage{graphicx}
\usepackage{booktabs,longtable,tabu} % Nice tables
\usepackage{multirow} % More control over tables

\usepackage{float}

\usepackage[full]{textcomp}

\usepackage[scaled=.98,sups,osf]{XCharter}% lining figures in math, osf in text
\usepackage[scaled=1.04,varqu,varl]{inconsolata}% inconsolata typewriter
\usepackage[type1]{cabin}% sans serif
\usepackage[charter,vvarbb,scaled=1.07]{newtxmath}
\usepackage[cal=boondoxo]{mathalfa}
\usepackage[T1]{fontenc}

\usepackage{bm} % boldmath must be called after the package

\graphicspath{{figures/}}

\newtheorem{theorem}{Theorem}[section]
\newtheorem{lemma}[theorem]{Lemma}

\newtheorem{proposition}[theorem]{Proposition}

\newtheorem{corollary}[theorem]{Corollary}

\theoremstyle{definition}

\newtheorem{definition}[theorem]{Definition}

\newtheorem{remark}[theorem]{Remark}

\numberwithin{equation}{section} 
\numberwithin{figure}{section}
\numberwithin{table}{section}

\floatstyle{plaintop}
\newfloat{recipe}{thp}{lor}
\floatname{recipe}{Recipe}
\numberwithin{recipe}{section}

\renewcommand{\phi}{\varphi}

\newcommand{\iunit}{\mathrm{i}}

\newcommand{\Id}{\mathbf{I}}

\providecommand{\mathbbm}{\mathbb} % In case we don't load bbm

\newcommand{\R}{\mathbbm{R}}

\newcommand{\C}{\mathbbm{C}}

\newcommand{\real}{\operatorname{Re}}

\newcommand{\diff}{\mathrm{d}}
\newcommand{\idiff}{\, \diff}

\newcommand{\norm}[1]{\left\Vert {#1} \right\Vert}

\newcommand{\triplenorm}[1]{{\left\vert\kern-0.25ex\left\vert\kern-0.25ex\left\vert #1
    \right\vert\kern-0.25ex\right\vert\kern-0.25ex\right\vert}}

\newcommand{\ind}{\mathbbm{1}} %indicatrix function

\DeclareMathOperator{\curl}{curl}

\def\XXint#1#2#3{{%
\setbox0=\hbox{$#1{#2#3}{\int}$}
\vcenter{\hbox{$#2#3$}}\kern-.5\wd0}}

\renewcommand{\leq}{\leqslant}
\renewcommand{\geq}{\geqslant}
\renewcommand{\subset}{\subseteq}

\newcommand{\Div}{\nabla\cdot}